\documentclass[12pt,fleqn, a4]{article}

\usepackage[french]{babel}
\usepackage{inputenc}
\usepackage[T1]{fontenc}

\usepackage{amssymb}
\usepackage{latexsym}
\usepackage{graphics}

\usepackage{graphics}
\usepackage{color}
\newcommand{\colval}{0.3}
\definecolor{colone}{gray}{\colval}

\newcommand{\dcb}{\begin{array}{lll}}
\newcommand{\dce}{\end{array}}
\newcommand{\ebe}{\begin{enumerate}\setlength{\baselineskip}{13pt}\setlength{\parskip}{5pt}}
\newcommand{\dbe}{\end{enumerate}}

\newcommand{\ibegin}{\begin{itemize}\setlength{\baselineskip}{19pt}\setlength{\parskip}{7pt}}
\newcommand{\iend}{\end{itemize}}
\newcommand{\ok}{\rule{4pt}{6pt}}
\newcommand{\desb}{\begin{description}}
\newcommand{\dese}{\end{description}}

\newtheorem{Thm}{Theorem}[section]
\newtheorem {Cor}[Thm]{Corollary}
\newtheorem {definition}[Thm]{Definition}%[section]
\newtheorem {pro}{Proposition}[Thm]
\newtheorem {Lemma}[Thm]{Lemma}
\newtheorem {rem}[Thm]{Remark}%[section]

\newtheorem {assumption}[Thm]{Assumption}%[section]

\newcommand {\bd}{\begin{definition}}
\newcommand {\ed}{\end{definition}}
\newcommand {\bpro}{\begin{pro}}
\newcommand {\epro}{\end{pro}}
\newcommand {\bl}{\begin{Lemma}}
\newcommand {\el}{\end{Lemma}}
\newcommand {\bcor}{\begin{Cor}}
\newcommand {\ecor}{\end{Cor}}
\newcommand {\brem }{\begin{rem} \rm }
\newcommand {\erem }{\end{rem}}
\newcommand{\bethe}{\begin{Thm}}
\newcommand{\ethe}{\end{Thm}}
\newcommand {\bassumption}{\begin{assumption}}
\newcommand {\eassumption}{\end{assumption}}

\def \ind{1\!\!1}
\def\cro#1{\langle #1\rangle}

\begin{document}

\begin{center}
\textbf{\Large Local solution method for the problem of enlargement of filtration}
\end{center}

\begin{center}
Shiqi Song

{\footnotesize Laboratoire Analyse et Probabilités\\
Université d'Evry Val D'Essonne, France\\
shiqi.song@univ-evry.fr}
\end{center}

\section{Introduction}

This is a handout made for the mini-course on stochastic calculus and its applications to mathematics of finance at the University of Tahar Moulay Saida, organized by Professor Abdeldjebbar Kandouci.

\subsection{The motivation}

The enlargement of filtration theory is a study of the semimartingales when the basic filtration changes. This theory provides particular techniques in the stochastic calculus. Such techniques have found numerous recent applications in the mathematics of finance, and these applications in its turn have been the motivation of recent efforts to put forward the known results of the theory and to develop new ones (see \cite{ACDJ, CJZ, FP, FJS,  JYC, JS3, KP1, KP2, LR, MY, N, protter, Song-splitting, Song-drift, XY, TXY, WYY}). To attend in this movements, we present here a methodology, that we will call the local solution method, introduced in \cite{SongThesis}. 

We will show that the local solution method is an effective and flexible method. The method has been checked in various examples that we find in the classical references. Three of them will be presented below, namely firstly the initial enlargement with the whole running supremum of a brownian motion, secondly the progressive enlargement with Emery's last passage time, and finally the filtration expanded by the future infimums of a linear diffusion. Beside of these, using the local solution method, we will give new proofs of three of the classical formulas, namely Jacod's formula, the progressive enlargement formula and the enlargement formula with honest time.

\subsection{Problem of enlargement of filtration} 

We now describe formally what we do below.

We consider a probability space $(\Omega,\mathcal{B},\mathbb{P})$ equipped with two filtrations $\mathbb{F}=(\mathcal{F}_t)_{t\geq 0}$ and $\mathbb{G}=(\mathcal{G}_t)_{t\geq 0}$ of sub-$\sigma$-algebras of $\mathcal{B}$, such that $\mathcal{F}_t\subset \mathcal{G}_t$ for $t\geq 0$. We assume that $\mathbb{F}$ and $\mathbb{G}$ satisfy the usual condition.

Let $S$ be a $\mathbb{F}$ local martingale. The problem of enlargement of filtration for $S$ is the question whether $S$ is a $\mathbb{G}$ semimartingale. If it is the case, $S$ will be a $\mathbb{G}$ special semimartingale. Then, the drift part of $S$ in $\mathbb{G}$ will be called the formula of enlargement of filtration, while the association of $S$ with its drift is called the drift operator.

If all $\mathbb{F}$ local martingales are $\mathbb{G}$ local martingales, we say that $\mathbb{F}$ is immersed in $\mathbb{G}$ and we say that the hypothesis$(H)$ is satisfied. If all $\mathbb{F}$ local martingales are $\mathbb{G}$ semimartingales, we say that the hypothesis$(H')$ is satisfied. See \cite{BY, jacod2, Jeulin80, BJR}

\subsection{Solving the problem}

A usual standard to solve problems of enlargement of filtration is to classify the problems into two categories, namely the initial enlargement or the progressive enlargement, and to apply one of the following methods, i.e., Jacod's criterion, Yor's kernel method in the first case; the formula for progressive enlargement of filtration, the formula for honest times in the second case (see \cite{jacod2, Jeulin80, JY2, yor, MY}).

The local solution method add into the above table the consideration of two other points : on the one hand, the notion of local solutions besides the classifications of initial or progressive enlargements, and on the other hand, the local search of solutions and their integration to form a global one.

The local solution method is issued from the following observation. When one has tried with problems of enlargement of filtration, one has all the experience that something goes through easily somewhere, but usually not always and not everywhere. One notice also that the formula of enlargement of filtration is such a property that its validity can be defined at each particular point and can be checked within a neighbourhood of the point, and in the same time, its validity has to be checked for every point. There exist distinctly a local aspect and a global aspect of the problem of enlargement of filtration.

\section{Basic vocabulary on semimartingale calculus}

\subsection{Stopping times and spaces of semimartingales}

This paper is based on semimartingale calculus. We give some precision on the terminology used below. For a complete guide of semimartingale calculus we refer to \cite{CS, DM, HWY, jacod, protter}. 

We work with the semimartingale space $\mathcal{H}^1$ as defined in \cite{emery, Jeulin80, meyer2}. We define $\mathcal{H}^1_{loc}$ as in \cite{jacod, HWY}. A sequence of semimartingales $(X_n)_{n\geq 1}$ is said to converge in $\mathcal{H}^1_{loc}$, if there exists a sequence $(T_m)_{m\geq 1}$ of stopping times tending to the infinity, such that, for every fixed $T_m$, the sequence of stopped processes $(X^{T_m})_{n\geq 1}$ converges in $\mathcal{H}^1$. 

When $X$ is a special semimartingale, $X$ can be decomposed into $X=M+A$, where $M$ is a local martingale and $A$ is a predictable process with finite variation. We will call $A$ the drift of $X$.

A process $X=(X_t)_{t\geq 0}$ (in some probability space with a filtration) is said integrable, if each $X_t$ is integrable. It is said càdlàg, if for almost all $\omega$, the function $t\rightarrow X_t(\omega)$ is right continuous and has left limit at every point. It is said of class$(D)$, if the family of $X_T$, where $T$ runs over the family of all finite stopping times, is a uniformly integrable family (cf. \cite[p.106]{protter} and also \cite{DM2, HWY}).

The semimartingale calculus depend on the reference probability measure and on the reference filtration. In this paper, where different probability measures or filtrations may be involved in a same computation, exponent such as $\cdot^{\mathbb{P}\cdot\mathbb{F}}$ will be used to indicate the reference probability or filtration. This indication may however be ignored, if no confusion exists. 

The optional (resp. predictable) projection of a process $X$ is denoted by $^{o}\!X$ or $^{\mathbb{P}\cdot\mathbb{F}-o}\!X$ (resp. ${}^{p}\!X$ or $^{\mathbb{P}\cdot\mathbb{F}-p}\!X$). The optional (resp. predictable) dual projection of a process $A$ with finite variation is denoted by $X^{o}$ or $X^{\mathbb{P}\cdot\mathbb{F}-o}$ (resp. $X^{p}$ or $X^{\mathbb{P}\cdot\mathbb{F}-p}$). 

The word "positive" means "non negative", and "increasing" means "non decreasing".

For a function $f$ and a $\sigma$-algebra $\mathcal{T}$, the expression $f\in\mathcal{T}$ will mean that $f$ is $\mathcal{T}$ measurable. 

Let $X,Y$ be random variables taking values respectively in measurable spaces $(E,\mathcal{E})$ and $(F,\mathcal{F})$. We use the notion of regular conditional distribution of $Y$ given $X$ in the sense of \cite{freedman, RW}. We notice that, if $\pi$ is a regular conditional distribution, if $H(x,y)$ is $\mathcal{E}\otimes\mathcal{F}$ measurable and integrable, we have (by monotone class theorem)$$
\int H(X,y) \pi(X,dy)=\mathbb{E}[H(X,Y)|X^{-1}(\mathcal{E})]
$$
almost surely.

Equalities between random variables are almost sure equalities.

\

\subsection{Elementary integral}

We present here the notion of the integrals of elementary functions (which is applicable to elementary processes in an evident way).

Let $f$ be a real càdlàg function defined on an interval $(a,b]$ with $-\infty\leq a<b\leq \infty$ (where we make the convention that $(a,\infty]=(a,\infty)$). For any real number $c$, denote by $f^c$ the function $t\in(a,b]\rightarrow f(t\wedge c)$. A function $h$ on $(a,b]$ is called left elementary function, if $h$ can be written in such a form as : $h=\sum_{i=0}^nd_i\ind_{(x_i,x_{i+1}]}$, where $a=x_0<x_1<\ldots<x_{n}<x_{n+1}=b$ ($n\in\mathbb{N}$) denotes an interval partition of $[a,b]$ and $\{d_0,d_1,\ldots,d_n\}$ are real numbers. We define $h\ind_{(a,b]}\centerdot f$ the elementary integral of $h$ with respect to $f$ on $(a,b]$ to be the following function: $$
t\in(a,b]\rightarrow (h\ind_{(a,b]}\centerdot f)(t) = \sum_{i=0}^nd_i(f^{x_{i+1}}-f^{x_i})(t).
$$
Then,

\bl\label{elem-integral} 
\ebe
\item
If $h$ has another representation $h=\sum_{i=0}^{m}e_i\ind_{(y_i,y_{i+1}]}$, we have$$
(h\ind_{(a,b]}\centerdot f) = \sum_{i=0}^{m}e_i(f^{y_{i+1}}-f^{y_i})
$$
\item
The elementary integral is bi-linear on the product set of the family all left elementary functions and of the family of all real càdlàg functions on $(a,b]$.
\item
If $g$ is another left elementary function on $(a,b]$, we have$$
(g\ind_{(a,b]}\centerdot (h\ind_{(a,b]}\centerdot f))
=(gh\ind_{(a,b]}\centerdot f)
$$
\item
For any real number $c$, $(h\ind_{(a,b]}\centerdot f)^c=h\ind_{(a,b]}\centerdot f^c=h\ind_{(a,c\wedge b]}\centerdot f$

\item
$\Delta_t(h\ind_{(a,b]}\centerdot f)=h(t)\ind_{a<t\leq b}\Delta_t f$ for $t\in\mathbb{R}$.
\dbe
\el

\textbf{Proof.}
\ebe
\item
Let $a=z_1<\ldots<z_k<z_{k+1}=b$ be a refinement of $a=x_0<x_1<\ldots<x_{n}<x_{n+1}=b$ and of $a=y_0<y_1<\ldots<y_{m}<y_{m+1}=b$. We note that $$
\dcb
\sum_{i=0}^nd_i(f^{x_{i+1}}-f^{x_i})&=&\sum_{i=0}^nh(x_{i+1})(f^{x_{i+1}}-f^{x_i}) \\
\sum_{i=0}^{m}e_i(f^{y_{i+1}}-f^{y_i})&=&\sum_{i=0}^{m}h(y_{i+1})(f^{y_{i+1}}-f^{y_i})
\dce
$$ 
Denote by respectively $G$ and $F$ the above two expressions (as functions on $(a,b]$). Then, for any points $s<t$ in some $(z_i,z_{i+1}]$, we have$$
G(t)-G(s)=F(t)-F(s)=h(z_{i+1})(f(t)-f(s))
$$
Since $F(s)=G(s)=0$ for $s\leq a$, $F$ and $G$ coincides.
\item
The bi-linearity is the consequence of the first property.
\item
By the bi-linearity, we need only to check the third property for $h=\ind_{(x,x']}$ and $g=\ind_{(y,y']}$, where $x<x',y<y'$ are points in $(a,b]$. We have
$$
\dcb
(g\ind_{(a,b]}\centerdot (h\ind_{(a,b]}\centerdot f))
&=&(g\ind_{(a,b]}\centerdot (f^{x'}-f^{x}))\\
&=&(f^{x'}-f^{x})^{y'}-(f^{x'}-f^{x})^{y}\\
&=&(f^{x'}-f^{x})^{y'\wedge x'}-(f^{x'}-f^{x})^{y\wedge x'}\\
&=&(f^{x'}-f^{x})^{(y'\wedge x')\vee x}-(f^{x'}-f^{x})^{(y\wedge x')\vee x}\\
&=&f^{(y'\wedge x')\vee x}-f^{x}-f^{(y\wedge x')\vee x}+f^{x}\\
&=&f^{(y'\wedge x')\vee x}-f^{(y\wedge x')\vee x}\\
\dce
$$
If $y>x'$ or $x>y'$, the above function is identically null just as $gh$ and $gh\ind_{(a,b]}\centerdot f$ do. Otherwise, it is equal to $$
f^{y'\wedge x'}-f^{y\vee x}
=gh\ind_{(a,b]}\centerdot f
$$
\item
This fourth property is clear for $h=\ind_{(s,t]}$ for $a\leq s<t\leq b$. By the bi-linearity, it is valid in general.
\item
It is clear for $h=\ind_{(s,t]}$. By the bi-linearity, it is true in general. 
\dbe
The lemma is proved. \ok

Let $\mathtt{B}_i=(a_i,b_i], i\geq 0$, be a sequence of non empty left intervals. 
Look at the union set $\cup_{i\geq 0}\mathtt{B}_i$. One of the following situations hold for a $x\in\mathbb{R}_+$:
\ebe
\item[-]
$x$ is in the interior of one of $\mathtt{B}_i$ for $i\geq 0$.
\item[-]
$x$ is in the interior of $\cup_{i\geq 0}\mathtt{B}_i$, but it is in the interior of no of the $\mathtt{B}_i$ for $i\geq 0$. In this case, $x$ is the right end of one of $\mathtt{B}_i$. There exists a $\epsilon>0$ such that $(x,x+\epsilon)\subset\cup_{i\geq 0}\mathtt{B}_i$, and for any $i\geq 0$, either $\mathtt{B}_i\subset (0,x]$, or $\mathtt{B}_i\subset (x,\infty)$.
\item[-]
$x$ is in $\cup_{i\geq 0}\mathtt{B}_i$, but it is not in its interior. Then, $x$ is the right end of one of $\mathtt{B}_i$ and there exists a sequence of points in $(\cup_{i\geq 0}\mathtt{B}_i)^c$ decreasing to $x$.
\item[-]
$x$ is in $(\cup_{i\geq 0}\mathtt{B}_i)^c$.
\dbe
Consider the right end points $b_i$. A point $b_i$ will be said of first type if $a_j<b_i<b_j$ for some $j$. It is of second type if it is not of first type, but $b_i=a_j$ for some $j$. It is of third type if it is not of first neither of second type.

\bl\label{third-type}
Let $f$ be a càdlàg function on $\mathbb{R}$. Let $a<b$. Suppose that $(a,b]\subset\cup_{i\geq 0}\mathtt{B}_i$. Suppose that the family of right end points of third type has only a finite number of accumulation points in $(a,b)$. Suppose that the limit of $\ind_{(a,b]}\ind_{\cup_{0\leq i\leq n}\mathtt{B}_i}\centerdot f$ exists when $n\uparrow\infty$ with respect to the uniform norm on compact intervals. Then, the limit is simply equal to $\ind_{(a,b]}\centerdot f$.
\el

\textbf{Proof.}
We denoted by $\ind_{(a,b]}\ind_{\cup_{0\leq i<\infty}\mathtt{B}_i}\centerdot f$ the limit of $\ind_{(a,b]}\ind_{\cup_{0\leq i\leq n}\mathtt{B}_i}\centerdot f$.
\ebe
\item
Suppose firstly that there exists no right end point of third type in $(a,b)$. Let $a<s<t<b$. Then, $[s,t]$ is contained in $(\cup_{i\geq 0}\mathtt{B}_i)^\circ$, where the exponent $^\circ$ denotes the interior of a set. Note that in the case we consider here, the right end points are interior points of some $(\cup_{0\leq i\leq n}\mathtt{B}_i)^\circ, n\geq 0$. So, $[s,t]\subset\cup_{n\geq 0}(\cup_{0\leq i\leq n}\mathtt{B}_i)^\circ$. There exists therefore a $N>0$ such that $[s,t]\subset\cup_{n\geq 0}(\cup_{0\leq i\leq n}\mathtt{B}_i)^\circ$. We have$$
\dcb
&&f_t-f_s=(\ind_{(s,t]}\centerdot f)_t\\
&=&(\ind_{(s,t]}\ind_{\cup_{0\leq i\leq N}\mathtt{B}_i}\centerdot f)_t
=\lim_{n\uparrow\infty}(\ind_{(s,t]}\ind_{\cup_{0\leq i\leq n}\mathtt{B}_i}\centerdot f)_t\\
&=&\lim_{n\uparrow\infty}(\ind_{(s,t]}\ind_{(a,b]}\ind_{\cup_{0\leq i\leq n}\mathtt{B}_i}\centerdot f)_t\\
&=&\lim_{n\uparrow\infty}(\ind_{(s,t]}\centerdot (\ind_{(a,b]}\ind_{\cup_{0\leq i\leq n}\mathtt{B}_i}\centerdot f))_t\\

&=&(\ind_{(s,b]}\centerdot (\ind_{(a,b]}\ind_{\cup_{0\leq i<\infty}\mathtt{B}_i}\centerdot f))_t\\
&=&(\ind_{(s,b]}\ind_{\cup_{0\leq i<\infty}\mathtt{B}_i}\centerdot f)_t-(\ind_{(s,b]}\ind_{\cup_{0\leq i<\infty}\mathtt{B}_i}\centerdot f)_s\\
\dce
$$
As $(\ind_{(s,b]}\ind_{\cup_{0\leq i<\infty}\mathtt{B}_i}\centerdot f)_s\rightarrow 0$ when $s\downarrow a$, we obtain$$
f_t-f_a=(\ind_{(s,b]}\ind_{\cup_{0\leq i<\infty}\mathtt{B}_i}\centerdot f)_t, \forall a<t<b.
$$
Now to obtain the result stated in the lemma, we need only to check that $\Delta_b(\ind_{(s,b]}\ind_{\cup_{0\leq i<\infty}\mathtt{B}_i}\centerdot f)=\Delta_b f$. Notice that there exists a $N>0$ such that $b\in\cup_{0\leq i<n}\mathtt{B}_i$ for all $n\geq N$. We have$$
\dcb
\Delta_b(\ind_{(s,b]}\ind_{\cup_{0\leq i<\infty}\mathtt{B}_i}\centerdot f)
&=&\lim_{n\uparrow\infty}\Delta_b(\ind_{(s,b]}\ind_{\cup_{0\leq i<n}\mathtt{B}_i}\centerdot f)\\
&=&\lim_{n\uparrow\infty}\ind_{b\in\cup_{0\leq i<n}\mathtt{B}_i}\Delta_b f\\
&=&\Delta_b f
\dce
$$
The lemma is proved when no right end point of third type exists.

\item
There exist a finite number of right end points of third type in $(a,b)$. Let $v_1<v_2<\ldots<v_k$ are the right end points of third type in $(a,b)$. Applying the preceding result, $$
\dcb
\lim_{n\uparrow\infty}\ind_{(a,b]}\ind_{\cup_{0\leq i\leq n}\mathtt{B}_i}\centerdot f
&=&\sum_{j=0}^k\lim_{n\uparrow\infty}\ind_{(v_j,v_{j+1}]}\ind_{\cup_{0\leq i\leq n}\mathtt{B}_i}\centerdot f \ (v_0=a,v_{k+1}=b)\\

&=&\sum_{j=0}^k\ind_{(v_j,v_{j+1}]}\centerdot f
=\ind_{(a,b]}\centerdot f
\dce
$$
The lemma is true in this second case.
\item
There exist an infinite number of right end points of third type in $(a,b)$, but $b$ is the only accumulation point of these right end points of third type. We have, for $a<t<b$,
$$
\dcb
(\ind_{(a,b]}\ind_{\cup_{0\leq i\leq \infty}\mathtt{B}_i}\centerdot f)_t
&=&(\lim_{n\uparrow\infty}\ind_{(a,b]}\ind_{\cup_{0\leq i\leq n}\mathtt{B}_i}\centerdot f)_t\\
&=&(\lim_{n\uparrow\infty}\ind_{(a,t]}\ind_{\cup_{0\leq i\leq n}\mathtt{B}_i}\centerdot f)_t\\
&=&(\ind_{(a,t]}\centerdot f)_t
=(\ind_{(a,b]}\centerdot f)_t
\dce
$$
As before, the two functions has the same jumps at $b$. The lemma is true in this third case.

\item
There exist an infinite number of right end points of third type in $(a,b)$, but $a$ is the only accumulation point of these right end points of third type. Let $a<s$.
$$
\dcb
&&(\ind_{(a,b]}\ind_{\cup_{0\leq i\leq \infty}\mathtt{B}_i}\centerdot f)-(\ind_{(a,b]}\ind_{\cup_{0\leq i\leq \infty}\mathtt{B}_i}\centerdot f)^s\\
&=&
\ind_{(s,b]}\centerdot(\ind_{(a,b]}\ind_{\cup_{0\leq i\leq \infty}\mathtt{B}_i}\centerdot f)\\
&=&\lim_{n\uparrow\infty}\ind_{(s,b]}\centerdot(\ind_{(a,b]}\ind_{\cup_{0\leq i\leq n}\mathtt{B}_i}\centerdot f)\\
&=&\lim_{n\uparrow\infty}\ind_{(s,b]}\ind_{(a,b]}\ind_{\cup_{0\leq i\leq n}\mathtt{B}_i}\centerdot f\\
&=&\lim_{n\uparrow\infty}\ind_{(s,b]}\ind_{\cup_{0\leq i\leq n}\mathtt{B}_i}\centerdot f\\

&=&\ind_{(s,b]}\centerdot f

\dce
$$
Since $(\ind_{(a,b]}\ind_{\cup_{0\leq i\leq \infty}\mathtt{B}_i}\centerdot f)^s$ tends to zero when $s\downarrow a$, the lemma is true in this fourth case.

\item
There exist an infinite number of right end points of third type in $(a,b)$, but $a,b$ are the only accumulation points of these right end points of third type. We have
$$
\dcb
\lim_{n\uparrow\infty}\ind_{(a,b]}\ind_{\cup_{0\leq i\leq n}\mathtt{B}_i}\centerdot f
&=&\lim_{n\uparrow\infty}\ind_{(a,a+\frac{b-a}{2}]}\ind_{\cup_{0\leq i\leq n}\mathtt{B}_i}\centerdot f
+\lim_{n\uparrow\infty}\ind_{(a+\frac{b-a}{2},b]}\ind_{\cup_{0\leq i\leq n}\mathtt{B}_i}\centerdot f\\

&=&\ind_{(a,a+\frac{b-a}{2}]}\centerdot f+\ind_{(a+\frac{b-a}{2},b]}\centerdot f
=\ind_{(a,b]}\centerdot f
\dce
$$

\item
There exist an infinite number of right end points of third type in $(a,b)$, but there exist only a finite number of accumulation point of these right end points of third type in $(a,b)$. Let $v_1<v_2<\ldots<v_k$ be the accumulation points in $(a,b)$. Applying the preceding result, $$
\dcb
\lim_{n\uparrow\infty}\ind_{(a,b]}\ind_{\cup_{0\leq i\leq n}\mathtt{B}_i}\centerdot f
&=&\sum_{j=0}^k\lim_{n\uparrow\infty}\ind_{(v_j,v_{j+1}]}\ind_{\cup_{0\leq i\leq n}\mathtt{B}_i}\centerdot f \ (v_0=a,v_{k+1}=b)\\

&=&\sum_{j=0}^k\ind_{(v_j,v_{j+1}]}\centerdot f
=\ind_{(a,b]}\centerdot f
\dce
$$
The lemma is proved. \ok
\dbe

\

\

\

\section{Pieces of semimartingales and their integration to form a global semimartingale}\label{noncovering}

We study the question of integrating a collection of pieces of semimartingales into a global semimartingale. The result is initiated in \cite{SongThesis} and achieved in \cite{Song99}.

\subsection{Assumptions}\label{pieces}

Let $S$ be a real càdlàg $\mathbb{G}$ adapted process. Let $\mathtt{B}$ be a $\mathbb{G}$ random left interval (i.e. $\mathtt{B}=(T,U]$ with $T,U$ being two $\mathbb{G}$ stopping times). We say that $S$ is a $\mathbb{G}$ (special) semimartingale on $\mathtt{B}$, if $\ind_\mathtt{B}\centerdot S$ is a $\mathbb{G}$ (special) semimartingale. We consider the following assumption.

\bassumption\label{partialcovering}
We suppose
\ebe
\item[i.]
$S$ is a special semimartingale in its natural filtration.
\item[ii.]
There exists a sequence of random left intervals $(\mathtt{B}_i)_{i\in\mathbb{N}}$ on each of which $S$ is a special semimartingale. In any bounded open interval $\subset\cup_{i\geq 0}\mathtt{B}_i$, there exist only a finite number of accumulation points of the right end points of third type in the sense of Lemma \ref{third-type}. 
\item[iii.]
There is a special semimartingale $\check{S}$ such that $\{S\neq \check{S}\}\subset \cup_{i\geq 0}\mathtt{B}_i$.
\dbe
\eassumption

Remark that, under the above assumption, denoting by $S^c$ the continuous martingale part of $S$ in its natural filtration, the bracket process $$
[S,S]_t=\cro{S^c,S^c}_t+\sum_{s\leq t}(\Delta_s S)^2, \ t\geq 0,
$$
is well-defined locally integrable increasing process in the natural filtration of $S$.

For any $\mathbb{G}$ stopping time $R$, let $d_R=\inf\{s\geq R: s\notin \cup_{i\geq 0}\mathtt{B}_i\}$. Either $d_{R}$ belongs or not to $\cup_{i\geq 0}\mathtt{B}_i$. But always $S_{d_R}=\check{S}_{d_R}$. Let$$
\mathtt{A}=\cup_{s\in\mathbb{Q}_+}(s, d_s],\ \mathtt{C}=\mathtt{A}\setminus \cup_{i\geq 0}\mathtt{B}_i
$$
$\mathtt{A}$ and $\mathtt{C}$ are $\mathbb{G}$ predictable set. $\mathtt{C}$ is in addition a thin set contained in $\cup_{s\in\mathbb{Q}_+}[d_s]$. Let $C$ be the thin process $$
C_t=\ind_{\{t\in\mathtt{C}\}}\Delta_t(S-\check{S}), t\geq 0
$$

We introduce the process
$g_{t} = \sup\{0\leq s < t; s \notin \cup_{n\geq 0}\mathtt{B}_n\},\ t > 0$ (when the set is empty, $g_{t}=0$). It is an increasing left continuous $\mathbb{G}$ adapted process, i.e. a predictable process. For any
$\epsilon > 0$, set
$\mathtt{A}_{\epsilon} = \{s\geq 0: s - g_{s} > \epsilon\}$. We can check that $\mathtt{A}_{\epsilon} \subseteq \mathtt{A}_{\delta}$ for $\delta<\epsilon$ and $\mathtt{A}=\cup_{\epsilon>0}\mathtt{A}_\epsilon$. Note that $0\notin \cup_{i\geq 0}\mathtt{B}_i$. Define successively for every $n\geq 0$ the $\mathbb{G}$ stopping times: ($d_{R_{-1}}=0$)$$
\dcb
R_{n}=\inf\{s\geq d_{R_{n-1}}: s\in \mathtt{A}_\epsilon\}\ n\geq 0.
\dce
$$
For any $s\in\mathtt{A}_\epsilon$, there exists a $s'<s$ such that $[s',s]\subset\mathtt{A}_\epsilon$. Therefore, $R_n\notin \mathtt{A}_\epsilon$ and $R_n<d_{R_n}$ if $R_n<\infty$. Moreover, $R_n\in \cup_{i\geq 0}\mathtt{B}_i$ and $d_{R_n}-g_{d_{R_n}}> \epsilon$ and $(d_{R_n},d_{R_n}+\epsilon)\cap\mathtt{A}_\epsilon=\emptyset$. Consequently $d_{R_{n}}+\epsilon\leq R_{n+1}$ and $\lim_{k\rightarrow\infty}R_k=\infty$. We can write $\mathtt{A}_\epsilon=\cup_{n\geq 1} (R_{n},d_{R_{n}}]$ and hence $\ind_{\mathtt{A}_\epsilon}$ is a left elementary process on any finite interval. 

For a càdlàg process $X$, set $\mathtt{j}(X)=\{t>0:X_{t-}>0>X_{t} \mbox{ or }X_{t-}\leq 0<X_{t}\}$. We introduce the process (eventually taking infinite values)	
$$
\dcb
A_t&=&\sum_{0<s\leq t}\ind_{\{s\in\mathtt{A}\}}\left[\ind_{(S-\check{S})_{s-}>0}(S-\check{S})^-_{s}+\ind_{(S-\check{S})_{s-}\leq 0}(S-\check{S})^+_{s}\right]\\
&=&\sum_{0<s\leq t}\ind_{\{s\in\mathtt{A}\}} \ind_{\mathtt{j}(S-\check{S})}(s)\ |S-\check{S}|_{s}\ t\geq 0.
\dce
$$
Under Assumption \ref{partialcovering}, $\ind_{\mathtt{B}_i}\centerdot S$ is a special semimartingale for any $i\geq 0$. We denote by $\chi^{\mathtt{B}_i}$ the drift of $\ind_{\mathtt{B}_i}\centerdot S$. It is clear that the two random measures $\mathsf{d}\chi^{\mathtt{B}_i}$ and $\mathsf{d}\chi^{\mathtt{B}_j}$ coincides on $\mathtt{B}_i\cap\mathtt{B}_j$ for $i,j\geq 0$. We can therefore define with no ambiguity a $\sigma$-finite (signed) random measure $\mathsf{d}\chi^\cup$ on $\cup_{i\geq 0}\mathtt{B}_i$ such that $\ind_{\mathtt{B}_i}\mathsf{d}\chi^\cup=\mathsf{d}\chi^{\mathtt{B}_i}$.

We will say that a signed random measure $\lambda$ on $\mathbb{R}_+$ has a distribution function, if $\int_0^t|d\lambda|_s<\infty$ for any $t\geq 0$. In this case, the process $t\geq 0 \rightarrow \lambda([0,t])$ is called the distribution function of $\lambda$.

\

\subsection{The results}

\bethe\label{semimartingaleproperty}
Suppose Assumption \ref{partialcovering}. For $S$ to be a semimartingale on the whole $\mathbb{R}_+$, it is necessary and it is sufficient that the random measures $\mathsf{d}\chi^\cup$ has a distribution function $\chi^\cup$ and the process $C$ is the jump process of a special semimartingale. 
\ethe

The necessity of this theorem is clear (recalling that $S$ is locally in class$(D)$). The sufficiency will be the consequence of the following lemmas.

\bl\label{UB}
(We assume only the two first conditions in Assumption \ref{partialcovering}.) Suppose that the random measure $\mathsf{d}\chi^\cup$ has a distribution function $\chi^\cup$. Then, $\ind_{\cup_{0\leq i\leq n}\mathtt{B}_i}\centerdot S, n\geq 0,$ converges in $\mathcal{H}^1_{\mathsf{loc}}$ to a semimartingale that we denote by $\ind_{\cup_{i\geq 0}\mathtt{B}_i}\centerdot S$. This semimartingale is special whose drift is $\chi^\cup$. 
\el

\textbf{Proof.} We use \cite[Corollaire(1.8)]{Jeulin80} to control the martingale part of $\ind_{\cup_{0\leq i\leq n}\mathtt{B}_i}\centerdot S$ by $\ind_{\cup_{0\leq i\leq n}\mathtt{B}_i}\centerdot [S,S]$, and the drift part of $\ind_{\cup_{0\leq i\leq n}\mathtt{B}_i}\centerdot S$ by $\chi^\cup$. We notice that $\ind_{\cup_{0\leq i\leq n}\mathtt{B}_i}$ is a left elementary process and hence the computation of $\ind_{\cup_{0\leq i\leq n}\mathtt{B}_i}\centerdot [S,S]$ is made by the quadratic variation which is independent of the filtration $\mathbb{G}$.\ok

\bl\label{BU+-}
Suppose that the distribution function $\chi^\cup$ exists and the process $C$ is the jump process of a special semimartingale. Then, the process $A$ is finite valued, and the three $\mathbb{G}$ semimartingales $\ind_{\mathtt{A}_\epsilon} \centerdot (S-\check{S})$, $\ind_{\mathtt{A}_\epsilon}\centerdot (S-\check{S})^+$ and $\ind_{\mathtt{A}_\epsilon}\centerdot (S-\check{S})^-$, $n\geq 0,$ converge in $\mathcal{H}^1_{\mathsf{loc}}$ to semimartingales that we denote respectively by $\ind_{\mathtt{A}} \centerdot (S-\check{S})$, $\ind_{\mathtt{A}}\centerdot (S-\check{S})^+$ and $\ind_{\mathtt{A}}\centerdot (S-\check{S})^-$. They are special semimartingales. We have $$
\dcb
\ind_{\mathtt{A}}\centerdot (S-\check{S})^+
=\ind_{\{(S-\check{S})_{-}>0\}}\ind_{\mathtt{A}} \centerdot (S-\check{S}) + \ind_{\mathtt{A}}\centerdot A+\frac{1}{2}\ind_{\mathtt{A}} \centerdot l^\cup\\
\\
\ind_{\mathtt{A}}\centerdot (S-\check{S})^-
=-\ind_{\{(S-\check{S})_{-}\leq 0\}}\ind_{\mathtt{A}} \centerdot (S-\check{S}) + \ind_{\mathtt{A}}\centerdot A+\frac{1}{2}\ind_{\mathtt{A}} \centerdot l^\cup
\dce
$$
for some continuous increasing process $l^\cup$ null at the origin which increases only on the set $\{t\in\mathtt{A}: (S-\check{S})_{t-}=0\}$
\el

\textbf{Proof.} 
Let the process $C$ be the jump process of $M+V$, where $M$ is a $\mathbb{G}$ local martingale and $V$ is a $\mathbb{G}$ predictable process with finite variation. By stopping $S,\check{S},M,V$ if necessary, we assume that $(S-\check{S})^*_\infty$ and $M^*$ are integrable, $(S-\check{S})_\infty$ exists, and the total variations of $\chi^\cup$ and of $V$ are integrable.

As the random measure $\mathsf{d}\chi^\cup$ has a distribution function $\chi^\cup$, the random measure $\ind_{(R,d_R]}\centerdot\mathsf{d}\chi^\cup$ also has a distribution function. The following limit computation is valid in $\mathcal{H}^1$: 
$$
\dcb
&&\ind_{(R,d_R]}\centerdot(\ind_{\cup_{i\geq 0}\mathtt{B}_i}\centerdot (S-\check{S}))\\
&=&\lim_{n\uparrow\infty}\ind_{(R,d_R]}\centerdot(\ind_{\cup_{0\leq i\leq n}\mathtt{B}_i}\centerdot (S-\check{S}))\\
&=&\lim_{n\uparrow\infty}(\ind_{(R,d_R]}\ind_{\cup_{0\leq i\leq n}\mathtt{B}_i})\centerdot (S-\check{S})\\
\dce
$$
Notice that the limit in $\mathcal{H}^1$ implies the limit with respect to the uniform convergence in compact intervals in probability. Lemma \ref{third-type} is applicable. For any $R<t<d_R$,$$
(\ind_{(R,d_R]}\centerdot(\ind_{\cup_{i\geq 0}\mathtt{B}_i}\centerdot (S-\check{S})))_t
=(\ind_{(R,d_R]}\centerdot (S-\check{S}))_t
$$
With Lemma \ref{elem-integral}, we compute the jumps at $d_R$ : $$
\Delta_{d_R}(\ind_{(R,d_R]}\centerdot(\ind_{\cup_{i\geq 0}\mathtt{B}_i}\centerdot (S-\check{S})))
=
\ind_{\{R<d_R\}}\ind_{\{d_R\in \cup_{i\geq 0}\mathtt{B}_i\}}\Delta_{d_R}(S-\check{S})
$$
These facts implies
$$
\dcb
\ind_{(R,d_R]}\centerdot (S-\check{S})=\ind_{(R,d_R]}\centerdot(\ind_{\cup_{i\geq 0}\mathtt{B}_i}\centerdot (S-\check{S}))+\ind_{\{R<d_R\}}\ind_{\{d_R\notin \cup_{i\geq 0}\mathtt{B}_i\}}\Delta_{d_R}(S-\check{S})\ind_{[d_R,\infty)}
\dce
$$
Note that $
[(d_{R})_{\{R<d_R,d_R\notin \cup_{i\geq 0}\mathtt{B}_i\}}]
=(R,d_R]\setminus \cup_{i\geq 0}\mathtt{B}_i
$
is a $\mathbb{G}$ predictable set. This means that the $\mathbb{G}$ stopping time $(d_{R})_{\{R<d_R,d_R\notin \cup_{i\geq 0}\mathtt{B}_i\}}$ is a predictable stopping time. Let $\hat{d}_R$ to denote this predictable stopping time. Consider the following jump process:
$$
\dcb
&&\ind_{\{R<d_R\}}\ind_{\{d_R\notin \cup_{i\geq 0}\mathtt{B}_i\}}\Delta_{d_R}(S-\check{S})\ind_{[d_R,\infty)}\\
&=&\ind_{\{R<d_R\}}\ind_{\{d_R\notin \cup_{i\geq 0}\mathtt{B}_i\}}C_{d_R}\ind_{[d_R,\infty)}\\

&=&\ind_{\{R<d_R\}}\ind_{\{d_R\notin \cup_{i\geq 0}\mathtt{B}_i\}}\Delta_{d_R}(M+V)\ind_{[d_R,\infty)}\\

&=&\Delta_{\hat{d}_{R}}(M+V)\ind_{[\hat{d}_{R},\infty)}\\

&=&\ind_{[\hat{d}_{R}]}\centerdot(M+V)\\
\dce
$$
Combining this equation with Lemma \ref{UB}, we see that $\ind_{(R,d_R]}\centerdot (S-\check{S})$ is a $\mathbb{G}$ special semimartingale whose drift is given by $\ind_{(R,d_R]}\centerdot\chi^\cup+\ind_{[\hat{d}_{R}]}\centerdot V$. 

Applying now the argument of the proof of Lemma \ref{UB}, we see that $\ind_{\mathtt{A}_\epsilon}\centerdot (S-\check{S})$ converges in $\mathcal{H}^1$ to a special semimartingale $\ind_{\mathtt{A}}\centerdot (S-\check{S})$.

Recall, for $t>0,$ $$
\ind_{(R,d_R]}\centerdot A_t = \sum_{R<s\leq d_R\wedge t} [1_{\{(S-\check{S})_{s-}>0\}}(S-\check{S})_{s}^{-}+1_{\{(S-\check{S})_{s-}\leq 0\}}(S-\check{S})_{s}^{+}]
= \sum_{R<s\leq d_R\wedge t,s\in\mathtt{j}(S-\check{S})} |S-\check{S}|_{s}.
$$
According to \cite[Chapter 9 \S 6]{HWY} applied to the semimartingale $(S-\check{S})^{d_R}\ind_{[R,\infty)}=\ind_{(R,d_R]}\centerdot (S-\check{S})+(S_R-\check{S}_R)\ind_{[R,\infty)}$,  
we know that the process 
$$
l^{(R,d_R]}_{t}(X) = 2[\ind_{(R,d_R]}\centerdot (S-\check{S})^{+}_t-\ind_{\{(S-\check{S})_{-}>0\}}\ind_{(R,d_R]} \centerdot (S-\check{S})_{t} - \ind_{(R,d_R]}\centerdot A_{t} ], \ t>0,
$$
is not decreasing, continuous and null at the origin, which increases only on the set $\{t\in(R,d_R]:(S-\check{S})_{t-}=0\}$. Note that, if we take another $\mathbb{G}$ stopping time $R'$, the random measure $\mathsf{d}l^{(R',d_{R'}]}$ coincide with the random measure $\mathsf{d}l^{(R,d_{R}]}$ on $(R,d_{R}]\cap(R',d_{R'}]$. Therefore, there exists a random measure $\mathsf{d}l^\cup$ on $\mathtt{A}$ such that $\ind_{(R,d_R]}\centerdot \mathsf{d}l^\cup=\mathsf{d}l^{(R,d_R]}$. (Note that $\mathsf{d}l^\cup$ is diffuse which does not charge $\mathtt{C}$.)

We have the following computation for any $\epsilon>0$:
$$
\dcb
&&\mathbb{E}[\int\ind_{\mathtt{A}_\epsilon}(s) (dl^\cup_s+2dA_s)]\\
&=&\mathbb{E}[\int\ind_{\cup_{n\geq 0}(R_n,d_{R_n}]}(s) (dl^\cup_s+2dA_s)]\\
&=&\lim_{k\uparrow\infty}\sum_{0\leq n\leq k}\mathbb{E}[l^{(R_n,d_{R_n}]}_{d_{R_n}}+2\ind_{(R_n,d_{R_n}]}\centerdot A_{d_{R_n}}]\\
&=&\lim_{k\uparrow\infty}\sum_{0\leq n\leq k}2\mathbb{E}[(S-\check{S})_{d_{R_n}}^{+}-(S-\check{S})_{R_n}^{+}-\ind_{\{(S-\check{S})_{-}>0\}}\ind_{(R_n,d_{R_n}]} \centerdot (S-\check{S})_{d_{R_n}} ]\\
&\leq&\lim_{k\uparrow\infty}\sum_{0\leq n\leq k}2\mathbb{E}[(S-\check{S})_{d_{R_n}}^{+}\ind_{\{R_n<\infty,d_{R_n}=\infty\}}+\ind_{\{(S-\check{S})_{-}>0\}}(\ind_{(R_n,d_{R_n}]} \centerdot |d\chi^\cup|_\infty +\ind_{[\hat{d}_{R}]}\centerdot |dV|_\infty)]\\
&\leq&2\mathbb{E}[\sum_{n\geq 0}(S-\check{S})_{d_{R_n}}^{+}\ind_{\{R_n<\infty,d_{R_n}=\infty\}}]+2\mathbb{E}[\int_0^\infty (|d\chi^\cup_s|+|dV_s|) ]\\
&\leq&2\mathbb{E}[(S-\check{S})^*_{\infty}]+2\mathbb{E}[\int_0^\infty (|d\chi^\cup_s|+|dV_s|) ]<\infty\\
\dce
$$
Here the last inequality is because there exist only one $n$ such that $R_n<\infty,d_{R_n}=\infty$. Let $\epsilon\rightarrow 0$. $\ind_{\mathtt{A}_\epsilon}(s)$ tends to $\ind_{\mathtt{A}}(s)$. We conclude that $
\mathbb{E}[\int_0^\infty \ind_{\mathtt{A}}(s) (dl^\cup_s+dA_s)]<\infty.
$
That means that $A$ is finite valued and $\mathsf{d}l^\cup$ have a distribution functions $l^\cup$. 

It is now straightforward to see that $
\ind_{\mathtt{A}_\epsilon}\centerdot (S-\check{S})^+
$
converge in $\mathcal{H}^1$ to a limit that we denote by $\ind_{\mathtt{A}}\centerdot (S-\check{S})^+$, which is equal to
$$
\ind_{\mathtt{A}}\centerdot (S-\check{S})^+
=\ind_{\{(S-\check{S})_{-}>0\}}\ind_{\mathtt{A}} \centerdot (S-\check{S}) + \ind_{\mathtt{A}}\centerdot A+\frac{1}{2}\ind_{\mathtt{A}} \centerdot l^\cup
$$ 
The first part of the lemma is proved.

The other part of the lemma can be proved similarly. Notice that we obtain the same random measure $\mathsf{d}l^\cup$ in the decomposition of $(S-\check{S})^-$\ok

\bl\label{Z-Ytheorem}
Suppose that the distribution function $\chi^\cup$ exists and the process $C$ is the jumps process of a $\mathbb{G}$ special semimartingale. Then, $S$ is a special semimartingale. More precisely, let \[
\begin{array}{lcl}
V^{+}&=&(S-\check{S})^{+}-(S-\check{S})_0^{+}-\ind_{\mathtt{A}}\centerdot(S-\check{S})^{+}\\
V^{-}&=&(S-\check{S})^{-}-(S-\check{S})_0^{-}-\ind_{\mathtt{A}}\centerdot(S-\check{S})^{-}
\end{array}
\]
$V^+, V^-$ are non decreasing locally integrable processes. They increase only outside of $\mathtt{A}$. We have :\[
\begin{array}{lll}
S&=&S_{0}+\ind_{\mathtt{A}}\centerdot (S-\check{S})+V^{+}-V^{-}+ \check{S}.
\end{array}
\]
\el

\textbf{Proof.} Let $X=(S-\check{S})^{+}$ and $X'=(S-\check{S})^{-}$. For a $\epsilon>0$, we compute the
difference :
\begin{eqnarray*}
&&X - X_{0}-\ind_{\mathtt{A}_{\epsilon}}\centerdot X = X - X^{0} - \sum_{n\geq 0} (X^{d_{R_{n}}}-X^{R_{n}})\\
&=& \sum_{n\geq 0} (X^{d_{R_{n}}}-X^{R_{n}})
+\sum_{n\geq 0} (X^{R_{n}}-X^{d_{R_{n-1}}}) - \sum_{n\geq 1} (X^{d_{R_{n}}}-X^{R_{n}})\\
&=& \sum_{n\geq 0} (X^{R_{n}}-X^{d_{R_{n-1}}})
\end{eqnarray*}
According to Lemma \ref{BU+-} $\ind_\mathtt{A}\centerdot X$ exists as the limit in $\mathcal{H}^1_{loc}$ of $\ind_{\mathtt{A}_{\epsilon}}\centerdot X$. Let $\epsilon$ tend to zero. The first term of the above identity tends to a process $V^+=X- X_{0} - \ind_{\mathtt{A}}\centerdot X$, uniformly on every compact interval in probability (in particular $V^+$ is càdlàg).

Consider the last term of the above identity. For $t>0$, let $N(t) = \sup\{k\geq 1; d_{R_{k}} \leq t\}$ ($N(t) = 0$ if the set is empty). Recall that, on the set $\{d_{R_{j}}<\infty\}$, $X_{d_{R_{j}}} = 0$. We have  
$$
\dcb
&&\sum_{n\geq 0} (X^{R_{n}}-X^{d_{R_{n-1}}})_t\\

&=&\sum_{n= 0}^{N(t)+1} X^{R_{n}}_{t}\\
&=&X^{R_{0}}_{t}+X^{R_{1}}_{t}+X^{R_{2}}_{t}+X^{R_{3}}_{t}+\ldots+X^{R_{N(t)}}_{t}+X^{R_{N(t)+1}}_{t}\\
&&\mbox{recall that $d_{R_{N(t)}}\leq t$ so that $R_{N(t)}\leq t$}\\
&=&X_{R_{0}}+X_{R_{1}}+X_{R_{2}}+X_{R_{3}}+\ldots+X_{R_{N(t)}}+X^{R_{N(t)+1}}_t\\
&=&\sum_{n=0}^{N(t)} X_{R_{n}} + X_{t}1_{\{d_{R_{N(t)}}\leq t < R_{N(t)+1}\}}+ X_{R_{N(t)+1}}1_{\{ R_{N(t)+1}\leq t\}}
\dce
$$ 
Notice that $X_{t}1_{\{d_{R_{N(t)}}\leq t < R_{N(t)+1}\}}=X_{t}1_{\{d_{R_{N(t)}}< t < R_{N(t)+1}\}}$ because $X_{d_{R_{N(t)}}}=0$. If $t \in \mathtt{A}$, for $\epsilon$ small enough, $t$
will belongs to
$\mathtt{A}_{\epsilon}$. As the interval $(d_{R_{N(t)}},R_{N(t)+1})$ is contained in the complementary of $\mathtt{A}_{\epsilon}$, we must have  $\ind_{\{d_{R_{N(t)}}< t
\leq R_{N(t)+1}\}} = 0$. If $t \notin \mathtt{A}$, by the Assumption \ref{partialcovering}, $X_{t} = 0$. In sum, for every $t>0$, there is a $\epsilon(t)$ such that, for $\epsilon < \epsilon(t)$, 
$$
\sum_{n\geq 0} (X^{R_{n}}-X^{d_{R_{n-1}}})_t
=\sum_{n=0}^{N(t)} X_{R_{n}} + X_{R_{N(t)+1}}\ind_{\{ R_{N(t)+1}\leq t\}}
=\sum_{n\geq 0} X_{R_{n}}\ind_{\{R_{n}\leq t\}}
$$
From this expression, we can write, for $0<s<t$,
$$
V^{+}_{t} - V^{+}_{s} 
=  \lim_{\epsilon \downarrow 0} (\sum_{n\geq 0} X_{R_{n}}\ind_{\{ R_{n}\leq t\}}-\sum_{n\geq 0} X_{R_{n}}\ind_{\{ R_{n}\leq s\}})
=  \lim_{\epsilon \downarrow 0} \sum_{n\geq 0} X_{R_{n}}\ind_{\{s<  R_{n}\leq t\}} \geq 0,
$$
i.e. $V^+$ is an increasing process. Moreover, for a fixed $a>0$, since $\mathtt{A}_a\subset \mathtt{A}_\epsilon$ for any $\epsilon<a$, since $R_n\notin\mathtt{A}_\epsilon$, $\ind_{\mathtt{A}_a}\centerdot V^+=0$. This argument shows that the random measure $\mathtt{d}V^+$ does not charge $\mathtt{A}$. Finally, as $X$ is locally in class$(D)$, $V^+$ is locally integrable. 

According to Lemma \ref{BU+-} 
$$
\ind_{\mathtt{A}}\centerdot X =\ind_{\{(S-\check{S})_{-}>0\}}\ind_{\mathtt{A}} \centerdot (S-\check{S}) +\ind_{\mathtt{A}}\centerdot A+\frac{1}{2}\ind_{\mathtt{A}}\centerdot l^\cup
$$
and, therefore,
$$
X=X_0+\ind_{\{(S-\check{S})_{-}>0\}}\ind_{\mathtt{A}} \centerdot (S-\check{S}) +\ind_{\mathtt{A}}\centerdot A+\frac{1}{2}\ind_{\mathtt{A}}\centerdot l^\cup+V^+
$$
In the same way we prove$$
X'=X'_0-\ind_{\{(S-\check{S})_{-}\leq 0\}}\ind_{\mathtt{A}} \centerdot (S-\check{S}) +\ind_{\mathtt{A}}\centerdot A+\frac{1}{2}\ind_{\mathtt{A}}\centerdot l^\cup+V^-
$$
Writing finally $S=X-X'+\check{S}$, we prove the theorem.  \ok

\

\section{A method to find local solutions}\label{alpha}

We now consider the problem of enlargement of filtration for a $\mathbb{F}$ semimartingale $X$. By a local solution, we mean the fact that, for a couple $T\leq U$ of $\mathbb{G}$ stopping times, $\ind_{(T,U]}\centerdot X$ is a $\mathbb{G}$ semimartingale. According to Theorem \ref{semimartingaleproperty}, if we success to collect enough local solutions, we will be able to solve the problem for $X$.

In this section we give a method to find local solutions. The result is taken from \cite{SongThesis}.

\subsection{The functions $h^u, u\geq 0$}\label{h-formula}

We need some topological properties behind the filtrations. We assume that $\Omega$ is a Polish space and $\mathcal{B}$ is its borel $\sigma$-algebra. We assume that there exists a filtration $\mathbb{F}^\circ=(\mathcal{F}^\circ_t)_{t\geq 0}$ of sub-$\sigma$-algebras of $\mathcal{B}$ such that $$
\mathcal{F}_t=\mathcal{F}^\circ_{t+}\vee\mathcal{N}, \ t\geq 0,
$$ 
where $\mathcal{N}$ is the family of $(\mathbb{P},\mathcal{B})$-negligible sets. We assume that there exists a measurable map $I$ from $\Omega$ into another polish space $E$ equipped with its borel $\sigma$-algebra and a filtration $\mathbb{I}=(\mathcal{I}_t)_{t\geq 0}$ of countably generated sub-borel-$\sigma$-algebras, such that $\mathcal{G}_{t}=\mathcal{G}^\circ_{t+}\vee \mathcal{N}, \ t\geq 0,$ where
$
\mathcal{G}^\circ_{t}=\mathcal{F}^\circ_{t} \vee I^{-1}(\mathcal{I}_t), t\geq 0.
$
We recall the following result (see \cite[Theorem 4.36]{HWY})

\bl\label{FF0}
For any $\mathbb{F}$ stopping time $T$, there exists a $\mathbb{F}^\circ_+$ stopping time $T^\circ$ such that $T=T^\circ$ almost surely. For any $A\in\mathcal{F}_T$, there exists a $A^\circ\in\mathcal{F}^\circ_{T^\circ+}$ such that $A\Delta A^\circ$ is negligible. For any $\mathbb{F}$ optional (resp. predictable) process $X$, there exists a $\mathbb{F}^\circ_+$ optional (resp. predictable) process $X^\circ$ such that $X$ and $X^\circ$ are indistinguishable. 

A similar result holds between the couple $\mathbb{G}$ and $\mathbb{G}^\circ_+=(\mathcal{G}^\circ_{t+})_{t\geq 0}$. 
\el

\brem
As a consequence of Lemma \ref{FF0}, we can speak about, for example, the predictable dual projection of an integrable increasing process with respect to the filtration $\mathbb{G}^\circ_+=(\mathcal{G}^\circ_{t+})_{t\geq 0}$ : this will mean that we compute the predictable dual projection in $\mathbb{G}$, then we take a version in $\mathbb{G}^\circ_+$.
\erem

Consider the product space $\Omega\times E$ equipped with its product $\sigma$-algebra and its product filtration $\mathbb{J}$ composed of 
$$
\mathcal{J}_t=\sigma\{A\times B: A\in\mathcal{F}^\circ_t, \ B\in\mathcal{I}_t\}, t\geq 0.
$$
We introduce the map $\phi$ on $\Omega$ such that
$\phi(\omega) = (\omega, I(\omega)) \in \Omega \times E$. Notice that, for $t\geq 0$, $\mathcal{G}^\circ_t=\phi^{-1}(\mathcal{J}_t)$. Therefore, for $C\in\mathcal{G}^\circ_{t+}$, there exist a sequence $(D_n)_{n\geq 1}$ of sets in respectively $\mathcal{J}_{t+\frac{1}{n}}$ such that $C=\phi^{-1}(D_n)$, which means $C=\phi^{-1}(\cap_{n\geq 1}\cup_{k\geq n}D_k)\in \phi^{-1}(\mathcal{J}_{t+})$. This observation yields the equality $\mathcal{G}^\circ_{t+}=\phi^{-1}(\mathcal{J}_{t+}), t\geq 0$.

We equip the product space $\Omega \times E$ with the the image probability $\mu$ of $\mathbb{P}$ by $\phi$. We introduce the identity map $\mathsf{i}$ on $\Omega$, the map $\zeta(\omega,x)=x$ and the map $\iota(\omega,x)=\omega$ for $(\omega,x)\in\Omega\times E$. For $t\geq 0$, let  
$$
\begin{array}{lll}
\pi_{t,F/F}(\omega, d\omega')
&=&
\mbox{regular conditional distribution of the map $\mathsf{i}$ under $\mathbb{P}$}
\\
&&\mbox{given the map $\mathsf{i}$ (itself) as a map valued in $(\Omega,\mathcal{F}^\circ_t)$}\\

\pi_{t,I/F}(\omega, dx')
&=&
\mbox{regular conditional distribution of $I$ under $\mathbb{P}$}
\\
&&\mbox{given the map $\mathsf{i}$ as a map valued in $(\Omega,\mathcal{F}^\circ_t)$}\\

\pi_{t,I/I}(x,dx')
&=&
\mbox{regular conditional distribution of $I$ under $\mathbb{P}$}
\\
&&\mbox{given $I$ (itself) as a map valued in $(E,\mathcal{I}_t)$}\\
\end{array}
$$

\brem
There exist situations where $\mathbb{F}^\circ$ is generated by a borel map $Y$ from $\Omega$ into a polish space $F$ equipped with a filtration $\mathbb{K}=(\mathcal{K}_t)_{t\geq 0}$ of sub-borel-$\sigma$-algebras on $F$, such that $\mathcal{F}^\circ_t=Y^{-1}(\mathcal{K}_t), t\geq 0$. Let $\check{\pi}_{t}(y,d\omega')$ be the regular conditional distribution of $\mathsf{i}$ under $\mathbb{P}$ given the map $Y$ considered as a map in the space $(F,\mathcal{K}_t)$. Then, $
\check{\pi}_t(Y(\omega),d\omega')
$
is a version of $\pi_{t,F/F}(\omega,d\omega')$, and its image measure by the map $I$ on $E$ is a version of $\pi_{t,I/F}(\omega,dx')$.
\erem

For a fixed $u \geq 0$, for any $t \geq u,$ let $h_{t}^{u}(\omega, x)$
to be a function on $\Omega \times E$ which is ${\mathcal{J}}_{t}$-measurable, such that, for $\mathbb{P}$-almost all $\omega$,
$$
h_{t}^{u}(\omega, x) = \left. {\frac{2 \pi_{t,I/F}(\omega, dx)}
{\pi_{t,I/F}(\omega, dx) +
\int \pi_{u,I/F}(\omega, dx'') \pi_{u,I/I}(x'', dx)}}\right|
_{\mathcal{I}_{t}},
$$

We introduce a family of probabilities $\nu^{u}$ on $\mathcal{J}_\infty$
indexed by
$u \geq 0$ determined by the equations :
\[
\int f(\omega, x) \nu^{u}(d\omega dx)
=
\mathbb{E}_{\mu}\left[ 
\int \pi_{u,F/F}(\iota, d\omega)
\int \pi_{u,I/I}(\zeta, dx) f(\omega, x)\right] 
\]
where $f(\omega,x)$ is a positive $\mathcal{J}_\infty$ measurable function. We notice that $\nu^u$ coincides with $\mu$ on $\mathcal{J}_u$, on $\iota^{-1}(\mathcal{F}^\circ_\infty)$, and on $\zeta^{-1}(\mathcal{I}_\infty)$. The kernel $\int \pi_{u,F/F}(\iota, d\omega)
\int \pi_{u,I/I}(\zeta, dx)$ is a regular conditional distribution of the identity map on $\Omega\times E$, under the probability
$\nu^{u}$ given $\mathcal{J}_u$.

\bl\label{martingaletransfer} 
Let $M$ to be a
$\mathbb{F}^\circ_+$-adapted càdlàg $\mathbb{P}$ integrable process. Then, $M$ is a
$(\mathbb{P},\mathbb{F})$-martingale on $[u,\infty)$, if and only if $M(\iota)$ is a $(\nu^u,\mathbb{J}_+)$ martingale on $[u,\infty)$.
\el

\textbf{Proof.} Note that, because of the càdlàg path property, $M$ is a
$(\nu^u,\mathbb{J}_+)$ martingale, if and only if
\[
\dcb
\mathbf{E}_{\nu^{u}}[M_{t}(\iota)\ind_A(\iota)\ind_B(\zeta)]
=\mathbf{E}_{\nu^{u}}[M_{s}(\iota)\ind_A(\iota)\ind_B(\zeta)] \ \ok
\dce 
\]
for any $t\geq s\geq u$ and for any $A\in\mathcal{F}^\circ_s$, $B\in\mathcal{I}_s$. This is equivalent to
\[
\dcb
&&\mathbf{E}_{\mu}\left[ 
\int \pi_{u,F/F}(\iota, d\omega)M_{t}(\omega)\ind_A(\omega)
\int \pi_{u,I/I}(\zeta, dx) \ind_B(x)\right]\\
&=&\mathbf{E}_{\mu}\left[ 
\int \pi_{u,F/F}(\iota, d\omega)M_{s}(\omega)\ind_A(\omega)
\int \pi_{u,I/I}(\zeta, dx) \ind_B(x)\right]\\
\dce 
\]
or to
\[
\dcb
&&\mathbf{E}_{\mathbb{P}}\left[ 
M_{t}\ind_A
\mathbb{E}[\int \pi_{u,I/I}(I, dx) \ind_B(x)|\mathcal{F}^\circ_u]\right]
=\mathbf{E}_{\mathbb{P}}\left[ 
M_{s}\ind_A
\mathbb{E}[\int \pi_{u,I/I}(I, dx) \ind_B(x)|\mathcal{F}^\circ_u]\right]\\
\dce 
\]
The last relation is true if and only if $M$ is a $(\mathbb{P},\mathbb{F})$-martingale on $[u,\infty)$. \ \ok

\bl\label{h+}
For $t
\geq u$, 
the function $h_{t}^{u}$ exists and it is a version of $\left.\frac{2\mu}{\mu+\nu^{u}}\right |_{{\mathcal{J}}_{t}}$ (respectively, $2-h_{t}^{u}=\left.\frac{2\nu^u}{\mu+\nu^{u}}\right |_{{\mathcal{J}}_{t}}$).
Consequently, the right limit version
$$
h_{t+}^{u} =
\lim_{\epsilon>0, \epsilon\rightarrow
0}\sup_{s \in \mathbb{Q}_{+}, u\leq s\leq u+\epsilon}  h(\omega, x),\ t\geq u,
$$
is well-defined and the process $h^u_+=(h^u_{t+})_{t\geq u}$ is a càdlàg $(\mu+\nu^{u},\mathbb{J}_+)$ uniformly integrable martingale on $[u,\infty)$. We have $0\leq h^u_{+}\leq 2$ and, if
$$
\dcb
\tau^{u}(\omega, x)
&=&
\inf\{t \geq u: h_{t+}^{u} = 2\}\\
\upsilon^{u}(\omega, x)
&=&
\inf\{t \geq u: h_{t+}^{u} = 0\}\\
\dce
$$ 
$h^u_{t+}=2$ for $t\in[\tau^u,\infty)$ and $h^u_{t+}=0$ for $t\in[\upsilon^u,\infty)$, $(\mu+\nu^u)$ almost surely. We have also $$
\mu[\upsilon^u<\infty]=0,\ \ \nu^u[\tau^u<\infty]=0.
$$
\el

\textbf{Proof.} We prove first of all two identities : Let $t\geq u$. Let $H\geq 0$ be $\mathcal{J}_{t}$ measurable,
\[
\begin{array}{lcl}
\mathbb{E}_{\mu}[H]&=&\mathbb{E}_\mathbb{P}[\int H(\mathsf{i},x)\pi_{t,I/F}(\mathsf{i}, dx)],\\
\mathbb{E}_{\nu^{u}}[H]&=&\mathbb{E}_{\mathbb{P}}[\int \pi_{u,I/F}(\mathsf{i}, dx'') \int H(\mathsf{i},x) \pi_{u,I/I}(x'', dx)].
\end{array}
\]
By monotone class theorem, it is enough to check them for $H=\ind_A(\iota)\ind_B(\zeta)$ where $A\in\mathcal{F}^\circ_t$ and $B\in\mathcal{I}_t$. But then, because$$
\ind_A(\iota)\int \ind_B(x)\pi_{t,I/F}(\iota, dx)
=\int \ind_A(\iota)\ind_B(x)\pi_{t,I/F}(\iota, dx),
$$ 
we have $$
\dcb
\mathbb{E}_{\mu}[\ind_A(\iota)\ind_B(\zeta)]
&=&\mathbb{E}_{\mathbb{P}}[\ind_A\ind_B(I)]\\
&=&\mathbb{E}_{\mathbb{P}}[\ind_A\mathbf{E}_{\mathbb{P}}[\ind_B(I)|\mathcal{F}_t]]\\
&=&\mathbb{E}_{\mathbb{P}}[\ind_A\int \ind_B(x)\pi_{t,I/F}(\cdot, dx)]\\
&=&\mathbb{E}_{\mathbb{P}}[\int \ind_A(\cdot)\ind_B(x)\pi_{t,I/F}(\cdot, dx)]
\dce
$$
and
$$
\dcb
\mathbb{E}_{\nu^{u}}[\ind_A(\iota)\ind_B(\zeta)]
&=&\mathbb{E}_{\mu}\left[ 
\int \pi_{u,F/F}(\iota, d\omega)
\int \pi_{u,I/I}(\zeta, dx) \ind_A(\omega)\ind_B(x)\right]\\
&=&\mathbb{E}_{\mathbb{P}}\left[ 
\int \pi_{u,F/F}(\mathsf{i}, d\omega)\ind_A(\omega)
\int \pi_{u,I/I}(I, dx) \ind_B(x)\right]\\
&=&\mathbb{E}_{\mathbb{P}}\left[ 
\int \pi_{u,F/F}(\mathsf{i}, d\omega)\ind_A(\omega)
\mathbb{E}_{\mathbb{P}}[\int \pi_{u,I/I}(I, dx) \ind_B(x)|\mathcal{F}_u]\right]\\
&=&\mathbb{E}_{\mathbb{P}}\left[\ind_A
\mathbb{E}_{\mathbb{P}}[\int \pi_{u,I/I}(I, dx) \ind_B(x)|\mathcal{F}_u]\right]\\
&=&\mathbb{E}_{\mathbb{P}}[\ind_A\int \pi_{u,I/F}(\mathsf{i}, dx'') \int \pi_{u,I/I}(x'', dx)\ind_B(x)]\\
&=&\mathbb{E}_{\mathbb{P}}[\int \pi_{u,I/F}(\mathsf{i}, dx'') \int \pi_{u,I/I}(x'', dx)\ind_A(\mathsf{i})\ind_B(x)] \ 
\dce
$$
which concludes the identities.

Let $h_{t}^{u}$ be the density function $\left.\frac{2\mu}{\mu+\nu^{u}}\right |_{{\mathcal{J}}_{t}}$. For $A\in\mathcal{F}^\circ_t$, $B\in\mathcal{I}_t$, we have$$
2\mathbb{E}_{\mu}[\ind_A\ind_B]=\mathbb{E}_{\mu}[\ind_A\ind_Bh^u_t]+\mathbb{E}_{\nu^u}[\ind_A\ind_Bh^u_t]
$$
Applying the two identities for $\mathbb{E}_{\mu}$ and for $\mathbb{E}_{\nu^u}$, the above relation becomes$$
\dcb
&&2\mathbb{E}_\mathbb{P}[\int \pi_{t,I/F}(\mathsf{i}, dx)\ind_A(\mathsf{i})\ind_B(x)]\\
&=&\mathbb{E}_\mathbb{P}[\int \pi_{t,I/F}(\mathsf{i}, dx)\ind_A(\mathsf{i})\ind_B(x)h^u_t(\mathsf{i},x)]+\mathbb{E}_{\mathbb{P}}[\int \pi_{u,I/F}(\mathsf{i}, dx'') \int  \pi_{u,I/I}(x'', dx)\ind_A(\mathsf{i})\ind_B(x)h^u_t(\mathsf{i},x)]\\
\dce
$$
or
$$
\dcb
&&2\mathbb{E}_\mathbb{P}[\ind_A(\mathsf{i})\int \pi_{t,I/F}(\mathsf{i}, dx)\ind_B(x)]\\

&=&\mathbb{E}_\mathbb{P}[\ind_A(\mathsf{i})\left(\int \pi_{t,I/F}(\mathsf{i}, dx)+\int \pi_{u,I/F}(\mathsf{i}, dx'') \int  \pi_{u,I/I}(x'', dx)\right)\ind_B(x)h^u_t(\mathsf{i},x)]
\dce
$$
When $A$ runs over all $\mathcal{F}^\circ_t$, the above equation becomes$$
2\int \pi_{t,I/F}(\mathsf{i}, dx)\ind_B(x)
=
\left(\int \pi_{t,I/F}(\mathsf{i}, dx)+\int \pi_{u,I/F}(\mathsf{i}, dx'') \int  \pi_{u,I/I}(x'', dx)\right)\ind_B(x)h^u_t(\mathsf{i},x)
$$
$\mathbb{P}$-almost surely. Since $\mathcal{I}_t$ is countably generated, we conclude$$
\left.2\int \pi_{t,I/F}(\mathsf{i}, dx)\right|_{\mathcal{I}_t}
=
\left.\left(\int \pi_{t,I/F}(\mathsf{i}, dx)+\int \pi_{u,I/F}(\mathsf{i}, dx'') \int  \pi_{u,I/I}(x'', dx)\right)h^u_t(\mathsf{i},x)\right|_{\mathcal{I}_t}
$$ 
$\mathbb{P}$-almost surely. This is the first part of the lemma.

The second statements on $h^u_+$ are valid, because $h^u_+$ and $2-h^u_+$ are positive $(\mu+\nu^{u},\mathbb{J}_+)$ supermartingales. For the last assertion, we have$$
\dcb
\mu[\upsilon^u<\infty]
&=&\mathbb{E}_{\frac{1}{2}(\mu+\nu^u)}[\ind_{\{\upsilon^u\leq t\}}h^u_{t+}]=0\\ 
\nu^u[\tau^u<\infty]
&=&\mathbb{E}_{\frac{1}{2}(\mu+\nu^u)}[\ind_{\{\tau^u\leq t\}}(2-h^u_{t+})]=0.\ \ok
\dce
$$

Along with the function $h^u$, we introduce two processes on $[u,\infty)$:
\begin{eqnarray*}
\alpha_{t}^{u}(\omega, x) 
&=&
\frac{h_{t+}^{u}(\omega, x)}
{2-h_{t+}^{u}(\omega, x)} \ind_{\{t < \tau^{u}(\omega, x)\}}\\
\beta_{t}^{u}(\omega, x) 
&=&
\frac{2-h_{t+}^{u}(\omega, x)}{h_{t+}^{u}(\omega, x)} \ind_{\{t < \upsilon^{u}(\omega, x)\}}
\end{eqnarray*}

\subsection{Results based on the process $\alpha^u$}

In this subsection we study the process $\alpha^u$. 

\bl\label{basicformulas}
Fix $u\geq 0$. Let $\rho$ to be a $\mathbb{J}_{+}$-stopping time such that $u\leq \rho\leq \tau^u$. We have, for any
$\mathbb{J}_{+}$-stopping time
$\kappa$ with $u\leq \kappa$, for any $A \in
{\mathcal{J}}_{\kappa+}$,$$
\mathbb{E}_{\mu}[\ind_A\ind_{\{\kappa < \rho\}}] = \mathbb{E}_{\nu^{u}}[\ind_A\ind_{\{\kappa < \rho\}}
\alpha_{\kappa}^{u}]
$$
Consequently, $\ind_{[0,u)}+\alpha^{u}\ind_{[u,\rho)}$ is a $(\nu^u,{\mathbb{J}}_{+})$ supermartingale. Moreover, $\alpha^u>0$ on $[u,\rho)$ under the probability $\mu$. For any positive $\mathbb{J}_+$ predictable process $H$, $$
\mathbb{E}_{\mu}[H_\kappa\ind_{\{\kappa < \rho\}}] = \mathbb{E}_{\nu^{u}}[H_\kappa\ind_{\{\kappa < \rho\}}
\alpha_{\kappa}^{u}]
$$
Suppose in addition that $(\alpha^u\ind_{[u,\rho)})^\kappa$ is of class $(D)$ under $\nu^u$. Let $V$ be the increasing $\mathbb{J}_+$ predictable process associated with the supermartingale $\ind_{[0,u)}+\alpha^{u}\ind_{[u,\rho)}$ (see \cite[p.115 Theorem 13]{protter} and Lemma \ref{FF0}). For any positive $\mathbb{J}_+$ predictable process $H$, we have$$
\mathbb{E}_{\mu}[H_{\rho}\ind_{\{u<\rho\leq \kappa\}}] =
\mathbb{E}_{\nu^{u}}[\int_{u}^{\kappa}
H_s dV_s]
$$
Let $B$ be a $\mathbb{J}_+$ predictable process with bounded total variation. We have$$
\mathbb{E}_{\mu}[B_{\kappa\wedge \rho}-B_u] =
\mathbb{E}_{\nu^{u}}[\int_{u}^{\kappa}
\alpha^u_{s-} dB_s]
$$
Consequently, $\mathbb{E}_{\mu}[\int_u^{\kappa}\ind_{\{u<s\leq \rho\}}\ind_{\{\alpha^u_{s-}=0\}}dB_s]=0$. Let $C$ be a $\mathbb{J}_+$ optional process. Suppose that the random measure $\mathsf{d}C$ on the open random interval $(u,\kappa\wedge \rho)$ has bounded total variation. For any bounded $\mathbb{J}_+$ predictable process $H$, We have$$
\mathbb{E}_{\mu}[\int_0^\infty H_s\ind_{(u,\kappa\wedge \rho)}(s) dC_s] =
\mathbb{E}_{\nu^{u}}[\int_{0}^\infty H_s\ind_{(u,\kappa\wedge \rho)}(s) \alpha^u_sdC_s]
$$
In particular, $$
(\ind_{(u,\kappa\wedge \rho)}\centerdot C)^{\mu\cdot\mathbb{J}_+-p}
=
\frac{1}{\alpha^u_-}\centerdot(\ind_{(u,\kappa\wedge \rho)}\alpha^u\centerdot C)^{\nu^u\cdot\mathbb{J}_+-p}
$$
\el

\textbf{Proof.} Note that $h^u_{\kappa+}$ is a version of $\left.\frac{2\mu}{\mu+\nu^u}\right|_{\mathcal{J}_{\kappa+}}$. From this relation, it results that, for any positive $\mathcal{J}_{\kappa+}$ measurable function $f$,$$
2\mathbb{E}_{\mu}[f \ind_{\{\kappa<\rho\}}]=\mathbb{E}_{\mu}[f \ind_{\{\kappa<\rho\}}h^u_{\kappa+}]+\mathbb{E}_{\nu^u}[f \ind_{\{\kappa<\rho\}}h^u_{\kappa+}]
$$
or equivalently
$$
\mathbb{E}_{\mu}[f (2-h^u_{\kappa+})\ind_{\{\kappa<\rho\}}]=\mathbb{E}_{\nu^u}[f (2-h^u_{\kappa+})\alpha_{\kappa}\ind_{\{\kappa<\rho\}}]
$$
This last identity is an equivalent form of the first formula of the lemma. To see the supermartingale property of $\ind_{[0,u)}+\alpha^{u}\ind_{[u,\rho)}$, it is enough to notice that $\mu=\nu^u$ on $\mathcal{J}_u$ and, for $A\in\mathcal{J}_{s+}$, $$
\dcb
\mathbb{E}_{\nu^u}[\ind_A\alpha^u_t\ind_{\{t<\rho\}}]
=\mathbb{E}_{\mu}[\ind_A\ind_{\{t<\rho\}}]
\leq\mathbb{E}_{\mu}[\ind_A\ind_{\{s<\rho\}}]
=\mathbb{E}_{\nu^u}[\ind_A\alpha^u_s\ind_{\{s<\rho\}}], \mbox{ for $u\leq s\leq t$}\\

\mathbb{E}_{\nu^u}[\ind_A\alpha^u_t\ind_{\{t<\rho\}}]
=\mathbb{E}_{\mu}[\ind_A\ind_{\{t<\rho\}}]
\leq\mathbb{E}_{\mu}[\ind_A\ind_{\{s<\rho\}}]
=\mathbb{E}_{\mu}[\ind_A]
=\mathbb{E}_{\nu^u}[\ind_A], \mbox{ for $0\leq s<u\leq  t$}\\

\mathbb{E}_{\nu^u}[\ind_A]
=\mathbb{E}_{\nu^u}[\ind_A], \mbox{ for $0\leq s<t<u$}
\dce
$$
The positivity of $\alpha^u$ on $[u,\rho)$ is the consequence of $\mu[\upsilon^u<\infty]=0$ by Lemma \ref{h+}.

The second formula of the lemma is a direct consequence of the first one. To prove the third formula of the lemma, we need only to check it on the processes $H$ of the form $\ind_A\ind_{(a,\infty)}$ with $0\leq a<\infty$ and $A\in\mathcal{J}_{a}$, and $\kappa\geq u$. We have 
$$
\dcb
\mathbb{E}_{\mu}[H_\rho\ind_{\{u<\rho\leq \kappa\}}]
&=&\mathbb{E}_{\mu}[\ind_A\ind_{\{a<\rho\}}\ind_{\{u<\rho\leq \kappa\}}]\\
&=&\mathbb{E}_{\mu}[\ind_A\ind_{\{a<\rho\}}\ind_{\{u<\rho\}}]
-\mathbb{E}_{\mu}[\ind_A\ind_{\{a<\rho\}}\ind_{\{\kappa<\rho\}}]\\
&=&\mathbb{E}_{\nu^u}[\ind_A\alpha^u_{a\vee u}\ind_{\{a<\rho\}}\ind_{\{u<\rho\}}]
-\mathbb{E}_{\nu^u}[\ind_A\alpha^u_{a\vee \kappa}\ind_{\{a<\rho\}}\ind_{\{\kappa<\rho\}}]\\

&=&\mathbb{E}_{\nu^u}[\ind_A\left(\alpha^u_{a\vee u}\ind_{\{a\vee u<\rho\}}-\alpha^u_{a\vee \kappa}\ind_{\{a\vee \kappa<\rho\}}\right)]
\\

&=&\mathbb{E}_{\nu^u}[\ind_A\left(V_ {a\vee \kappa}-V_ {a\vee u}\right)]\\

&=&\mathbb{E}_{\nu^u}[\ind_A\int_u^\kappa \ind_{(a,\infty)}(s)dV_s]\\

&=&\mathbb{E}_{\nu^u}[\int_u^\kappa H_sdV_s]
\dce
$$
For the fourth formula, we write$$
\dcb
\mathbb{E}_{\mu}[B_{\kappa\wedge\rho}-B_u] 
&=&\mathbb{E}_{\mu}[(B_{\kappa}-B_u)\ind_{\{\kappa<\rho\}}]
+\mathbb{E}_{\mu}[(B_{\rho}-B_u)\ind_{\{\kappa\geq \rho\}}]\\
&=&\mathbb{E}_{\nu^u}[(B_{\kappa}-B_u)\alpha^u_\kappa\ind_{\{\kappa<\rho\}}]
+\mathbb{E}_{\mu}[(B_{\rho}-B_u)\ind_{\{u<\rho\leq \kappa\}}]\\
&=&\mathbb{E}_{\nu^u}[(B_{\kappa}-B_u)\alpha^u_\kappa\ind_{\{\kappa<\rho\}}]
+\mathbb{E}_{\nu^u}[\int_u^\kappa(B_{s}-B_u)dV_s]\\
&=&\mathbb{E}_{\nu^{u}}[\int_{u}^{\kappa}
\alpha^u_{s-} dB_s]
\dce
$$
Finally, for the last formulas, let $K=H\ind_{(u,\kappa\wedge \rho)}\centerdot C$. We note that $K_\kappa=K_{\kappa\wedge \rho-}$ and $K_s=0, \forall s\leq u$. Applying the preceding formulas, we can write $$
\dcb
\mathbb{E}_{\mu}[\int_0^\infty H_s\ind_{(u,\kappa\wedge \rho)}(s) dC_s] 
&=&\mathbb{E}_{\mu}[K_{\kappa\wedge \rho-}]\\
&=&\mathbb{E}_{\nu^{u}}[K_{\kappa-}\ind_{\{\kappa<\rho\}}\alpha^u_\kappa+\int_{u}^\kappa K_{s-}dV_s]\\
&=&\mathbb{E}_{\nu^{u}}[K_{\kappa}\ind_{\{\kappa<\rho\}}\alpha^u_\kappa+\int_{u}^\kappa K_{s-}dV_s]\\
&=&\mathbb{E}_{\nu^{u}}[\int_u^\kappa K_{s-}d(\alpha^u\ind_{[u,\rho)})_s+\int_u^\kappa (\alpha^u\ind_{[u,\rho)})_{s}dK_s+\int_{u}^\kappa K_{s-}dV_s]\\
&=&\mathbb{E}_{\nu^{u}}[\int_u^\kappa (\alpha^u\ind_{[u,\rho)})_{s}dK_s]\\
&&\mbox{ because $(\alpha^u\ind_{[u,\rho)})^\kappa$ is of class$(D)$}\\
&=&\mathbb{E}_{\nu^{u}}[\int_0^\infty \alpha^u_s\ind_{(u,\kappa\wedge \rho)}(s)H_sdC_s]

\dce
$$
We can write again
$$
\dcb
\mathbb{E}_{\mu}[\int_0^\infty H_s d(\ind_{(u,\kappa\wedge \rho)}\centerdot C)^{\mu\cdot\mathbb{J}_+-p}_s] 
&=&\mathbb{E}_{\mu}[\int_0^\infty H_s \ind_{\{\alpha^u_{s-}>0\}} d(\ind_{(u,\kappa\wedge \rho)}\centerdot C)^{\mu\cdot\mathbb{J}_+-p}_s]\\

&=&\mathbb{E}_{\mu}[\int_0^\infty H_s\ind_{\{\alpha^u_{s-}>0\}}\ind_{(u,\kappa\wedge \rho)}(s) dC_s] \\
&=&\mathbb{E}_{\nu^{u}}[\int_0^\infty H_s\ind_{\{\alpha^u_{s-}>0\}}\alpha^u_s\ind_{(u,\kappa\wedge \rho)}(s)dC_s]\\
&=&\mathbb{E}_{\nu^{u}}[\int_0^\infty H_s\ind_{\{\alpha^u_{s-}>0\}}d(\alpha^u\ind_{(u,\kappa\wedge \rho)}\centerdot C)^{\nu^u\cdot\mathbb{J}_+-p}_s]\\

&=&\mathbb{E}_{\nu^{u}}[\int_0^\infty \alpha^u_{s-}\frac{1}{\alpha^u_{s-}}H_s\ind_{\{\alpha^u_{s-}>0\}}d(\alpha^u\ind_{(u,\kappa\wedge \rho)}\centerdot C)^{\nu^u\cdot\mathbb{J}_+-p}_s]\\

&=&\mathbb{E}_{\mu}[\int_0^\infty H_s\frac{1}{\alpha^u_{s-}}d(\alpha^u\ind_{(u,\kappa\wedge \rho)}\centerdot C)^{\nu^u\cdot\mathbb{J}_+-p}_s]\ \ok

\dce
$$

\bethe
\label{theorem:our-formula}
Let $M$ be a
bounded
$(\mathbb{P},\mathbb{F})$ martingale (assumed to be $\mathbb{F}^\circ_+$ adapted). Let $u\geq 0$. Let $\rho$ to be a $\mathbb{J}_{+}$-stopping time such that $u\leq \rho\leq \tau^u$. Let $\langle M(\iota),\alpha^{u}\ind_{[u,\rho)} \rangle^{\nu^u\cdot\mathbb{J}_+}$ denote the $\mathbb{J}_+$-predictable
bracket of the couple of $M(\iota)$ and $\ind_{[0,u)}+\alpha^{u}\ind_{[u,\rho)}$, computed under the probability $\nu^{u}$ in
the filtration ${\mathbb{J}}_{+}$. Then, there exists an increasing sequence $(\eta_n)_{n\geq 1}$ of $\mathbb{J}_+$ stopping times such that $\sup_{n\geq 1}\eta_n(\phi)\geq \rho(\phi)$ $\mathbb{P}$-almost surely, and for every $n\geq 1$, $M$ is a special
$\mathbb{G}$-semimartingale on the random left interval 
$(u,\eta_n\circ\phi]$ such that the process 
\[
(M^{\eta_n(\phi)}-M^{u}) - 
\ind_{(u,\eta_n(\phi)]}\frac{1}{\alpha^u(\phi)_{-}}\centerdot \langle M(\iota),
\alpha^{u}\ind_{[u,\rho)}
\rangle^{\nu^u\cdot\mathbb{J}_+}\circ \phi - 
\ind_{\{u<\rho(\phi)\leq \eta_n(\phi)\}}\Delta_{\rho(\phi)}M \ind_{[\rho(\phi),\infty)}
\]
is a $(\mathbb{P},\mathbb{G})$ local martingale.

If $\ind_{(u,\rho(\phi)]}\frac{1}{\alpha^u(\phi)_{-}}\centerdot \langle M(\iota),
\alpha^{u}\ind_{[u,\rho)}
\rangle^{\nu^u\cdot\mathbb{J}_+}\circ \phi$ has a distribution function under $\mathbb{P}$, the above decomposition formula remains valid if we replace $\eta_n$ by $\rho$.
\ethe

\brem
Notice that, according to \cite[Proposition II.2.4]{SongThesis}, different version of $\langle M(\iota),
\alpha^{u}\ind_{[u,\rho)}
\rangle^{\nu^u\cdot\mathbb{J}_+}$ under $\nu^u$ give rise of indistinguishable versions of $\langle M(\iota),
\alpha^{u}\ind_{[u,\rho)}
\rangle^{\nu^u\cdot\mathbb{J}_+}\circ \phi$ \ on $(u,\tau^u(\phi)]$ under $\mathbb{P}$. 
\erem

\textbf{Proof.} 
In this proof, instead of
$\nu^{u}$ (resp. $\tau^{u}, \upsilon^u, \alpha^{u}, \langle M(\iota),
\alpha^{u}\ind_{[u,\rho)}
\rangle^{\nu^u\cdot\mathbb{J}_+}$), we shall write
$\nu$ (resp.
$\tau, \upsilon,
\alpha, \langle M(\iota),
\alpha\ind_{[u,\rho)}
\rangle$).
For each integer $n
\geq 1$, set $$
\dcb
\gamma^+_{n} = \inf\{t \geq u:
\alpha_{t}>n\}\ \mbox{ and }\
\gamma^-_{n} = \inf\{t \geq u:
\alpha_{t}<1/n\},\ n\geq 1.
\dce
$$ 
Let $(\gamma^\circ_n)_{n\geq 1}$ be a sequence of $\mathbb{J}_{+}$-stopping times tending to $\infty$ under $\nu$ such that $\gamma^\circ_n\geq u$ for every $n\geq 1$, and $\cro{M(\iota),\alpha\ind_{[u,\rho)}}^{\gamma^\circ_n}$ has bounded total variation. 
Let $\eta_n=\gamma^+_{n}\wedge\gamma^-_{n}\wedge\gamma^\circ_{n}\wedge\rho$. Note that $\gamma^+_\infty=\tau$ on $\{\gamma^+_\infty<\infty\}$, $\gamma^-_\infty=\tau\wedge \upsilon$ under $\mu+\nu$. We have $$
\dcb
\mathbb{P}[\eta_\infty(\phi)<\rho(\phi)]
=\mu[\eta_\infty<\rho]
=\mathbb{E}_\nu[\ind_{\{\eta_\infty<\rho\}}\alpha_{\eta_\infty}]
=\mathbb{E}_\nu[\ind_{\{\tau\wedge\upsilon<\rho\}}\alpha_{\tau\wedge\upsilon}]
=0.
\dce
$$

Fix an $n\geq 1$. Let $\kappa$ be a $\mathbb{J}_{+}$-stopping time such that $u
\leq \kappa \leq \eta_{n}$. The stopping time $\kappa$ satisfies therefore the condition of Lemma \ref{basicformulas}. Recall also Lemma \ref{martingaletransfer} according to which $M(\iota)$ is a $(\nu^u,\mathbb{J}_+)$ martingale on $[u,\infty)$. 

Set $M^\flat=M-\Delta_{\rho(\phi)}M\ind_{[\rho(\phi),\infty)}$. Applying Lemma \ref{basicformulas}, we can write
\[
\begin{array}{lll}
&&\mathbb{E}_{\mathbb{P}}[(M^\flat_{\kappa(\phi)}- M^\flat_{u})]
\\
&=&\mathbb{E}_{\mu}[(M^\flat(\iota)_\kappa- M^\flat(\iota)_u)]
\\ 
&=& \mathbb{E}_{\mu}[(M(\iota)_\kappa- M(\iota)_u)\ind_{\{\kappa < \rho\}}] +
\mathbb{E}_{\mu}[(M^\flat(\iota)_\kappa- M^\flat(\iota)_u)1_{\{\kappa = \rho\}}] 
\\
&=&
\mathbb{E}_{\nu}[(M(\iota)_\kappa- M(\iota)_u)\alpha_{\kappa}\ind_{\{\kappa < \rho\}}] 
+
\mathbb{E}_{\mu}[(M(\iota)_{\rho-}- M(\iota)_u)\ind_{\{u<\kappa=\rho\}}]\\
&=&
\mathbb{E}_{\nu}[(M(\iota)_\kappa- M(\iota)_u)\alpha_{\kappa}\ind_{\{\kappa < \rho\}}] 
+
\mathbb{E}_{\nu}[\int_{u}^{\kappa}(M(\iota)_{-}- M(\iota)_{u})
dV]\\
&=&
\mathbb{E}_{\nu}[\int_{u}^{\kappa}
d\langle M(\iota), \alpha\ind_{[u,\rho)} \rangle_s]\\
&=&
\mathbb{E}_{\mu}[\int_{u}^{\kappa} \frac{1}{\alpha_{s-}}
d\langle M(\iota), \alpha\ind_{[u,\rho)} \rangle_s]\\
&&\mbox{ The above process is well-defined}\\
&&\mbox{ and it has bounded total variation, because $\frac{1}{\alpha_-}\leq n$ on $(u,\kappa]$}\\

&=&
\mathbb{E}_{\mathbb{P}}[(\int_{u}^{\kappa(\phi)} \frac{1}{\alpha(\phi)_{s-}}
d(\langle M(\iota), \alpha\ind_{[u,\rho)} \rangle\circ\phi)_s]
\dce
\]
This identity on $\mathbf{E}_{\mathbb{P}}[(M^\flat_{\kappa(\phi)}- M^\flat_{u})]$ implies that  (cf. \cite[Theorem 4.40]{HWY}) 
\[
(M^{\eta_n(\phi)}-M^{u}) - 
\ind_{(u,\eta_n(\phi)]}\frac{1}{\alpha(\phi)_{-}}\centerdot \langle M(\iota),
\alpha\ind_{[u,\rho)}
\rangle\circ \phi - 
\ind_{\{u<\rho(\phi)\leq \eta_n(\phi)\}}\Delta_{\rho(\phi)}M \ind_{[\rho(\phi),\infty)}
\]
is a $(\mathbb{P},\mathbb{G})$ local martingale. The first part of the theorem is proved. 

Suppose now that $\ind_{(u,\rho(\phi)]}\frac{1}{\alpha(\phi)_{-}}\centerdot \langle M(\iota),
\alpha\ind_{[u,\rho)}
\rangle\circ \phi$ has a distribution function $C=(C_t)_{t\geq 0}$ under $\mathbb{P}$. It is clear that $\ind_{(u,\eta_n(\phi)]}\centerdot M^\flat$ converges in $\mathcal{H}^1(\mathbb{P},\mathbb{G})$. Note that $\eta_\infty=\rho$ under $\mu$. Look at $\lim_{n\uparrow\infty}(M^\flat_{\eta_n}-M^\flat_u)$. If $\eta_n<\rho$ for all $n\geq 1$, this limit is equal to $(M^\flat_{\rho-}-M^\flat_u)=(M^\flat_{\rho}-M^\flat_u)$. If $\eta_n=\rho$ for some $n\geq 1$, this limit is equal to $(M^\flat_{\rho}-M^\flat_u)$ also. Look at the limit $$
\lim_{n\uparrow\infty}\ind_{(u,\eta_n(\phi)]}\frac{1}{\alpha(\phi)_{-}}\centerdot \langle M(\iota),
\alpha\ind_{[u,\rho)}
\rangle\circ \phi
=\ind_{\cup_{n\geq 1}(u,\eta_n(\phi)]}\frac{1}{\alpha(\phi)_{-}}\centerdot \langle M(\iota),
\alpha\ind_{[u,\rho)}
\rangle\circ \phi
$$
We note that $(u,\rho]\setminus \cup_{n\geq 1}(u,\eta_n]=[\rho_{\{\forall n\geq 1, \eta_n<\rho\}}]$. The time $\rho_{\{\forall n\geq 1, \eta_n<\rho\}}$ is a $\mathbb{J}_+$ predictable stopping time. So, on the set $\{\rho_{\{\forall n\geq 1, \eta_n<\rho\}}<\infty\}$, $$
\dcb
&&\ind_{(u,\rho]\setminus \cup_{n\geq 1}(u,\eta_n]}\centerdot\langle M(\iota), \alpha\ind_{[u,\rho)} \rangle\\
&=&\Delta_{\rho_{\{\forall n\geq 1, \eta_n<\rho\}}}\langle M(\iota), \alpha\ind_{[u,\rho)} \rangle\\

&=&\mathbb{E}_{\nu}[\Delta_{\rho_{\{\forall n\geq 1, \eta_n<\rho\}}} M(\iota) \Delta_{\rho_{\{\forall n\geq 1, \eta_n<\rho\}}}(\alpha\ind_{[u,\rho)})|\mathcal{J}_{\rho_{\{\forall n\geq 1, \eta_n<\rho\}}-}]\\

&=&\mathbb{E}_{\nu}[\ind_{\{\forall n\geq 1, \eta_n<\rho<\infty\}}\Delta_{\rho} M(\iota) \Delta_{\rho}(\alpha\ind_{[u,\rho)})|\mathcal{J}_{\rho_{\{\forall n\geq 1, \eta_n<\rho\}}-}]\\

&=&\mathbb{E}_{\nu}[\ind_{\{\forall n\geq 1, \eta_n<\rho<\infty\}}\Delta_{\rho} M(\iota) (-\alpha_{\rho-})|\mathcal{J}_{\rho_{\{\forall n\geq 1, \eta_n<\rho\}}-}]\\

&=&\mathbb{E}_{\nu}[\Delta_{\rho_{\{\forall n\geq 1, \eta_n<\rho\}}} M(\iota) (-\alpha_{\rho_{\{\forall n\geq 1, \eta_n<\rho\}}-})|\mathcal{J}_{\rho_{\{\forall n\geq 1, \eta_n<\rho\}}-}]\ind_{\{\rho_{\{\forall n\geq 1, \eta_n<\rho\}}<\infty\}}\\

&=&\mathbb{E}_{\nu}[\Delta_{\rho_{\{\forall n\geq 1, \eta_n<\rho\}}} M(\iota) |\mathcal{J}_{\rho_{\{\forall n\geq 1, \eta_n<\rho\}}-}](-\alpha_{\rho_{\{\forall n\geq 1, \eta_n<\rho\}}-})\ind_{\{\rho_{\{\forall n\geq 1, \eta_n<\rho\}}<\infty\}}\\

&=&0\
\mbox{ because $\rho_{\{\forall n\geq 1, \eta_n<\rho\}}$ is predictable and $M$ is a bounded martingale}

\dce
$$
We have actually
$$
\lim_{n\uparrow\infty}\ind_{(u,\eta_n(\phi)]}\frac{1}{\alpha(\phi)_{-}}\centerdot \langle M(\iota),
\alpha\ind_{[u,\rho)}
\rangle\circ \phi
=\ind_{(u,\rho(\phi)]}\frac{1}{\alpha(\phi)_{-}}\centerdot \langle M(\iota),
\alpha\ind_{[u,\rho)}
\rangle\circ \phi
$$
The theorem is proved.
\ok

The following lemma will also be useful.

\bl\label{alpha} Let $u \geq 0$. 
Let $\lambda$ be a ${\mathbb{J}_+}$-stopping time $\geq u$. Let $B$ be  a
$\mathcal{J}_{u}$-measurable set. If  for any $t\geq u$,  there is a
$\mathcal{J}_{t}$-measurable non negative function $l_{t}$ such that, $\mathbb{P}$-almost surely on
$\omega$
\[
\dcb
&&\int \pi_{t,I/F}(\omega, dx) f(x) \ind_{B}(\omega,x)\ind_{\{t \leq  \lambda(\omega,x)\}} \\
&=& \int \pi_{u,I/F}(\omega, dx'') \int\pi_{u,I/I}(x'', dx) f(x)\ind_{B}(\omega,x)\ind_{\{t \leq  \lambda(\omega,x)\}} l_{t}(\omega,x), 
\dce
\]
for all bounded $\mathcal{I}_{t}$-measurable function $f$.
Then, \[\tau^{u}\ind_{B} \geq \lambda \ind_{B}, \ \  
\alpha_{t}^{u}\ind_{B}\ind_{\{t < \lambda\}} = \limsup_{t'\in\mathbb{Q}_+,t'\downarrow t} l_{t'}\ind_{B}\ind_{\{t < \lambda\}}, t\geq u.\]
\el

Note that, if $\lambda$ is a $\mathbb{J}$ stopping time, we can replace $\ind_{\{t\leq\lambda\}}$ par $\ind_{\{t<\lambda\}}$ in the assumption.

\textbf{Proof.} We have
$$
\dcb
&&\left(\int \pi_{t,I/F}(\omega, dx)
+ \int \pi_{u,I/F}(\omega, dx'') \int\pi_{u,I/I}(x'', dx)\right) f(x)\ind_{B}(\omega,x)\ind_{\{t \leq \lambda(\omega,x)\}}\\
&=&\int \pi_{u,I/F}(\omega, dx'') \int\pi_{u,I/I}(x'', dx) f(x)\ind_{B}(\omega,x)\ind_{\{t \leq \lambda(\omega,x)\}} (l_{t}(\omega,x)+1)
\dce
$$
It follows that (noting that $\{t \leq \lambda\}\in\mathcal{J}_t$)
$$
\dcb
&&\left(\int \pi_{t,I/F}(\omega, dx)
+ \int \pi_{u,I/F}(\omega, dx'') \int\pi_{u,I/I}(x'', dx)\right) f(x)\ind_{B}(\omega,x)\ind_{\{t \leq \lambda(\omega,x)\}}h^u_t(\omega,x)\\
&=&2\int \pi_{t,I/F}(\omega, dx) f(x) \ind_{B}(\omega,x)\ind_{\{t \leq \lambda(\omega,x)\}}\\
&=& \int \pi_{u,I/F}(\omega, dx'') \int\pi_{u,I/I}(x'', dx) f(x)\ind_{B}(\omega,x)\ind_{\{t \leq \lambda(\omega,x)\}} 2l_{t}(\omega,x)\\
&=& \int \pi_{u,I/F}(\omega, dx'') \int\pi_{u,I/I}(x'', dx) f(x)\ind_{B}(\omega,x)\ind_{\{t \leq \lambda(\omega,x)\}} \frac{2l_{t}(\omega,x)}{l_{t}(\omega,x)+1}(l_{t}(\omega,x)+1)\\
&=&\left(\int \pi_{t,I/F}(\omega, dx)
+ \int \pi_{u,I/F}(\omega, dx'') \int\pi_{u,I/I}(x'', dx)\right) f(x)\ind_{B}(\omega,x)\ind_{\{t \leq \lambda(\omega,x)\}}\frac{2l_{t}(\omega,x)}{l_{t}(\omega,x)+1}
\dce
$$
The above identity means
$$
h_{t}^{u}(\omega,x) \ind_{B}(\omega,x)\ind_{\{t \leq \lambda(\omega,x)\}} = 
\frac{2l_{t}(\omega,x)}{l_{t}(\omega,x)+1}\ind_{B}(\omega,x)\ind_{\{t \leq \lambda(\omega,x)\}} < 2.
$$ 
This is enough to conclude the lemma. \ok

\

\subsection{Results based on the process $\beta^u$}

In this subsection we study the process $\beta^u$. First of all, as $\mu[\upsilon^u<\infty]=0$, $\beta^u=\frac{1}{\alpha^u}\ind_{[u,\tau^u)}$ under $\mu$.

\bl\label{betaformulas}
Fix $u\geq 0$. Let $\rho$ to be a $\mathbb{J}_{+}$-stopping time such that $u\leq \rho\leq \upsilon^u$. We have, for any
$\mathbb{J}_{+}$-stopping time
$\kappa$ with $u\leq \kappa \leq \rho$, for any $A \in
{\mathcal{J}}_{\kappa+}$,$$
\mathbb{E}_{\nu^u}[\ind_A\ind_{\{\kappa < \rho\}}] = \mathbb{E}_{\mu}[\ind_A\ind_{\{\kappa < \rho\}}
\beta_{\kappa}^{u}]
$$
Consequently, $\ind_{[0,u)}+\beta^{u}\ind_{[u,\rho)}$ is a $(\mu,{\mathbb{J}}_{+})$ supermartingale. 
\el

This lemma can be proved in the same way we have proved Lemma \ref{basicformulas}.

\bethe
Let $u\geq 0$. Let $M$ be a positive $(\mathbb{P},\mathbb{F})$ local martingale (assumed to be $\mathbb {F}^\circ_+$ adapted). Then, $M\beta^u(\phi)$ is a $(\mathbb{P},\mathbb{G})$ supermartingale on $[u,\infty)$. In particular, if, for a $\mathbb{G}$ stopping time $T\geq u$, $\beta^u_T(\phi)>0$ $\mathbb{P}$-almost surely, $M$ is a $(\mathbb{P},\mathbb{G})$ semimartingale on $(u,T]$.
\ethe

\textbf{Proof.} By stopping $M$ with $\mathbb{F}^\circ_+$ stopping time, we can assume that $M$ is in class$(D)$. We have, for $t>s>u$, for $A\in\mathcal{J}^\circ_{s+}$, 
$$
\dcb
&&\mathbb{E}_{\mathbb{P}}[\ind_A(\phi)M_{t}\beta^u(\phi)_t]
\\
&=&\mathbb{E}_{\mathbb{P}}[\ind_A(\phi)M_{t}\beta^u(\phi)_t\ind_{\{t < \upsilon^u(\phi)\}}]\ \mbox{ (see Lemma \ref{h+})}
\\
&=&\mathbb{E}_{\mu}[\ind_AM(\iota)_{t}\beta^u_t\ind_{\{t < \upsilon^u\}}]
\\
&=&
\mathbb{E}_{\nu^u}[\ind_AM(\iota)_t\ind_{\{t < \upsilon^u\}}]\\
&\leq&
\mathbb{E}_{\nu^u}[\ind_AM(\iota)_t\ind_{\{s < \upsilon^u\}}]\\
&=&
\mathbb{E}_{\nu^u}[\ind_AM(\iota)_s\ind_{\{s < \upsilon^u\}}]\ \mbox{ according to Lemma \ref{martingaletransfer}}\\
&=&\mathbb{E}_{\mu}[\ind_AM_{s}(\iota)\beta^u_s\ind_{\{s < \upsilon^u\}}]\\
&=&\mathbb{E}_{\mathbb{P}}[\ind_A(\phi)M_{s}\beta^u(\phi)_s\ind_{\{s < \upsilon^u(\phi)\}}]\\
&=&\mathbb{E}_{\mathbb{P}}[\ind_A(\phi)M_{s}\beta^u(\phi)_s]
\\
\dce
$$
This means $$
\mathbb{E}_{\mathbb{P}}[M_{t}\beta^u(\phi)_t|\mathcal{G}^\circ_{s+}]
\leq M_{s}\beta^u(\phi)_s
$$ 
Notice that the above computation shows also that $M_{t}\beta^u(\phi)_t$ is $\mathbb{P}$ integrable. Applying Lemma \ref{FF0}, we prove the first part of the theorem. 

To see the second part of the theorem, we need only to note that $\frac{1}{\beta^u(\phi)}$ is a $(\mathbb{P},\mathbb{G})$ semimartingale on $(u,T]$ according to \cite[Proposition 2.9]{choulli-stricker}. \ok

\

\section{The classical results}

In this section, we work on a probability space $(\Omega,\mathcal{B},\mathbb{P})$ satisfying the topological assumption (i.e. the existence of $\mathbb{F}^\circ$ behind the filtration $\mathbb{F}$) as it is described in the subsection \ref{h-formula}.  

One of the motivation to introduce the local solution method was to provide a unified interpretation of the classical formulas. That is what is done in this section. We employ the same notations as in the preceding section.

\subsection{Jacod's criterion}\label{jcrit}

We consider a random variable $\xi$ taking values in a polish space $E$. We consider the enlargement of filtration given by $$
\mathcal{G}_t=\cap_{s>t}(\mathcal{F}_s\vee\sigma(\xi)),\ t\geq 0
$$
It is the so-called initial enlargement of filtration. In \cite{jacod2}, the following assumption is introduced

\bassumption($A'$)
For any $t\geq 0$, the conditional law of $\xi$ given the $\sigma$-algebra $\mathcal{F}_t$ is absolutely continuous with respect to the (unconditional) law $\lambda$ of $\xi$
\eassumption
It is proved then (cf. \cite[Lemme(1.8)]{jacod2}) that, if $q=q_t(x)=q_t(\omega,x)$ denotes the density function of the conditional law of $\xi$ given the $\sigma$-algebra $\mathcal{F}_t$ with respect to $\lambda(dx)$, there exists a jointly measurable version of $q$ such that $t\rightarrow q_t(\omega,x)$ is càdlàg, and for every $x$, the map $(\omega,t)\rightarrow q_t(\omega,x)$ is a $(\mathbb{P},\mathbb{F})$ martingale. 

\bethe\label{jacod-criterion} (\cite[Théorème(2.1)]{jacod2})
For every bounded $\mathbb{F}$ martingale $M$, it is a $\mathbb{G}$ special semimartingale whose drift is given by the formula:$$
\frac{1}{q_{-}(\xi)}\centerdot \left.(\cro{q(x),M}^{\mathbb{P}\cdot\mathbb{F}})\right|_{x=\xi}
$$ 
\ethe

Let us see how the above formula can be obtained by the local solution method. We have assumed the existence of $\mathbb{F}^\circ$ behind the filtration $\mathbb{F}$. Note that, if Assumption($A'$) holds, since $\mathcal{F}^\circ_t\subset\mathcal{F}_t, t\geq 0,$ the same assumption holds with respect to the filtration $\mathbb{F}^\circ$. Let $p=p_t(x)=p_t(\omega,x)$ denote the density function of the conditional law of $\xi$ given the $\mathcal{F}^\circ_t$ with respect to $\lambda(dx)$. We also assume that $q_t(\omega,x)$ is $\mathbb{J}_+$ adapted by Lemma \ref{FF0}.

Define the process $I$ by $I_t=\xi, \forall t\geq 0$ ($I$ is a constant process taking values in the space $C(\mathbb{R}_+,E)$ of all continuous function on $\mathbb{R}_+$ taking value in $E$, equipped with the natural filtration). Define the filtration $\mathbb{G}^\circ$ as in subsection \ref{h-formula}. Then, the filtrations $\mathbb{G}$ and $\mathbb{G}^\circ$ are related in the way : $\mathcal{G}_t=\mathcal{G}^\circ_{t+}\vee\mathcal{N}, t\geq 0,$ as in subsection \ref{h-formula}.

Let $u\geq 0$. We have clearly $\pi_{u,I/I}(y'',dy)=\delta_{y''}(dy)$. For $x\in E$, denote by $\mathfrak{c}_x$ the constant function at value $x$. The process $I$ is written as $I=\mathfrak{c}_\xi$. For $t\geq u$, for any bounded borel function $F$ on $I$, we have
$$
\dcb
\int \pi_{t,I/F}(\mathsf{i},dy)F(y)

=\mathbb{E}[F(I)|\mathcal{F}^\circ_t]

=\mathbb{E}[F(\mathfrak{c}_{\xi})|\mathcal{F}^\circ_t]
=\int F(\mathfrak{c}_{x})p_t(x)\lambda(dx)
\dce
$$
Therefore,
$$
\dcb
&&\int \pi_{t,I/F}(\mathsf{i},dy)F(y)\ind_{\{p_u(y_0)>0\}}
\\
&=&\int F(\mathfrak{c}_{x})\ind_{\{p_u((\mathfrak{c}_{x})_0)>0\}}p_t(x)\lambda(dx)

\\
&=&\int F(\mathfrak{c}_{x})\frac{p_t((\mathfrak{c}_{x})_0)}{p_u((\mathfrak{c}_{x})_0)}\ind_{\{p_u((\mathfrak{c}_{x})_0)>0\}}p_u(x)\lambda(dx)\\
\\
&=&\int \pi_{u,I/F}(\mathsf{i},dy)F(y)\frac{p_t(y_0)}{p_u(y_0)}\ind_{\{p_u(y_0)>0\}}\\
\\
&=&\int \pi_{u,I/F}(\mathsf{i},dy'')\int \pi_{u,I/I}(y'',dy)F(y)\frac{p_t(y_0)}{p_u(y_0)}\ind_{\{p_u(y_0)>0\}}
\dce
$$
By Lemma \ref{alpha}, we conclude that $\tau^u(\omega,,y)\ind_{\{p_u(\omega,y_0)>0\}}=\infty\ind_{\{p_u(\omega,y_0)>0\}}$ under $(\mu+\nu^u)$ and $$
\alpha^u_t(\omega,y)\ind_{\{p_u(\omega,y_0)>0\}}=\ind_{\{p_u(\omega,y_0)>0\}}\frac{p_{t+}(\omega,y_0)}{p_u(\omega,y_0)}, \ t\geq u
$$ 
From now on, let $u=0$. By the first formula in Lemma \ref{basicformulas}, the measure
$
\ind_{\{p_0(y_0)>0\}}\frac{p_{t+}(y_0)}{p_0(y_0)}d\nu^0
$
coincides with 
$\ind_{\{p_0(\iota,\zeta_0)>0\}}d\mu$ on $\mathcal{J}_{t+}$. By the second identity in the proof of Lemma \ref{h+}, we can write, for $A\in\mathcal{F}^\circ_t, B\in\mathcal{I}_t$ 
$$
\dcb
&&\int\ind_A(\omega)\ind_B(y)\ind_{\{p_0(\omega,y_0)>0\}}\frac{q_{t}(\omega,y_0)}{p_0(\omega,y_0)}d \nu^{0}(d\omega,dy)\\
&=&\mathbb{E}_{\mathbb{P}}[\ind_A(\mathsf{i})\int \pi_{0,I/F}(\mathsf{i}, dy'') \int \pi_{0,I/I}(y'', dy)\ind_B(y)\ind_{\{p_0(\mathsf{i},y_0)>0\}}\frac{q_{t}(\mathsf{i},y_0)}{p_0(\mathsf{i},y_0)}]\\

&=&\mathbb{E}_{\mathbb{P}}[\ind_A(\mathsf{i})\int \pi_{0,I/F}(\mathsf{i}, dy) \ind_B(y)\ind_{\{p_0(\mathsf{i},y_0)>0\}}\frac{q_{t}(\mathsf{i},y_0)}{p_0(\mathsf{i},y_0)}]\\

&=&\mathbb{E}_{\mathbb{P}}[\ind_A(\mathsf{i})\int  \ind_B(\mathfrak{c}_x)\ind_{\{p_0(\mathsf{i},(\mathfrak{c}_x)_0)>0\}}\frac{q_{t}(\mathsf{i},(\mathfrak{c}_x)_0)}{p_0(\mathsf{i},(\mathfrak{c}_x)_0)}p_0(\mathsf{i},x)\lambda(dx)]\\

&=&\mathbb{E}_{\mathbb{P}}[\ind_A(\mathsf{i})\int  \ind_B(\mathfrak{c}_x)\ind_{\{p_0(\mathsf{i},x)>0\}}q_t(\mathsf{i},x)\lambda(dx)]\\

&=&\mathbb{E}_{\mathbb{P}}[\ind_A(\mathsf{i})\mathbb{E}[ \ind_B(\mathfrak{c}_\xi)\ind_{\{p_0(\mathsf{i},\xi)>0\}}|\mathcal{F}_t]]\\

&=&\mathbb{E}_{\mathbb{P}}[\ind_A(\mathsf{i}) \ind_B(I)\ind_{\{p_0(\mathsf{i},\xi)>0\}}]\\

&=&\mathbb{E}_{\mu}[\ind_A(\iota) \ind_B(\zeta)\ind_{\{p_0(\iota,\zeta_0)>0\}}]\\
\dce
$$
On the other hand, with the martingale property of $q(x)$, we check that $\ind_{\{p_0(\omega,y_0)>0\}}\frac{q_{t}(\omega,y_0)}{p_0(\omega,y_0)}, t\geq 0$, is a $(\nu^u,\mathbb{J}_{+})$ martingale. The above identity remains valid, if we replace $t$ by $t+\epsilon, \epsilon>0,$ and if we replace $\ind_A(\iota) \ind_B(\zeta)$ by a $\ind_C, C\in\mathcal{J}_{t+}$. Let then $\epsilon\downarrow 0$. We prove thus that the two measures $\ind_{\{p_0(y_0)>0\}}\frac{p_{t+}(y_0)}{p_0(y_0)}d\nu^0$ and $\ind_{\{p_0(y_0)>0\}}\frac{q_{t}(y_0)}{p_0(y_0)}d\nu^0$ coincide all with 
$\ind_{\{p_0(\iota,\zeta_0)>0\}}d\mu$ on $\mathcal{J}_{t+}$. Consequently $\ind_{\{p_0(y_0)>0\}}\frac{p_{t+}(y_0)}{p_0(y_0)}=\ind_{\{p_0(y_0)>0\}}\frac{q_{t}(y_0)}{p_0(y_0)}$ $\nu^0$-almost surely.

Note that $\rho=\infty\ind_{\{p_0(y_0)>0\}}$ is a $\mathbb{J}_+$ stopping time with $0\leq \rho\leq \tau^0$. Using the product structure of the random measure $\pi_{0,F/F}(\omega,d\omega')\otimes\pi_{0,I/I}(y,dy')$, we check that $\ind_{\{p_0(y_0)>0\}}\frac{1}{p_0(y_0)}\left(\cro{M(\iota),q(x)}^{\mathbb{P}\cdot\mathbb{F}}\right)_{x=y_0}$ is a version of $\cro{M(\iota),\alpha^0\ind_{[0,\rho)}}^{\nu^0\cdot\mathbb{J}_+}$. As $$
\mathbb{E}[\ind_{\{p_0(\xi)=0\}}]
=\mathbb{E}[\mathbb{E}[\ind_{\{p_0(\xi)=0\}}|\mathcal{F}^\circ_0]]
=\mathbb{E}[\int \ind_{\{p_0(x)=0\}}p_0(x)\lambda(dx)]
=0
$$
we have $\rho(\phi)=\infty$ $\mathbb{P}$-almost surely.

Now Theorem \ref{theorem:our-formula} is applicable. We conclude that $M^{\infty}-M^0$ is a $(\mathbb{P},\mathbb{G})$ special semimartingale whose drift is given by
$$
\ind_{(0,\rho(\phi)]}\frac{1}{\alpha^0_-(\phi)}\centerdot \left(\cro{M(\iota),\alpha^0\ind_{[0,\rho)}}^{\nu^0\cdot\mathbb{J}_+}\right)\circ\phi
=
\frac{1}{q_{-}(\xi)}\centerdot \left.(\cro{q(x),M}^{\mathbb{P}\cdot\mathbb{F}})\right|_{x=\xi}
$$
This proves Theorem \ref{jacod-criterion}.

\subsection{Progressive enlargement of filtration}\label{progressive}

Consider a positive random variable $\mathfrak{t}$ (assumed not identically null). We introduce the process$$
\dcb
I= \ind_{(\mathfrak{t},\infty)}\in\mathtt{D}^g(\mathbb{R}_+,\mathbb{R})
\dce
$$
where $\mathtt{D}^g(\mathbb{R}_+,\mathbb{R})$ is the space of all càglàd functions. We equip this space with the distance $\mathfrak{d}(x,y), x,y\in \mathtt{D}^g(\mathbb{R}_+,\mathbb{R})$, as the Skorohod distance of $x_+,y_+$ in the space $\mathtt{D}(\mathbb{R}_+,\mathbb{R})$ of càdlàg functions. Let $\mathbb{G}$ be the filtration of the $\sigma$-algebras $$
\mathcal{G}_t=\cap_{s>t}(\mathcal{F}\vee\sigma(I_u:u\leq s)),\ t\geq 0
$$
It is the so-called progressive enlargement of filtration. In \cite[Proposition(4.12)]{Jeulin80}, it is proved :

\bethe\label{avant-tau}
For all bounded $(\mathbb{P},\mathbb{F})$ martingale $M$, the stopped process $M^\mathfrak{t}$ is a special $(\mathbb{P},\mathbb{G})$ semimartingale whose drift is given by $$
\ind_{(0,\mathfrak{t}]}\frac{1}{Z_-}\centerdot (\cro{N,M}^{\mathbb{P}\cdot \mathbb{F}}+\mathsf{B}^M)
$$
where $Z$ is the Azéma supermartingale $\mathbb{E}[t<\mathfrak{t}|\mathcal{F}_t], t\geq 0$, $N$ is the martingale part of $Z$, and $\mathsf{B}^M$ is the $(\mathbb{P},\mathbb{F})$ dual project of the jump process $\ind_{\{0<\mathfrak{t}\}}\Delta_\mathfrak{t}M \ind_{[\mathfrak{t},\infty)}$.
\ethe

We call the formula in this theorem the formula of progressive enlargement of filtration. Let us look at how the above formula can be established with the local solution method. Firstly recall that the existence of the filtration $\mathbb{F}^\circ$ is assumed. Notice that this filtration plays only the role of an intermediary. We can choose suitable $\mathbb{F}^\circ$ according to the situation. By Lemma \ref{FF0}, we can assume that $\widetilde{Z}_t=\mathbb{E}[\ind_{\{t\leq\mathfrak{t}\}}|\mathcal{F}_t], t\geq 0,$ is a borel function. We can then suppose that the stopped process $\widetilde{Z}^t$ is $\mathcal{F}^\circ_t$ measurable.

Let $\mathcal{I}$ be the natural filtration on $\mathtt{D}^g(\mathbb{R}_+,\mathbb{R})$ generated by the coordinates. For $x\in \mathtt{D}^g(\mathbb{R}_+,\mathbb{R})$, let $\rho(x)=\inf\{s\geq 0: x_s\neq x_0\}$. We consider $\rho$ as a function on $\Omega\times\mathtt{D}(\mathbb{R}_+,\mathbb{R})$ in identifying $\rho(x)$ with $\rho(\zeta(\omega,x))$.

Note that, for $t\geq 0$, $I^t=0$ on $\{t\leq \mathfrak{t}\}$ so that $$
\sigma(I_s:s\leq t)\cap\{\mathfrak{t}\geq t\}
=\{\emptyset,\Omega\}\cap\{\mathfrak{t}\geq  t\}
$$
We now compute $\alpha^u$ for $u=0$. For any bounded function $f\in\mathcal{I}_t$, $f(I)\ind_{\{t\leq \mathfrak{t}\}}=f(0)\ind_{\{t\leq \mathfrak{t}\}}$. Hence, $$
\dcb
&&\int\pi_{0,I/I}(I,dx)f(x)\ind_{\{t\leq\rho(x)\}}
=\mathbb{E}[f(I)\ind_{\{t\leq \mathfrak{t}\}}|\sigma(I_0)]\\
\\
&=&f(0)\mathbb{E}[\ind_{\{t\leq \mathfrak{t}\}}|\sigma(I_0)]
=f(0)\mathbb{P}[t\leq \mathfrak{t}]
=f(0)\mathbb{P}[t\leq \mathfrak{t}](1-I_0)
=\int f(x)\delta_0(dx)\ \mathbb{P}[t\leq \mathfrak{t}](1-I_0)
\dce
$$
We have then$$
\dcb
&&\int\pi_{0,I/F}(\mathsf{i},dx'')\int\pi_{0,I/I}(x'',dx)f(x)\ind_{\{t\leq \rho(x)\}}\\
&=&\int\pi_{0,I/F}(\mathsf{i},dx'')f(0)\mathbb{P}[t\leq \mathfrak{t}](1-x''_0)\\
&=&f(0)\mathbb{P}[t\leq\mathfrak{t}]\mathbb{E}[1-I_0|\mathcal{F}^\circ_0]\\
&=&f(0)\mathbb{P}[t\leq\mathfrak{t}]\\
\dce
$$
As for $\pi_{t,I/F}$, note that, by supermartingale property, if $Z_0=0$, $Z_t=0$. Hence, we have
$$
\dcb
&&\int\pi_{t,I/F}(\mathsf{i},dx)f(x)\ind_{\{t\leq\rho(x)\}}\\
&=&\mathbb{E}[f(I)\ind_{\{t\leq\mathfrak{t}\}}|\mathcal{F}^\circ_t]\\
&=&f(0)\mathbb{E}[\ind_{\{t\leq\mathfrak{t}\}}|\mathcal{F}^\circ_t]\\
&=&f(0)\widetilde{Z}_t\\
&=&\frac{\widetilde{Z}_t}{\mathbb{P}[t\leq\mathfrak{t}]} f(0)\mathbb{P}[t<\mathfrak{t}]\\
&=&\int\pi_{0,I/F}(\mathsf{i},dx'')\int\pi_{0,I/I}(x'',dx)f(x)\ind_{\{t\leq\rho(x)\}}\frac{\widetilde{Z}_t}{\mathbb{P}[t\leq\mathfrak{t}]}\\
\dce
$$
By Lemma \ref{alpha}, we conclude that $\tau^0\geq \rho$ under $(\mu+\nu^0)$ and $$
\alpha^0_t(\omega,x)\ind_{\{t<\rho(x)\}}=\lim_{s\downarrow t}\frac{\widetilde{Z}_{s}}{\mathbb{P}[s\leq\mathfrak{t}]}\ind_{\{t<\rho(x)\}}
=\frac{Z_{t}}{\mathbb{P}[t<\mathfrak{t}]}\ind_{\{t<\rho(x)\}}, t\geq 0,
$$ 
according to \cite[p.63]{Jeulin80}. 

We now compute $\cro{M(\iota),\alpha^0\ind_{[0,\rho)}}^{\nu^0\cdot\mathbb{J}_+}$. We recall that $M$ and $Z$ are $(\nu^0,\mathbb{J}_+)$ semimartingales. Let $t_0=\inf\{s: \mathbb{P}[s<\mathfrak{t}]=0\}$. Then, $$
\int\pi_{0,I/I}(I,dx)\ind_{\{t_0<\rho(x)\}}
=\mathbb{E}[\ind_{\{t_0<\mathfrak{t}\}}|\sigma(I_0)]
=\mathbb{P}[t_0<\mathfrak{t}]=0
$$
We conclude that, under $\nu^0$, the process $v_t=\frac{1}{\mathbb{P}[t<\mathfrak{t}]}\ind_{\{t<\rho(x)\}}$ is well-defined $\mathbb{J}_+$ adapted with finite variation. By integration by parts formula, we write$$
d(Z_tv_t)=v_{t-}dZ_t+Z_{t-}dv_t+d[Z,v]_t
$$ 
Therefore, $$
d[M(\iota),Zv]_t
=v_{t-}d[M(\iota),Z]_t+Z_{t-}d[M(\iota),v]_t+d[M(\iota),[Z,v]]_t
$$
We deal successively the three terms at the right hand side of the above identity. Let $0\leq s<t$. Let $A\in\mathcal{F}^\circ_s$ and $B\in\mathcal{I}_s$. In the following computation we will consider the functions on $\Omega$ and respectively on $\mathtt{D}(\mathbb{R}_+,\mathbb{R})$ as functions on the product space $\Omega\times \mathtt{D}(\mathbb{R}_+,\mathbb{R})$. We give up the notations such as $\iota$ or $\zeta$. We will not denote all $\omega$ and $x$. Begin with the last term. $$
\dcb
&&\mathbb{E}_{\nu^0}[\ind_A\ind_B\left([M,[Z,v]]_t-[M,[Z,v]]_s\right)]\\

&=&\mathbb{E}[\int\pi_{0,I/I}(I,dx)\int\pi_{0,F/F}(\mathsf{i},d\omega)\ind_A\ind_B\left([M,[Z,v]]_t-[M,[Z,v]]_s\right)]\\

&=&\mathbb{E}[\int\pi_{0,I/I}(I,dx)\ind_B\int\pi_{0,F/F}(\mathsf{i},d\omega)\ind_A\sum_{s<u\leq t}\Delta_uM \Delta_uZ\Delta_uv]\\

&=&\mathbb{E}[\int\pi_{0,I/I}(I,dx)\ind_B\int\pi_{0,F/F}(\mathsf{i},d\omega)\ind_A\int_s^t\Delta_uvd[M,Z]_u]\\

&=&\mathbb{E}[\int\pi_{0,F/F}(\mathsf{i},d\omega)\ind_A\int_s^td[M,Z]_u\int\pi_{0,I/I}(I,dx)\ind_B\Delta_uv]\\

\dce
$$
Look at the integral $\int\pi_{0,I/I}(I,dx)\ind_B\Delta_uv$ (noting the $t\geq u>s\geq 0$). $$
\dcb
\int\pi_{0,I/I}(I,dx)\ind_B\Delta_uv
&=&\mathbb{E}[\ind_B(I)\left(\frac{1}{\mathbb{P}[u<\mathfrak{t}]}\ind_{u<\mathfrak{t}}-\frac{1}{\mathbb{P}[u\leq \mathfrak{t}]}\ind_{u\leq \mathfrak{t}}\right)|\sigma(I_0)]\\

&=&\mathbb{E}[\ind_B(0)\left(\frac{1}{\mathbb{P}[u<\mathfrak{t}]}\ind_{u<\mathfrak{t}}-\frac{1}{\mathbb{P}[u\leq \mathfrak{t}]}\ind_{u\leq \mathfrak{t}}\right)|\sigma(I_0)]\\

&=&\ind_B(0)\mathbb{E}[\frac{1}{\mathbb{P}[u<\mathfrak{t}]}\ind_{u<\mathfrak{t}}-\frac{1}{\mathbb{P}[u\leq \mathfrak{t}]}\ind_{u\leq \mathfrak{t}}|\sigma(I_0)]\\

&=&\ind_B(0)\mathbb{E}[\frac{1}{\mathbb{P}[u<\mathfrak{t}]}\ind_{u<\mathfrak{t}}-\frac{1}{\mathbb{P}[u\leq \mathfrak{t}]}\ind_{u\leq \mathfrak{t}}]\\
&=&0
\dce
$$
Consequently$$
\mathbb{E}_{\nu^0}[\ind_A\ind_B\left([M,[Z,v]]_t-[M(\iota),[Z,v]]_s\right)]
=0.
$$
Consider the second term.
$$
\dcb
&&\mathbb{E}_{\nu^0}[\ind_A\ind_B\int_s^tZ_{u-}d[M,v]_u]\\

&=&\mathbb{E}[\int\pi_{0,I/I}(I,dx)\int\pi_{0,F/F}(\mathsf{i},d\omega)\ind_A\ind_B\int_s^tZ_{u-}d[M,v]_u]\\

&=&\mathbb{E}[\int\pi_{0,I/I}(I,dx)\ind_B\int\pi_{0,F/F}(\mathsf{i},d\omega)\ind_A\sum_{s<u\leq t}Z_{u-}\Delta_uv\Delta_uM]\\

&=&\mathbb{E}[\int\pi_{0,I/I}(I,dx)\ind_B\int\pi_{0,F/F}(\mathsf{i},d\omega)\ind_A\int_s^tZ_{u-}\Delta_uvdM_u]\\
&&\mbox{ because under $\pi_{0,F/F}(\mathsf{i},d\omega)$ $v$ is a deterministic process}\\
&=&0
\ \mbox{ because $M$ is martingale under $\pi_{0,F/F}(\mathsf{i},d\omega)$}\\

\dce
$$
Consider the first term
$$
\dcb
&&\mathbb{E}_{\nu^0}[\ind_A\ind_B\int_s^tv_{u-}d[M,Z]_u]\\

&=&\mathbb{E}[\int\pi_{0,I/I}(I,dx)\int\pi_{0,F/F}(\mathsf{i},d\omega)\ind_A\ind_B\int_s^tv_{u-}d[M,Z]_u]\\

&=&\mathbb{E}[\int\pi_{0,I/I}(I,dx)\ind_B\int\pi_{0,F/F}(\mathsf{i},d\omega)\ind_A\int_s^tv_{u-}d[M,Z]_u]\\

&=&\mathbb{E}[\int\pi_{0,I/I}(I,dx)\ind_B\int\pi_{0,F/F}(\mathsf{i},d\omega)\ind_A\int_s^tv_{u-}d\cro{M,Z}^{\mathbb{P}\cdot\mathbb{F}}_u]\\
&&\mbox{ because under $\pi_{0,F/F}(\mathsf{i},d\omega)$ $v$ is a deterministic process}\\
&&\mbox{ while $M$ and $Z$ has a same behavior as under $(\mathbb{P},\mathbb{F})$}\\
\dce
$$
Combining these three terms, we obtain
$$
\dcb
&&\mathbb{E}_{\nu^0}[\ind_A\ind_B\int_s^td[M,Zv]_u]\\

&=&\mathbb{E}[\int\pi_{0,I/I}(I,dx)\ind_B\int\pi_{0,F/F}(\mathsf{i},d\omega)\ind_A\int_s^tv_{u-}d\cro{M,Z}^{\mathbb{P}\cdot\mathbb{F}}_u]\\

&=&\mathbb{E}_{\nu^0}[\ind_A\ind_B\int_s^t v_{u-} d\cro{M,N}^{\mathbb{P}\cdot\mathbb{F}}_u]\\

\dce
$$
By a usual limit argument, we conclude from this identity that $$
\cro{M,Zv}^{\nu^0\cdot\mathbb{J}_+}=v_-\centerdot\cro{M,N}^{\mathbb{P}\cdot\mathbb{F}}
$$
We can then compute 
$$
\dcb
d\cro{M,\alpha^0\ind_{[0,\rho)}}^{\nu^0\cdot\mathbb{J}_+}_t
&=&d\cro{M,Zv}^{\nu^0\cdot\mathbb{J}_+}_t\\

&=& v_{t-} d\cro{M,N}^{\mathbb{P}\cdot\mathbb{F}}_t\\

&=& \frac{1}{\mathbb{P}[t\leq \mathfrak{t}]}\ind_{\{t\leq \rho\}}d\cro{M,N}^{\mathbb{P}\cdot\mathbb{F}}_t\\

&=&\frac{Z_{t-}}{\mathbb{P}[t\leq \mathfrak{t}]} \ind_{\{t\leq \rho\}}\frac{1}{Z_{t-}}d\cro{M,N}^{\mathbb{P}\cdot\mathbb{F}}_u\\

&=&\alpha^0_{t-}\frac{1}{Z_{t-}}d\cro{M,N}^{\mathbb{P}\cdot\mathbb{F}}_u
\dce
$$
and we obtain
$$
\ind_{(0,\rho(\phi)]}\frac{1}{\alpha^0_-(\phi)}\centerdot \left(\cro{M(\iota),\alpha^0\ind_{[0,\rho)}}^{\nu^0\cdot\mathbb{J}_+}\right)\circ\phi
=\ind_{(0,\mathfrak{t}]}\frac{1}{Z_-}\centerdot \cro{M,N}^{\mathbb{P}\cdot\mathbb{F}}
$$

Theorem \ref{theorem:our-formula} is applicable. $\ind_{(0,\rho(\phi)]}\centerdot M$ is then the sum of a $(\mathbb{P},\mathbb{G})$ local martingale and the following process
$$
\dcb
&&\ind_{(0,\rho(\phi)]}\frac{1}{\alpha^0_-(\phi)}\centerdot \left(\cro{M(\iota),\alpha^0\ind_{[0,\rho)}}^{\nu^0\cdot\mathbb{J}_+	}\right)\circ\phi
+\ind_{\{0<\rho(\phi)\}}\Delta_{\rho(\phi)}M\ind_{[\rho(\phi),\infty)}\\

&=&
\ind_{(0,\mathfrak{t}]}\frac{1}{Z_-}\centerdot \cro{M,N}^{\mathbb{P}\cdot\mathbb{F}}
+\ind_{\{0<\mathfrak{t}\}}\Delta_{\mathfrak{t}}M\ind_{[\mathfrak{t},\infty)}
\dce
$$
This being done, to obtain the drift of $M$ in $\mathbb{G}$, we need only to compute the $(\mathbb{P},\mathbb{G})$ dual projection of the process$$
\ind_{\{0<\rho(\phi)\}}\Delta_{\rho(\phi)}M\ind_{[\rho(\phi),\infty)}
=\ind_{\{0<\mathfrak{t}\}}\Delta_{\mathfrak{t}}M\ind_{[\mathfrak{t},\infty)}
$$
which is given by $\ind_{(0,\mathfrak{t}]}\frac{1}{Z_-}\centerdot \mathsf{B}^M$, according to \cite[Lemme 4.]{JY2}. Theorem \ref{avant-tau} is proved.

\subsection{Honest time}\label{honest}

We consider the same situation as in the subsection \ref{progressive}. But we assume in addition that $\mathfrak{t}$ is a honest time (see \cite[Proposition (5.1)]{Jeulin80}). This means that there exists an increasing $\mathbb{F}$ adapted left continuous process $A$ such that $A_t\leq t\wedge \mathfrak{t}$ and $\mathfrak{t}=A_t$ on $\{\mathfrak{t}<t\}$. Let $\omega\in\{\mathfrak{t}\leq t\}$. Then, for any $\epsilon>0$, $\omega\in\{\mathfrak{t}< t+\epsilon\}$. So, $\mathfrak{t}(\omega)=A_{t+\epsilon}(\omega)$. Let $\epsilon\downarrow 0$. We see that $\mathfrak{t}(\omega)=A_{t+}(\omega)$. In fact, we can take$$
A_t=\sup\{0\leq s<t: \widetilde{Z}_s=1\},\ t>0,
$$
where $\widetilde{Z}={}^{\mathbb{P}\cdot\mathbb{F}-o}(\ind_{[0,\mathfrak{t}]})$.

We will now establish the formula of enlargement of filtration for honest time, using the local solution method. As in subsection \ref{progressive}, we assume that $A$ and $\widetilde{Z}$ are $\mathbb{F}^\circ$ adapted. The formula before $\mathfrak{t}$ has already been proved in Theorem \ref{avant-tau}. Let us consider the formula after $\mathfrak{t}$.

Fix $u>0$. For $t\geq u$, $$
\sigma(I_s:s\leq t)\cap\{\mathfrak{t}< t\}
=\sigma(\{\mathfrak{t}< s\}:s\leq t)\cap\{\mathfrak{t}< t\}
=\sigma(\mathfrak{t})\cap\{\mathfrak{t}< t\}
$$
So, for any bounded function $f\in\mathcal{I}_t$,$$
\dcb
&&\int\pi_{u,I/I}(I,dx)f(x)\ind_{\{\rho(x)< u\}}
=\mathbb{E}[f(I)\ind_{\{\mathfrak{t}< u\}}|\sigma(I_s:s\leq u)]\\
&=&\mathbb{E}[f(I^u)\ind_{\{\mathfrak{t}< u\}}|\sigma(I_s:s\leq  u)]\\
&=&f(I^u)\ind_{\{\mathfrak{t}< u\}}=f(I)\ind_{\{\rho(I)< u\}}\\
&=&\int f(x)\ind_{\{\rho(x)< u\}} \delta_{I}(dx)
\dce
$$
As for $\pi_{t,I/F}$, we have$$
\dcb
&&\int\pi_{t,I/F}(I,dx)f(x)\ind_{\{\rho(x)< u\}}
=\mathbb{E}[f(I)\ind_{\{\mathfrak{t}< u\}}|\mathcal{F}^\circ_t]\\
&=&\mathbb{E}[f(\ind_{[\mathfrak{t},\infty)})\ind_{\{\mathfrak{t}< u\}}|\mathcal{F}^\circ_t]\\
&=&\mathbb{E}[f(\ind_{[A_u,\infty)})\ind_{\{\mathfrak{t}< u\}}|\mathcal{F}^\circ_t]\\
&=&f(\ind_{[A_u,\infty)})\mathbb{E}[\ind_{\{\mathfrak{t}< u\}}|\mathcal{F}^\circ_t]\\
\dce
$$
Let $m^u_t=\mathbb{E}[\ind_{\{\mathfrak{t}< u\}}|\mathcal{F}^\circ_t], t\geq u$. We continue
$$
\dcb
&&\int\pi_{t,I/F}(I,dx)f(x)\ind_{\{\rho(x)< u\}}\ind_{\{m^u_u>0\}}\\
&=&f(\ind_{[A_u,\infty)})m^u_t\ind_{\{m^u_u>0\}}\\
&=&f(\ind_{[A_u,\infty)})\frac{m^u_t}{m^u_u}\ind_{\{m^u_u>0\}}m^u_u\\
&=&\int\pi_{u,I/F}(I,dx)f(x)\frac{m^u_t}{m^u_u}\ind_{\{m^u_u>0\}}\ind_{\{\rho(x)< u\}}\\
&=&\int\pi_{u,I/F}(I,dx'')\int\pi_{u,I/I}(x'',dx)f(x)\frac{m^u_t}{m^u_u}\ind_{\{m^u_u>0\}}\ind_{\{\rho(x)< u\}}\\
\dce
$$
By Lemma \ref{alpha}, we conclude that $\tau^u\ind_{\{\rho<u\}}\ind_{\{m^u_u>0\}}= \infty\ind_{\{\rho<u\}}\ind_{\{m^u_u>0\}}$ under $(\mu+\nu^u)$ and $$
\alpha^u_t(\omega,x)\ind_{\{m^u_u(\omega)>0\}}\ind_{\{\rho(x)<u\}}=\frac{m^u_{t+}(\omega)}{m^u_u(\omega)}\ind_{\{m^u_u(\omega)>0\}}\ind_{\{\rho(x)<u\}}, t\geq u
$$ 
It is straightforward to compute $\ind_{\{m^u_u(\iota)>0\}}\ind_{\{\rho<u\}}\ind_{(u,\infty)}\centerdot\cro{M(\iota),\alpha^u}^{\nu^u\cdot\mathbb{J}_+}$ which is given by$$
\ind_{\{m^u_u(\iota)>0\}}\ind_{\{\rho<u\}}\ind_{(u,\infty)}\centerdot\cro{M(\iota),\alpha^u}^{\nu^u\cdot\mathbb{J}_+}
=\ind_{\{\rho<u\}}\frac{1}{m^u_u(\iota)}\ind_{\{m^u_u(\iota)>0\}}\ind_{(u,\infty)}\centerdot\cro{M,m^u_+}^{\mathbb{P}\cdot\mathbb{F}}(\iota)
$$
Let us compute $\cro{M,m^u_+}^{\mathbb{P}\cdot\mathbb{F}}$. We begin with $m^u_+$. It is a bounded $(\mathbb{P},\mathbb{F})$ martingale. For $t\geq u$, $$
\dcb
m^u_{t+}&=&\mathbb{E}[\ind_{\{\mathfrak{t}< u\}}|\mathcal{F}_t]\\
&=&\mathbb{E}[\ind_{\{\mathfrak{t}< u\}}\ind_{\{\mathfrak{t}\leq  t\}}|\mathcal{F}_t]\\
&=&\mathbb{E}[\ind_{\{A_{t+}< u\}}\ind_{\{\mathfrak{t}\leq t\}}|\mathcal{F}_t]\\
&=&\ind_{\{A_{t+}< u\}}\mathbb{E}[\ind_{\{\mathfrak{t}\leq t\}}|\mathcal{F}_t]\\
&=&\ind_{\{A_{t+}< u\}}(1-Z_t)\\
&&\mbox{ or }\\
&=&\mathbb{E}[\ind_{\{\mathfrak{t}< u\}}\ind_{\{\mathfrak{t}<  t\}}|\mathcal{F}_t]\\
&=&\mathbb{E}[\ind_{\{A_{t}< u\}}\ind_{\{\mathfrak{t}< t\}}|\mathcal{F}_t]\\
&=&\ind_{\{A_{t}< u\}}\mathbb{E}[\ind_{\{\mathfrak{t}< t\}}|\mathcal{F}_t]\\
&=&\ind_{\{A_{t}< u\}}(1-\widetilde{Z}_t)\\

\dce
$$
Note that, if $t\in\{\widetilde{Z}=1\}$, $A_{t+\epsilon}\geq t$ for any $\epsilon>0$ so that $A_{t+}\geq t$. But $A_t\leq t$ so that $A_{t+}=t$. Let $T_u=\inf\{s\geq u: A_{s+}\geq u\}$. On the set $\{A_{u+}=u\}$, $T_u=u$. On the set $\{A_{u+}<u\}$, $T_u=\inf\{s\geq u: s\in\{\widetilde{Z}=1\}\}$, where we have $A_t=A_u$ for $u<t<T_u$, and $A_{T_u+}=T_u>u$ if $T_u<\infty$. Note that, according to [Jeulin Lemme(4.3)] $\{\widetilde{Z}=1\}$ is a closed set. It follows that $T_u\in\{\widetilde{Z}=1\}$ whenever $T_u<\infty$. Therefore, on the set $\{T_u<\infty\}$, $$
1-Z_{T_u}=\widetilde{Z}_{T_u}-Z_{T_u}={}^{\mathbb{P}\cdot\mathbb{F}-o}(\ind_{[\mathfrak{t}]})_{T_n}=\mathbb{E}[\ind_{\{\mathfrak{t}=T_u\}}|\mathcal{F}_{T_u}]
$$ 
Set $w^u_t=\ind_{\{A_{t+}< u\}}=\ind_{\{t<T_u\}}$ (For the last equality, it is evident if $t<u$. It is already explained for $t=u$. For $t>u$, the condition $A_{t+}<u$ is equivalent to say $A_{u+}<u$ and $A_{t+\epsilon}<u$ for some $\epsilon>0$, which implies $T_u\geq t+\epsilon$. Inversely, if $T_u> t$, we will have $A_{u+}<u$ so that $(u-\epsilon,u]\cap\{\widetilde{Z}=1\}=\emptyset$ for some $\epsilon>0$, because $\{\widetilde{Z}=1\}$ is a closed set. This means $A_{t+}\leq u-\epsilon$.). We compute for $t>u$$$
d(w^u(1-Z))_t=-w^u_{t-}dZ_t-(1-Z_{T_u-})\ind_{\{T_u\leq t\}}+\Delta_{T_u}Z\ind_{\{T_u\leq t\}}
=-w^u_{t-}dZ_t-(1-Z_{T_u})\ind_{\{T_u\leq t\}}.
$$
Consequently
$$
\dcb
d[M,m^u_+]_t
&=&-w^u_{t-}d[M,Z]_t-\Delta_{T_u}M(1-Z_{T_u-})\ind_{\{T_u\leq t\}}+\Delta_{T_u}M\Delta_{T_u}Z\ind_{\{T_u\leq t\}}\\

&=&-w^u_{t-}d[M,Z]_t-\Delta_{T_u}M(1-Z_{T_u})\ind_{\{T_u\leq t\}}
\dce
$$
Let $v$ denote the $(\mathbb{P},\mathbb{F})$ predictable dual projection of the increasing process $\Delta_{T_u}M(1-Z_{T_u})\ind_{\{u<T_u\leq t\}}, t>u$. Let $\widetilde{A}^{\mathfrak{t}}$ (resp. $A^{\mathfrak{t}}$) denotes the $(\mathbb{P},\mathbb{F})$ optional (resp. predictable) dual projection of $\ind_{\{0<\mathfrak{t}\}}\ind_{[\mathfrak{t},\infty)}$. For any bounded $\mathbb{F}$ predictable process $H$, we have
$$
\dcb
\mathbb{E}[\int_u^\infty H_udv_u]
&=&\mathbb{E}[H_{T_u}\Delta_{T_u}M(1-Z_{T_u})\ind_{\{u<T_u<\infty\}}]\\
&=&\mathbb{E}[H_{T_u}\Delta_{T_u}M\mathbb{E}[\ind_{\{\mathfrak{t}=T_u\}}|\mathcal{F}_{T_u}]\ind_{\{u<T_u<\infty\}}]\\
&=&\mathbb{E}[\mathbb{E}[H_{T_u}\Delta_{T_u}M\ind_{\{\mathfrak{t}=T_u\}}|\mathcal{F}_{T_u}]\ind_{\{u<T_u<\infty\}}]\\
&=&\mathbb{E}[\mathbb{E}[H_{\mathfrak{t}}\Delta_{\mathfrak{t}}M\ind_{\{\mathfrak{t}=T_u\}}|\mathcal{F}_{T_u}]\ind_{\{u<T_u<\infty\}}]\\

&=&\mathbb{E}[H_{\mathfrak{t}}\Delta_{\mathfrak{t}}M\ind_{\{\mathfrak{t}=T_u\}}\ind_{\{u<T_u<\infty\}}]\\

&=&\mathbb{E}[H_{\mathfrak{t}}\Delta_{\mathfrak{t}}M\ind_{\{\mathfrak{t}=T_u\}}\ind_{\{u<\mathfrak{t}<\infty\}}]\\

&=&\mathbb{E}[H_{\mathfrak{t}}\Delta_{\mathfrak{t}}M\ind_{\{\mathfrak{t}\leq T_u\}}\ind_{\{u<\mathfrak{t}<\infty\}}]\\
&&\mbox{ because $T_u\in\{\widetilde{Z}=1\}$ and $\mathfrak{t}$ is the last instant of the set $\{\widetilde{Z}=1\}$}\\
&=&\mathbb{E}[\int_0^\infty H_{s}\Delta_{s}M\ind_{(u,T_u]}(s)d(\ind_{[\mathfrak{t},\infty)})_s]\\

&=&\mathbb{E}[\int_0^\infty H_{s}\Delta_{s}M\ind_{(u,T_u]}(s)d\widetilde{A}^{\mathfrak{t}}_s]\\

&=&\mathbb{E}[\int_0^\infty H_{s}\Delta_{s}M\ind_{(u,T_u]}(s)d(\widetilde{A}^{\mathfrak{t}}-A^{\mathfrak{t}})_s]\\
&&\mbox{ because $M$ is a bounded martingale and $A^{\mathfrak{t}}$ is predictable}\\

&=&\mathbb{E}[\int_u^{T_u} H_{s}d[M,\widetilde{A}^{\mathfrak{t}}-A^{\mathfrak{t}}]_s]\\

&=&\mathbb{E}[\int_u^{T_u} H_{s}d\cro{M,\widetilde{A}^{\mathfrak{t}}-A^{\mathfrak{t}}}^{\mathbb{P}\cdot\mathbb{F}}_s]\\
\dce
$$
Putting these results together, we can write
$$
\dcb
$$
&&\ind_{\{m^u_u(\iota)>0\}}\ind_{\{\rho(x)<u\}}\ind_{(u,\tau^u]}\frac{1}{\alpha^u_-}\centerdot \cro{M(\iota),\alpha^u}^{\nu^u\cdot\mathbb{J}_+}\\
&=&
\ind_{\{\rho(x)<u\}}\left(\ind_{(u,\infty)}\frac{m^u_u}{m^u_-} \frac{1}{m^u_u}\ind_{\{m^u_u>0\}}\centerdot\cro{M,m^u_+}^{\mathbb{P}\cdot\mathbb{F}}\right)(\iota)\\

&=&
\ind_{\{\rho(x)<u\}}\left(\ind_{(u,\infty)} \ind_{\{m^u_u>0\}}\frac{1}{m^u_-}\centerdot\cro{M,m^u_+}^{\mathbb{P}\cdot\mathbb{F}}\right)(\iota)\\

&=&
\ind_{\{\rho(x)<u\}}\left(\ind_{(u,\infty)} \ind_{\{m^u_u>0\}}\frac{1}{m^u_-}\centerdot(-\ind_{(u,T_u]}\centerdot\cro{M,N}^{\mathbb{P}\cdot\mathbb{F}}-\ind_{(u,T_u]}\centerdot\cro{M,\widetilde{A}^{\mathfrak{t}}-A^{\mathfrak{t}}}^{\mathbb{P}\cdot\mathbb{F}})\right)(\iota)\\

&=&
-\ind_{\{\rho(x)<u\}}\left(\ind_{(u,\infty)} \ind_{\{m^u_u>0\}}\ind_{(u,T_u]}\frac{1}{m^u_-}\centerdot(\cro{M,N}^{\mathbb{P}\cdot\mathbb{F}}+\cro{M,\widetilde{A}^{\mathfrak{t}}-A^{\mathfrak{t}}}^{\mathbb{P}\cdot\mathbb{F}})\right)(\iota)\\

&=&
-\ind_{\{\rho(x)<u\}}\left( \ind_{\{m^u_u>0\}}\ind_{(u,T_u]}\frac{1}{1-Z_-}\centerdot(\cro{M,N}^{\mathbb{P}\cdot\mathbb{F}}+\cro{M,\widetilde{A}^{\mathfrak{t}}-A^{\mathfrak{t}}}^{\mathbb{P}\cdot\mathbb{F}})\right)(\iota)
\dce
$$
Pull this identity back to $\Omega$, we have
$$
\dcb
&&\ind_{\{m^u_u>0\}}\ind_{\{\rho(\phi)<u\}}\ind_{(u,\tau^u(\phi)]}\frac{1}{\alpha^u_-(\phi)}\centerdot \cro{M(\iota),\alpha^u}^{\nu^u\cdot\mathbb{J}_+}(\phi)\\

&=&-\ind_{\{\mathfrak{t}<u\}} \ind_{\{m^u_u>0\}}\ind_{(u,T_u]}\frac{1}{1-Z_-}\centerdot(\cro{M,N}^{\mathbb{P}\cdot\mathbb{F}}+\cro{M,\widetilde{A}^{\mathfrak{t}}-A^{\mathfrak{t}}}^{\mathbb{P}\cdot\mathbb{F}})\\

&=&-\ind_{\{\mathfrak{t}<u\}} \ind_{\{1-\widetilde{Z}_u>0\}}\ind_{(u,T_u]}\frac{1}{1-Z_-}\centerdot(\cro{M,N}^{\mathbb{P}\cdot\mathbb{F}}+\cro{M,\widetilde{A}^{\mathfrak{t}}-A^{\mathfrak{t}}}^{\mathbb{P}\cdot\mathbb{F}})\\

&=&-\ind_{\{\mathfrak{t}<u\}} \ind_{(u,\infty)}\frac{1}{1-Z_-}\centerdot(\cro{M,N}^{\mathbb{P}\cdot\mathbb{F}}+\cro{M,\widetilde{A}^{\mathfrak{t}}-A^{\mathfrak{t}}}^{\mathbb{P}\cdot\mathbb{F}})
\dce
$$
because on the set $\{\mathfrak{t}<u\}$, $1-\widetilde{Z}_u>0$ and $T_u=\infty$.

Now we can apply Theorem \ref{theorem:our-formula}, according to which (noting that, according to \cite[Lemme(4.3) ii) and p.63]{Jeulin80}, $1-Z_->0$ on $(\mathfrak{t},\infty)$) $\ind_{\{\rho(\phi)<u\}}\ind_{\{m^u_u>0\}}\ind_{(u,\tau^u(\phi)]}\centerdot M$ is the sum of a $(\mathbb{P},\mathbb{G})$ local martingale and the following process with finite variation: 
$$
\ind_{\{\rho(\phi)<u\}}\ind_{\{m^u_u>0\}}\left(\ind_{(u,\tau^u(\phi)]}\frac{1}{\alpha^u_-(\phi)}\centerdot \left(\cro{M(\iota),\alpha^u}^{\nu^u\cdot\mathbb{J}_+}\right)\circ\phi
+\ind_{\{u<\tau^u(\phi)\}}\Delta_{\tau^u(\phi)}M\ind_{[\tau^u(\phi),\infty)}\right)
$$
In particular, noting that, on the set $\{m^u_u>0, \rho(\phi)<u\}=\{\mathfrak{t}<u\}$, $\tau(\phi)=\infty$ so that $$
\ind_{\{0<\tau^u(\phi)\}}\Delta_{\tau^u(\phi)}M\ind_{[\tau^u(\phi),\infty)}=0,
$$
the above formula implies that $\ind_{\{\mathfrak{t}<u\}}\ind_{(u,\infty)}\centerdot M$ is a $(\mathbb{P},\mathbb{G})$ special semimartingale whose drift is given by 
$$
\ind_{\{\mathfrak{t}<u\}}\ind_{(u,\infty)}\frac{1}{\alpha^u_-(\phi)}\centerdot \left(\cro{M(\iota),\alpha^u}^{\nu^u\cdot\mathbb{J}_+}\right)\circ\phi
=
-\ind_{\{\mathfrak{t}<u\}} \ind_{(u,\infty)}\frac{1}{1-Z_-}\centerdot(\cro{M,N}^{\mathbb{P}\cdot\mathbb{F}}+\cro{M,\widetilde{A}^{\mathfrak{t}}-A^{\mathfrak{t}}}^{\mathbb{P}\cdot\mathbb{F}})
$$

For $u\in\mathbb{Q}_+^*$, set $\mathtt{B}_u=\{\mathfrak{t}<u\}\cap(u,\infty)$. The $\mathtt{B}_u$ are $\mathbb{G}$ predictable left intervals. For every $u\in\mathbb{Q}_+^*$, $\ind_{\mathtt{B}_u}\centerdot M$ is a $(\mathbb{P,\mathbb{G}})$ special semimartingale with the above given drift. The first and the second conditions in Assumption \ref{partialcovering} are satisfied. 

The random measure generated by the above drifts on $\cup_{u\in\mathbb{Q}_+^*}\mathtt{B}_u=(\mathfrak{t},\infty)$ is$$
-\ind_{(\mathfrak{t},\infty)} \frac{1}{1-Z_-}\centerdot(\cro{M,N}^{\mathbb{P}\cdot\mathbb{F}}+\cro{M,\widetilde{A}^{\mathfrak{t}}-A^{\mathfrak{t}}}^{\mathbb{P}\cdot\mathbb{F}})
$$
Knowing the ${}^{\mathbb{P}\cdot\mathbb{F}-p}(\ind_{(\mathfrak{t},\infty)})=(1-Z_-)$, it can be checked that this random measure has a distribution function. Lemma \ref{UB} is applicable. It yields that $\ind_{(\mathfrak{t},\infty)}\centerdot M=\ind_{\cup_{u\in\mathbb{Q}_+^*}\mathtt{B}_u}\centerdot M$ is a $(\mathbb{P,\mathbb{G}})$ special semimartingale with the drift
$$
-\ind_{(\mathfrak{t},\infty)} \frac{1}{1-Z_-}\centerdot(\cro{M,N}^{\mathbb{P}\cdot\mathbb{F}}+\cro{M,\widetilde{A}^{\mathfrak{t}}-A^{\mathfrak{t}}}^{\mathbb{P}\cdot\mathbb{F}})
$$
We have thus proved the formula of enlargement of filtration for the honest time $\mathfrak{t}$ (see \cite[Théorème(5.10), Lemme(5.17)]{Jeulin80}).

\bethe\label{honest}
Let $\mathbb{G}$ be the progressive enlargement of $\mathbb{F}$ with a honest time $\mathfrak{t}$. For all bounded $(\mathbb{P},\mathbb{F})$ martingale $M$, it is a special $(\mathbb{P},\mathbb{G})$ semimartingale whose drift is given by $$
\ind_{(0,\mathfrak{t}]}\frac{1}{Z_-}\centerdot (\cro{N,M}^{\mathbb{P}\cdot \mathbb{F}}+\mathsf{B}^M)
-\ind_{(\mathfrak{t},\infty)} \frac{1}{1-Z_-}\centerdot(\cro{M,N}^{\mathbb{P}\cdot\mathbb{F}}+\cro{M,\widetilde{A}^{\mathfrak{t}}-A^{\mathfrak{t}}}^{\mathbb{P}\cdot\mathbb{F}})
$$
\ethe

\subsection{Girsanov theorem and the formula of enlargement of filtration}

We have proved the classical formulas with the local solution method. Compared with the initial proofs of the classical formulas, the local solution method is made on a product space which bear more structure and provides some interesting details, especially when we compute $\alpha^u$ and $\cro{M(\iota),\alpha^u\ind_{[u,\rho)}}^{\nu^u\cdot\mathbb{J}_+}$.

The local solution method shows another interesting point : the Girsanov's theorem in the enlargement of filtration. When one look at the classical examples of enlargement of filtration, one notes immediately some resemblance between the enlargement of filtration formula and the Girsanov theorem formula. It was a question at what extend these two themes are related. Clearly, Jacod's criterion is an assumption of type of Girsanov theorem. As for the progressive enlargement of filtration, \cite{yoeurp} gives a proof of the formula using Girsanov theorem in a framework of Föllmer's measure. However, the use of Föllmer's measure makes that interpretation in \cite{yoeurp} inadequate to explain the phenomenon in the other examples such as the case of honest time. In fact, not all examples of enlargement of filtration give a formula of type of Girsanov theorem (see Section \ref{futur-infermum}). However, locally, if we can compute the function $\alpha^u$ as in Section \ref{alpha}, we will have a formula on the interval $(u,\tau^u(\phi)]$ which is of type of Girsanov theorem, because Theorem \ref{theorem:our-formula} is a generalized Girsanov theorem (cf. \cite{kunita, yoeurp}). The local computation of $\alpha^u$ is possible in the case a honest time as it is shown in subsection \ref{honest}. That is also the case in most of the classical examples.

\section{A different proof of \cite[Théorème(3.26)]{Jeulin80}}

Théorème(3.26) in \cite{Jeulin80} is a result on initial enlargement of filtration. It is interesting to note that this result was considered as an example where general methodology did not apply. We will give a proof using the local solution method. 

Let $\Omega$ be the space $C_0(\mathbb{R}_+,\mathbb{R})$ of all continuous real functions on $\mathbb{R}_+$ null at zero. The space $\Omega$ is equipped with the Wiener measure $\mathbb{P}$. Let $X$ be the identity function on $\Omega$, which is therefore a linear brownian motion issued from the origin. Set $U_t=\sup_{s\leq t}X_s$. Let $\mathbb{F}^\circ$ be the natural filtration of $X$. Define the process $I$ by $I_t=U, \forall t\geq 0$ ($I$ is a constant process taking values in the space $\Omega$). Define the filtration $\mathbb{G}^\circ$ as in Subsection \ref{h-formula}, as well as the filtrations $\mathbb{F}$ and $\mathbb{G}$. 

For $z\in\Omega$, denote by $z_t$ its coordinate at $t\geq 0$ and by $\mathfrak{s}(z)$ the function $\sup_{v\leq t}z_v, t\geq 0$. Let $\rho_u(z)=\inf\{s\geq u: \mathfrak{s}_s(z)>\mathfrak{s}_u(z)\}$. Let $\mathfrak{T}_u=\rho_u(U)$ ($\rho_u$ is considered as a map without randomness, while $\mathfrak{T}_u$ is a random variable). For $x\in C(\mathbb{R}_+,\Omega)$, $x_t$ denotes the coordinate of $x$ at $t$. Let $\mathfrak{l}$ is the identity function on $\mathbb{R}_+$. 

\bethe\label{sup-adjonction}\cite[Théorème(3.26)]{Jeulin80}
Let $M$ be a $(\mathbb{P},\mathbb{F})$ martingale. Then, $M$ is a $\mathbb{G}$ semimartingale, if and only if the $\int_0^t\frac{1}{U_s-X_s}|d\cro{M,X}^{\mathbb{P}\cdot\mathbb{F}}_s|<\infty$ for all $t>0$. In this case, the drift of $M$ is given by $$
-\frac{1}{U-X}\left(1-\frac{(U-X)^2}{\mathfrak{T}-\mathfrak{l}}\right)\centerdot\cro{M,X}^{\mathbb{P}\cdot\mathbb{F}}
$$ 
\ethe

We use the notation $\mathfrak{c}_\xi$ introduced in subsection \ref{jcrit}. Note that $I=\mathfrak{c}_{U}$ and hence $
\sigma(I)_\infty=\sigma(U)=\sigma(I)_0.
$
We have $\pi_{u,I/I}(x'',dx)=\ind_{\{x''\in E\}}\delta_{x''}(dx)$. To compute $\pi_{t,I/F}$ we introduce $B$ to be the process $B_t=X_{\mathfrak{T}_u+t}-X_{\mathfrak{T}_u}$ for $t\geq 0$. We have the decomposition, $t\geq 0$,$$
U_t=U_t\ind_{\{t<\mathfrak{T}_u\}}+U_{t-\mathfrak{T}_u+\mathfrak{T}_u}\ind_{\{\mathfrak{T}_u\leq t\}}
=U_t^{u}\ind_{\{t<\mathfrak{T}_u\}}+(\mathfrak{s}_{t-\mathfrak{T}_u}(B)+U_u)\ind_{\{\mathfrak{T}_u\leq t\}}
=U_t^{u}+\mathfrak{s}_{t-\mathfrak{T}_u}(B)\ind_{\{\mathfrak{T}_u\leq t\}}
$$
This yields$$
\sigma(U)=\sigma(U^u, \mathfrak{T}_u, \mathfrak{s}(B))
$$
We note that $\mathfrak{s}(B)$ is independent of $\mathcal{F}_{\mathfrak{T}_u}$. Therefore, for any bounded borel function $f$ on $C(\mathbb{R}_+,\Omega)$, there exists a suitable borel function $G$ such that $$
\mathbb{E}_{\mathbb{P}}[f(\mathfrak{c}_U)|\mathcal{F}_{\mathfrak{T}_u}]
=G(U^u,\mathfrak{T}_u)
$$
Now we can compute the function $h^u_t$ for $t\geq u$. 
$$
\dcb
&&\int\pi_{t,I/F}(X,dx)f(x)\ind_{t\leq\rho_u(x_0)}\\ 
&=& \mathbb{E}_{\mathbb{P}}[f(I)\ind_{t\leq \mathfrak{T}_u}|\mathcal{F}_t^\circ]\\
&=&\mathbb{E}_{\mathbb{P}}[f(\mathfrak{c}_U)\ind_{t\leq \mathfrak{T}_u}|\mathcal{F}_t^\circ]\\

&=&\mathbb{E}_{\mathbb{P}}[\mathbb{E}_{\mathbb{P}}[f(\mathfrak{c}_U)|\mathcal{F}_{\mathfrak{T}_u}]\ind_{t\leq \mathfrak{T}_u}|\mathcal{F}_t^\circ]\\

&=&\mathbb{E}_{\mathbb{P}}[G(U^u, \mathfrak{T}_u)\ind_{t\leq \mathfrak{T}_u}|\mathcal{F}_t^\circ]\\

&=&\mathbb{E}_{\mathbb{P}}[G(U^u, T_{U_u}\circ\theta_t+t)\ind_{t\leq \mathfrak{T}_u}|\mathcal{F}_t^\circ]\\ &&\mbox{where $\theta$ is the usual translation operator and $T_a=\inf\{s\geq 0: X_s=a\}$}\\

&=&\int_0^\infty G(U^u, v+t)\frac{1}{\sqrt{2\pi}}\frac{(U_u-X_t)}{\sqrt{v}^3}e^{-\frac{(U_u-X_t)^2}{2v}}dv\ind_{t\leq \mathfrak{T}_u}\\

&&\mbox{ see \cite[p.80]{KS} for the density function of $T_a$}\\
&=&\int_0^\infty G(U^u, v+t)\frac{1}{\sqrt{2\pi}}\frac{(U_t-X_t)}{\sqrt{v}^3}e^{-\frac{(U_t-X_t)^2}{2v}}dv\ind_{t\leq \mathfrak{T}_u}\\

&=&\int_t^\infty G(U^u, s)\frac{1}{\sqrt{2\pi}}\frac{(U_t-X_t)}{\sqrt{s-t}^3}e^{-\frac{(U_t-X_t)^2}{2(s-t)}}ds\ind_{t\leq \mathfrak{T}_u}\\
\dce
$$
and, for a positive borel function $g$ on $\mathbb{R}_+$,
$$
\dcb
&&\int\pi_{u,I/F}(X,dx)f(x)g(\rho_u(x_0))\\ 
&=& \mathbb{E}_{\mathbb{P}}[f(I)g(\mathfrak{T}_u)|\mathcal{F}_u^\circ]\\
&=&\mathbb{E}_{\mathbb{P}}[f(\mathfrak{c}_U)g(\mathfrak{T}_u)|\mathcal{F}_u^\circ]\\

&=&\mathbb{E}_{\mathbb{P}}[G(U^u,\mathfrak{T}_u)g(\mathfrak{T}_u)|\mathcal{F}_u^\circ]\\

&=&\mathbb{E}_{\mathbb{P}}[G(U^u,T_{U_u}\circ\theta_u+u)g(T_{U_u}\circ\theta_u+u)|\mathcal{F}_u^\circ]\\

&=&\int_0^\infty G(U^u, v+u)g(v+u)\frac{1}{\sqrt{2\pi}}\frac{(U_u-X_u)}{\sqrt{v}^3}e^{-\frac{(U_u-X_u)^2}{2v}}dv\\

&=&\int_u^\infty G(U^u, s)g(s)\frac{1}{\sqrt{2\pi}}\frac{(U_u-X_u)}{\sqrt{s-u}^3}e^{-\frac{(U_u-X_u)^2}{2(s-u)}}ds\\
\dce
$$
Set $$
\Psi(y,a)=\frac{1}{\sqrt{2\pi}}\frac{y}{\sqrt{a}^3}e^{-\frac{y^2}{2a}},\ y>0, a>0
$$
The above two identities yields$$
\dcb
&&\int\pi_{t,I/F}(X,dx)f(x)\ind_{t\leq \rho_u(x_0)}\\
&=&\int_t^\infty G(U^u, s)\Psi(U_t-X_t,s-t)ds\ind_{t\leq \mathfrak{T}_u}\\

&=&\int_u^\infty G(U^u, s)\ind_{t\leq s}\Psi(U_t-X_t,s-t)ds\ind_{t\leq \mathfrak{T}_u}\\

&=&\int_u^\infty G(U^u, s)\ind_{t\leq s}\frac{\Psi(U_t-X_t,s-t)}{\Psi(U_u-X_u,s-u)}\frac{1}{\sqrt{2\pi}}\frac{(U_u-X_u)}{\sqrt{s-u}^3}e^{-\frac{(U_u-X_u)^2}{2(s-u)}}ds\ind_{t\leq \mathfrak{T}_u}\\

&=&
\int\pi_{u,I/F}(X,dx)f(x)\ind_{t\leq \rho_u(x_0)}\frac{\Psi(U_t-X_t,\rho_u(x_0)-t)}{\Psi(U_u-X_u,\rho_u(x_0)-u)}\ind_{t\leq \mathfrak{T}_u}\\

&=&
\int\pi_{u,I/F}(X,dx'')\int\pi_{u,I/I}(x'',dx)f(x)\ind_{t\leq \rho_u(x_0)}\frac{\Psi(U_t-X_t,\rho_u(x_0)-t)}{\Psi(U_u-X_u,\rho_u(x_0)-u)}\ind_{t\leq \mathfrak{T}_u}
\dce
$$
Let us identify $X(\iota)$ with $X$. According to Lemma \ref{alpha}, $\tau^u\geq \rho_u(x_0)$ and $$
\alpha^u_t\ind_{\{t<\rho_u(x_0)\}}
=\frac{\Psi(U_t-X_t,\rho_u(x_0)-t)}{\Psi(U_u-X_u,\rho_u(x_0)-u)}\ind_{\{t<\mathfrak{T}_u\}}\ind_{\{t<\rho_u(x_0)\}}
$$
Under the probability $\nu^u$, $X$ is a $\mathbb{J}_+$ brownian motion. Note that $\rho_u(x_0)$ is in $\mathcal{J}_0$ and $\mathfrak{T}_u$ is a stopping time in a brownian filtration. So, $\rho_u(x_0)\wedge \mathfrak{T}_u$ is a $\mathbb{J}_+$ predictable stopping time. The martingale part of $\alpha^u$ on $(u,\rho_u(x_0)\wedge \mathfrak{T}_u)$ is then given by$$
\dcb
&&-\frac{1}{\Psi(U_u-X_u,\rho_u(x_0)-u)}\frac{\partial}{\partial y}\Psi(U_t-X_t,\rho_u(x_0)-t) dX_t\\
&=&
-\frac{1}{\Psi(U_u-X_u,\rho_u(x_0)-u)}\Psi(U_t-X_t,\rho_u(x_0)-t)\frac{1}{U_t-X_t}(1-\frac{(U_t-X_t)^2}{\rho_u(x_0)-t}) dX_t\\
&=&
-\alpha^u_t\frac{1}{U_t-X_t}(1-\frac{(U_t-X_t)^2}{\rho_u(x_0)-t}) dX_t
\dce
$$
The continuity of the processes makes that the $\mathbb{J}_+$ predictable bracket between $M$ and $X$ under $\nu^u$ is the same as that in $(\mathbb{P},\mathbb{F})$. Theorem \ref{theorem:our-formula} is applicable. Hence, $\ind_{(u,\mathfrak{T}_u]}\centerdot M$ is a $(\mathbb{P},\mathbb{G})$ special semimartingale with drift given by $$
- \ind_{(u,\mathfrak{T}_u]}\frac{1}{U-X}(1-\frac{(U-X)^2}{\mathfrak{T}_u-\mathfrak{l}}) \centerdot \cro{M,X}
$$
Note that the above random measure has a distribution function, according to the argument in \cite[p.54]{Jeulin80}. Note also that we can replace $\frac{(U-X)^2}{\mathfrak{T}_u-\mathfrak{l}}$ by $\frac{(U-X)^2}{\mathfrak{T}_\cdot-\mathfrak{l}}$ in the above formula, because for $t\in(u,\mathfrak{T}_u]$, $\mathfrak{T}_u=\mathfrak{T}_t$.

Now we consider the random left intervals $\mathtt{B}_u=(u,\mathfrak{T}_u]$ for $u\in\mathbb{Q}_+$. The drift of $\ind_{(u,\mathfrak{T}_u]}\centerdot M$ generate a random measure $\mathsf{d}\chi^\cup$ on $\cup_{u\in\mathbb{Q}_+}(u,\mathfrak{T}_u]$, which is clearly diffuse. Note that $\cup_{u\in\mathbb{Q}_+}(u,\mathfrak{T}_u]\setminus \{U>x\}$ is a countable set. Hence $$
d\chi^\cup_t
=\ind_{\{U_t>X_t\}}\frac{1}{U_t-X_t}(1-\frac{(U_t-X_t)^2}{\mathfrak{T}_t-t}) d\cro{M,X}_t
$$ 
According to Lemma \ref{UB}, if $\mathsf{d}\chi^\cup$ has a distribution function, $\ind_{\cup_{u\in\mathbb{Q}_+}(u,\mathfrak{T}_u]}\centerdot M$ is a $(\mathbb{P},\mathbb{G})$ special semimartingale whose drift is the distribution function of $\mathsf{d}\chi^\cup$. We note that the set $\cup_{u\in\mathbb{Q}_+}(u,\mathfrak{T}_u]$ is a $\mathbb{F}$ predictable set so that the stochastic integral $\ind_{\cup_{u\in\mathbb{Q}_+}(u,\mathfrak{T}_u]}\centerdot M$ in $(\mathbb{P},\mathbb{G})$ is the same as that in $(\mathbb{P},\mathbb{F})$. But then, since $\ind_{(\cup_{u\in\mathbb{Q}_+}(u,\mathfrak{T}_u])^c}\centerdot\cro{X,X}=\ind_{\{U-X=0\}}\centerdot\cro{X,X}=0$ (recalling Levy's theorem that $U-X$ is the absolute value of a brownian motion), $\ind_{(\cup_{u\in\mathbb{Q}_+}(u,\mathfrak{T}_u])^c}\centerdot M=0$ according to the predictable representation theorem under Wiener measure. This means that $M\ (=\ind_{\cup_{u\in\mathbb{Q}_+}(u,\mathfrak{T}_u]}\centerdot M)$ itself is a $(\mathbb{P},\mathbb{G})$ special semimartingale with drift the distribution function of $\mathsf{d}\chi^\cup$. 

Inversely, if $M$ is a $(\mathbb{P},\mathbb{G})$ special semimartingale, clearly $\mathsf{d}\chi^\cup$ has a distribution function.

To achieve the proof of Theorem \ref{sup-adjonction}, we note that, according to \cite[Lemme(3.22)]{Jeulin80} $\mathsf{d}\chi^\cup$ has a distribution function if and only if $
\int_0^t\frac{1}{U_s-X_s}\left|d\cro{M,X}^{\mathbb{P}\cdot\mathbb{F}}_s\right|<\infty, \ t\geq 0.
$
We recall that $X$ does not satisfy that condition.

\

\section{A last passage time in \cite{JYC}}

The computation of $\alpha^u$ is based on absolute continuity of $\pi_{t,I/F}$ with respect to $\int \pi_{u,I/F}(\omega,dx'')\pi_{u,I/I}(x'',dx)$. The random time studied in this section is borrowed from \cite{JYC} proposed by Emery. We will compute its $\alpha^u$. This computation will show a good example how the reference measure can be random and singular.

\subsection{The problem setting}

Let $W$ be a brownian motion. Let $\xi=\sup\{0\leq t\leq 1: W_1-2W_t=0\}$. We consider the progressive enlargement of filtration with the random time $\xi$. According to \cite{JYC}, the Azéma supermartingale of $\xi$ is given by $$
\mathbb{P}[t<\xi|\mathcal{F}_t]
=1-h(\frac{|W_t|}{\sqrt{1-t}})
$$
where $
h(y)=\sqrt{\frac{2}{\pi}}\int_0^y s^2e^{-\frac{s^2}{2}}ds.
$
This will give the formula of enlargement of filtration on the random interval $[0,\xi]$ as it is shown in Theorem \ref{avant-tau}. So, we consider here only the enlargement of filtration formula on the time interval $(\xi,\infty)$. We consider this problem for $W$.

\subsection{Computing $\tau^u,\alpha^u$ for $0\leq u<1$}

We assume that $\mathbb{F}^\circ$ is the natural filtration of $W$. We introduce $
I= \ind_{[\xi,\infty)}\in\mathtt{D}(\mathbb{R}_+,\mathbb{R}).
$
Let $\mathcal{I}$ be the natural filtration on $\mathtt{D}(\mathbb{R}_+,\mathbb{R})$. For $x\in \mathtt{D}(\mathbb{R}_+,\mathbb{R})$, let $\rho(x)=\inf\{s\geq 0: x_s=1\}$. In the following computations, we will write explicitly variable $x\in \mathtt{D}(\mathbb{R_+},\mathbb{R})$, but omit writing $\omega$. We will not write the variables of the functions $\alpha^u, \tau^u$.

Fix $0\leq u<1$. For $u\leq t< 1$, for any bounded function $f\in\mathcal{I}_t$,$$
\dcb
&&\int\pi_{t,I/I}(I,dx)f(x)\ind_{\{\rho(x)\leq t\}}
=f(I)\ind_{\{\xi\leq t\}}
=\int f(x)\ind_{\{\rho(x)\leq t\}} \delta_{I}(dx)
\dce
$$
Also,
$$
\dcb
&&\int\pi_{t,I/F}(W,dx)f(x)\ind_{\{\rho(x)\leq u\}}
=\mathbb{E}[f(I)\ind_{\{\xi\leq u\}}|\mathcal{F}^\circ_t]
=\mathbb{E}[f(\ind_{[\xi,\infty)})\ind_{\{\xi\leq u\}}|\mathcal{F}^\circ_t]\\
\dce
$$
Let us compute the conditional law of $\xi$ given $\mathcal{F}^\circ_t$. It will be based on the following relation : for $0\leq a\leq u$, 
$$
\dcb
\{\xi< a\}
=\{\inf_{a\leq s\leq 1}W_s> \frac{W_1}{2}> 0\}
\cup \{\sup_{a\leq s\leq 1}W_s< \frac{W_1}{2}< 0\}
\dce
$$
We introduce some notations : $$
\dcb
\varkappa^{(t)}_a&=&\inf_{a\leq s\leq t}W_s\\ 
\theta^{(t)}_a&=&\sup_{a\leq s\leq t}W_s\\
g_t(c,b)&=&\ind_{\{c<b\}}\ind_{\{0<b\}}\frac{2(2b-c)}{\sqrt{2\pi t}^3}\exp\{-\frac{(2b-c)^2}{2t}\}\\
h_t(e,d,b)&=&\int_e^dg_t(c,b)dc =\int_e^{b\wedge d}\ind_{\{0<b\}}\frac{2(2b-c)}{\sqrt{2\pi t}^3}\exp\{-\frac{(2b-c)^2}{2t}\}dc\\
k_{t}(z,y)&=&\int_y^zg_{t}(y,b)db=\int_y^z\frac{2(2b-y)}{\sqrt{2\pi t}^3}\exp\{-\frac{(2b-y)^2}{2t}\}\ind_{\{0<b\}}db\\
\Gamma_{1-t}&=&\sup_{t\leq s<1}(-(W_s-W_t))\\
\Upsilon_{1-t}&=&\sup_{t\leq s<1}(W_s-W_t)
\dce
$$
Now we compute
$$
\dcb
\inf_{a\leq s\leq 1}W_s=\inf_{a\leq s\leq t}W_s\wedge \inf_{t\leq s\leq 1}W_s
=\inf_{a\leq s\leq t}W_s\wedge (W_t+\inf_{t\leq s\leq 1}(W_s-W_t))\\
=\inf_{a\leq s\leq t}W_s\wedge (W_t-\sup_{t\leq s\leq 1}(-(W_s-W_t)))
=\varkappa^{(t)}_a\wedge (W_t-\Gamma_{1-t})
\dce
$$
and similarly
$$
\sup_{a\leq s\leq 1}W_s=\sup_{a\leq s\leq t}W_s\vee \sup_{t\leq s\leq 1}W_s
=\sup_{a\leq s\leq t}W_s\vee (W_t+\sup_{t\leq s\leq 1}(W_s-W_t))
=\theta^{(t)}_a\vee (W_t+\Upsilon_{1-t})
$$
Hence,
$$
\dcb
&&\mathbb{E}[\ind_{\{\xi< a\}}|\mathcal{F}^\circ_t]\\
&=&\mathbb{P}[\varkappa^{(t)}_a\wedge (W_t-\Gamma_{1-t})> \frac{W_1-W_t+W_t}{2}> 0|\mathcal{F}^\circ_t]\\
&&+\mathbb{P}[\theta^{(t)}_a\vee (W_t+\Upsilon_{1-t})< \frac{W_1-W_t+W_t}{2}< 0|\mathcal{F}^\circ_t]\\

&=&\mathbb{P}[\varkappa^{(t)}_a\wedge (W_t-\Gamma_{1-t})> \frac{-(-(W_1-W_t))+W_t}{2}> 0|\mathcal{F}^\circ_t]\\
&&+\mathbb{P}[\theta^{(t)}_a\vee (W_t+\Upsilon_{1-t})< \frac{W_1-W_t+W_t}{2}< 0|\mathcal{F}^\circ_t]\\

&=&\int_{-\infty}^\infty dc\int_0^\infty db\ g_{1-t}(c,b)
\left[\ind_{\{\varkappa^{(t)}_a\wedge (W_t-b)> \frac{-c+W_t}{2}> 0\}}
+
\ind_{\{\theta^{(t)}_a\vee (W_t+b)< \frac{c+W_t}{2}< 0\}}\right]\\
&&\mbox{ cf. \cite[Proposition 8.1]{KS}}\\

&=&\int_{-\infty}^\infty dc\int_0^\infty db\ g_{1-t}(c,b)
\left[\ind_{\{W_t-2\varkappa^{(t)}_a\wedge (W_t-b)< c< W_t\}}
+
\ind_{\{-W_t+2\theta^{(t)}_a\vee (W_t+b)< c< -W_t\}}\right]\\
&=&\int_0^\infty db\ h_{1-t}(W_t-2\varkappa^{(t)}_a\wedge (W_t-b), W_t, b)
+\int_0^\infty db\ h_{1-t}(2\theta^{(t)}_a\vee (W_t+b)-W_t, -W_t, b)
\dce
$$
(Here the operation $\wedge$ is prior to the multiplication.) Note that the last expression is a continuous function in $a$. It is also equal to $\mathbb{E}[\ind_{\{\xi\leq a\}}|\mathcal{F}^\circ_t]$. $1-t>0$ being fixed, we can compute the differential with respect to $a\leq u$ under the integral. Note also $$
\frac{\partial}{\partial e}h_{1-t}(e,d,b)=0,\ \forall e>b\wedge d
$$
We have
$$
\dcb
&&d_a\mathbb{E}[\ind_{\{\xi\leq a\}}|\mathcal{F}^\circ_t]\\
&=&\int_0^\infty g_{1-t}(W_t-2\varkappa^{(t)}_a\wedge (W_t-b),b)
\ind_{\{W_t-2\varkappa^{(t)}_a\wedge (W_t-b)\leq b\wedge W_t\}}\ind_{\{\varkappa^{(t)}_a\leq W_t-b\}}db\ 2 d_a\varkappa^{(t)}_a\\
&&+
\int_0^\infty g_{1-t}(2\theta^{(t)}_a\vee (W_t+b)-W_t,b)
\ind_{\{2\theta^{(t)}_a\vee (W_t+b)-W_t\leq b\wedge(-W_t)\}}\ind_{\{b+W_t\leq \theta^{(t)}_a\}}db\ 2 (-d_a\theta^{(t)}_a)\\

&=&\ind_{\{\varkappa^{(t)}_a\geq 0\}}\int_{W_t-2\varkappa^{(t)}_a}^{W_t-\varkappa^{(t)}_a} g_{1-t}(W_t-2\varkappa^{(t)}_a,b)
db\ 2 d_a\varkappa^{(t)}_a\\
&&+
\ind_{\{\theta^{(t)}_a\leq 0\}}\int_{2\theta^{(t)}_a-W_t}^{\theta^{(t)}_a-W_t} g_{1-t}(2\theta^{(t)}_a-W_t,b)
db\ 2 (-d_a\theta^{(t)}_a)\\

&=&\ind_{\{\varkappa^{(t)}_a\geq 0\}}k_{1-t}(W_t-\varkappa^{(t)}_a,W_t-2\varkappa^{(t)}_a)\ 2 d_a\varkappa^{(t)}_a\\
&&+
\ind_{\{\theta^{(t)}_a\leq 0\}}k_{1-t}(\theta^{(t)}_a-W_t,2\theta^{(t)}_a-W_t)
db\ 2 (-d_a\theta^{(t)}_a)
\dce
$$
This yields
$$
\dcb
&&\mathbb{E}[f(\ind_{[\xi,\infty)})\ind_{\{\xi\leq u\}}|\mathcal{F}^\circ_t]\\

&=&\int_0^uf(\ind_{[a,\infty)})\ind_{\{\varkappa^{(t)}_a\geq 0\}}k_{1-t}(W_t-\varkappa^{(t)}_a,W_t-2\varkappa^{(t)}_a)\ 2 d_a\varkappa^{(t)}_a\\
&&+
\int_0^uf(\ind_{[a,\infty)}) \ind_{\{\theta^{(t)}_a\leq 0\}}k_{1-t}(\theta^{(t)}_a-W_t,2\theta^{(t)}_a-W_t)\ 2 (-d_a\theta^{(t)}_a)
\dce
$$
For $a< u\leq t<1$,$$
\dcb
\varkappa^{(t)}_a=\inf_{a\leq s\leq t}W_s=\varkappa^{(u)}_a\wedge \varkappa^{(t)}_u,\ \
\theta^{(t)}_a=\sup_{a\leq s\leq t}W_s=\theta^{(u)}_a\vee \theta^{(t)}_u\\
\dce
$$
So,
$$
\dcb
\ind_{\{a\leq u\}}d_a\varkappa^{(t)}_a
=\ind_{\{a\leq u\}}\ind_{\{\varkappa^{(u)}_a\leq \varkappa^{(t)}_u\}}d_a\varkappa^{(u)}_a,\ \

\ind_{\{a\leq u\}}d_a\theta^{(t)}_a
=\ind_{\{a\leq u\}}\ind_{\{\theta^{(u)}_a\geq \theta^{(t)}_u\}}d_a\theta^{(u)}_a
\dce
$$
Notice also $$
\dcb
d_a\varkappa^{(u)}_a=\ind_{\{W_a=\varkappa^{(u)}_a\}}d_a\varkappa^{(u)}_a,\ \
d_a\theta^{(u)}_a=\ind_{\{W_a=\theta^{(u)}_a\}}d_a\theta^{(u)}_a
\dce
$$
Let $v\in[u,1)$. We continue the computation for $u\leq t$ : 
$$
\dcb
&&\int\pi_{t,I/F}(W,dx)f(x)\ind_{\{\rho(x)\leq u\}}\ind_{\{t\leq v\}}=\mathbb{E}[f(\ind_{[\xi,\infty)})\ind_{\{\xi\leq u\}}\ind_{\{t\leq v\}}|\mathcal{F}^\circ_t]\\
	
&=&\int_0^uf(\ind_{[a,\infty)})\ind_{\{\varkappa^{(t)}_a\geq 0\}}k_{1-t}(W_t-\varkappa^{(t)}_a,W_t-2\varkappa^{(t)}_a)\ 2 d_a\varkappa^{(t)}_a\ind_{\{t\leq v\}}\\
&&+
\int_0^uf(\ind_{[a,\infty)})\ind_{\{\theta^{(t)}_a\leq 0\}} k_{1-t}(\theta^{(t)}_a-W_t,2\theta^{(t)}_a-W_t)\ 2 (-d_a\theta^{(t)}_a)\ind_{\{t\leq v\}}\\
\\
&=&\int_0^uf(\ind_{[a,\infty)})k_{1-t}(W_t-\varkappa^{(t)}_a,W_t-2\varkappa^{(t)}_a)\ind_{\{\varkappa^{(t)}_a\geq 0,\varkappa^{(u)}_a\leq \varkappa^{(t)}_u,W_a=\varkappa^{(u)}_a\}}\ 2 d_a\varkappa^{(u)}_a\ind_{\{t\leq v\}}\\
&&+
\int_0^uf(\ind_{[a,\infty)}) k_{1-t}(\theta^{(t)}_a-W_t,2\theta^{(t)}_a-W_t)\  \ind_{\{\theta^{(t)}_a\leq 0,\theta^{(u)}_a\geq \theta^{(t)}_u,W_a=\theta^{(u)}_a\}}2(-d_a\theta^{(u)}_a)\ind_{\{t\leq v\}}\\
\\
&=&\int_0^uf(\ind_{[a,\infty)})k_{1-t}(W_t-\varkappa^{(t)}_a,W_t-2\varkappa^{(t)}_a)\ind_{\{0\leq \varkappa^{(u)}_a\leq \varkappa^{(t)}_u,W_a=\varkappa^{(u)}_a\}}\ 2 d_a\varkappa^{(u)}_a\ind_{\{t\leq v\}}\\
&&+
\int_0^uf(\ind_{[a,\infty)})k_{1-t}(\theta^{(t)}_a-W_t,2\theta^{(t)}_a-W_t)\  \ind_{\{0\geq \theta^{(u)}_a\geq \theta^{(t)}_u,W_a=\theta^{(u)}_a\}}2(-d_a\theta^{(u)}_a)\ind_{\{t\leq v\}}\\
\\
&=&\int_0^uf(\ind_{[a,\infty)})\\
&&\left[k_{1-t}(W_t-\varkappa^{(t)}_a,W_t-2\varkappa^{(t)}_a)\ind_{\{0\leq \varkappa^{(u)}_a\leq \varkappa^{(t)}_u,W_a=\varkappa^{(u)}_a\}}
+k_{1-t}(\theta^{(t)}_a-W_t,2\theta^{(t)}_a-W_t)\  \ind_{\{0\geq \theta^{(u)}_a\geq \theta^{(t)}_u,W_a=\theta^{(u)}_a\}}\right]\\ 
&&(2 d_a\varkappa^{(u)}_a+2(-d_a\theta^{(u)}_a))\ind_{\{t\leq v\}}\\
\\
&=&\int_0^uf(\ind_{[a,\infty)})\\
&&\left[\frac{k_{1-t}(W_t-\varkappa^{(t)}_a,W_t-2\varkappa^{(t)}_a)}{k_{1-u}(W_u-\varkappa^{(u)}_a,W_u-2\varkappa^{(u)}_a)}\ind_{\{0\leq \varkappa^{(u)}_a\leq \varkappa^{(t)}_u,W_a=\varkappa^{(u)}_a\}}+
\frac{k_{1-t}(\theta^{(t)}_a-W_t,2\theta^{(t)}_a-W_t)}{k_{1-u}(\theta^{(u)}_a-W_u,2\theta^{(u)}_a-W_u)}  \ind_{\{0\geq \theta^{(u)}_a\geq \theta^{(t)}_u,W_a=\theta^{(u)}_a\}}\right]\\
&&\left[k_{1-u}(W_u-\varkappa^{(u)}_a,W_u-2\varkappa^{(u)}_a)\ind_{\{0\leq \varkappa^{(u)}_a,W_a=\varkappa^{(u)}_a\}}
+
k_{1-u}(\theta^{(u)}_a-W_u,2\theta^{(u)}_a-W_u) \ind_{\{0\geq \theta^{(u)}_a,W_a=\theta^{(u)}_a\}}\right]\\ 
&&(2 d_a\varkappa^{(u)}_a+2(-d_a\theta^{(u)}_a))\ind_{\{t\leq v\}}\\
\\
&=&\int\pi_{u,I/F}(W,dx)f(x)\ind_{\{\rho(x)\leq u\}}\Xi(u,t,W^t,\rho(x))\ind_{\{t\leq v\}}\\

&=&\int\pi_{u,I/F}(W,dx'')\int\pi_{t,I/I}(x'',dx)f(x)\ind_{\{\rho(x)\leq u\}}\ind_{\{t\leq v\}}\Xi(u,t,W^t,\rho(x))
\dce
$$
where $\Xi(u,t,W^t,a)$ is the function:$$
\dcb
\left[\frac{k_{1-t}(W_t-\varkappa^{(t)}_a,W_t-2\varkappa^{(t)}_a)}{k_{1-u}(W_u-\varkappa^{(u)}_a,W_u-2\varkappa^{(u)}_a)}\ind_{\{0\leq \varkappa^{(u)}_a\leq \varkappa^{(t)}_u,W_a=\varkappa^{(u)}_a\}}
+\frac{k_{1-t}(\theta^{(t)}_a-W_t,2\theta^{(t)}_a-W_t)}{k_{1-u}(\theta^{(u)}_a-W_u,2\theta^{(u)}_a-W_u)}  \ind_{\{0\geq \theta^{(u)}_a\geq \theta^{(t)}_u,W_a=\theta^{(u)}_a\}}\right]
\dce
$$
Applying Lemma \ref{alpha}, $\tau^u\ind_{\{\rho(x)\leq u\}}\geq v\ind_{\{\rho(x)\leq u\}}$ and, for $u\leq t<v$,
$$
\dcb
&&\alpha^u_t\ind_{\{\rho(x)\leq u\}}\ind_{\{0\leq \varkappa^{(u)}_{\rho(x)}\leq \varkappa^{(t)}_u,W_{\rho(x)}=\varkappa^{(u)}_{\rho(x)}\}}\\

&=&\frac{k_{1-t}(W_t-\varkappa^{(t)}_{\rho(x)},W_t-2\varkappa^{(t)}_{\rho(x)})}{k_{1-u}(W_u-\varkappa^{(u)}_{\rho(x)},W_u-2\varkappa^{(u)}_{\rho(x)})}\ind_{\{\rho(x)\leq u\}}\ind_{\{0\leq \varkappa^{(u)}_{\rho(x)}\leq \varkappa^{(t)}_u,W_{\rho(x)}=\varkappa^{(u)}_{\rho(x)}\}}
\dce
$$
and
$$
\dcb
&&\alpha^u_t\ind_{\{\rho(x)\leq u\}}\ind_{\{0\geq \theta^{(u)}_{\rho(x)}\geq \theta^{(t)}_u,W_{\rho(x)}=\theta^{(u)}_{\rho(x)}\}}\\

&=&\frac{k_{1-t}(\theta^{(t)}_{\rho(x)}-W_t,2\theta^{(t)}_{\rho(x)}-W_t)}{k_{1-u}(\theta^{(u)}_{\rho(x)}-W_u,2\theta^{(u)}_{\rho(x)}-W_u)}\ind_{\{\rho(x)\leq u\}}\ind_{\{0\geq \theta^{(u)}_{\rho(x)}\geq \theta^{(t)}_u,W_{\rho(x)}=\theta^{(u)}_{\rho(x)}\}}
\dce
$$

\subsection{The brackets computations}

Let us consider the semimartingale decomposition of the last expressions, which are composed of a part determined by the function $k_{1-t}$ and another part with indicator functions. Recall, for $z>0$,
$$
\dcb
k_{1-t}(z,y)
&=&\int_y^zg_{1-t}(y,b)db=\int_y^z\frac{2(2b-y)}{\sqrt{2\pi (1-t)}^3}\exp\{-\frac{(2b-y)^2}{2(1-t)}\}\ind_{\{0<b\}}db\\
&=&\frac{1}{\sqrt{2\pi (1-t)}^3}(1-t)\int_{y\vee 0}^z\frac{2(2b-y)}{(1-t)}\exp\{-\frac{(2b-y)^2}{2(1-t)}\}db\\

&=&\frac{1}{\sqrt{2\pi (1-t)}^3}(1-t)\int_{y\vee 0}^z\exp\{-\frac{(2b-y)^2}{2(1-t)}\}d\left(\frac{(2b-y)^2}{2(1-t)}\right)\\

&=&\frac{1}{\sqrt{2\pi (1-t)}^3}(1-t)\left.(-\exp\{-\frac{(2b-y)^2}{2(1-t)}\})\right|_{y\vee 0}^z\\
&=&\frac{1}{\sqrt{2\pi (1-t)}^3}(1-t)\left(-\exp\{-\frac{(2z-y)^2}{2(1-t)}\}+\exp\{-\frac{y^2}{2(1-t)}\}\right)
\dce
$$
Hence
$$
\dcb
&&k_{1-t}(W_t-\varkappa^{(t)}_{\rho(x)},W_t-2\varkappa^{(t)}_{\rho(x)})\\
&=&\frac{1}{\sqrt{2\pi (1-t)}^3}(1-t)\left(-\exp\{-\frac{W_t^2}{2(1-t)}\}+\exp\{-\frac{(W_t-2\varkappa^{(t)}_{\rho(x)})^2}{2(1-t)}\}\right)\\
\dce
$$
and
$$
\dcb
&&k_{1-t}(\theta^{(t)}_{\rho(x)}-W_t,2\theta^{(t)}_{\rho(x)}-W_t)\\
&=&\frac{1}{\sqrt{2\pi (1-t)}^3}(1-t)\left(-\exp\{-\frac{W_t^2}{2(1-t)}\}+\exp\{-\frac{(2\theta^{(t)}_{\rho(x)}-W_t)^2}{2(1-t)}\}\right)\\
\dce
$$
As $W=W(\iota)$ is $(\nu^u,\mathbb{J}_+)$ brownian motion, for the above processes, their $(\nu^u,\mathbb{J}_+)$ martingale part on $\{\rho\leq u\}\cap(u,1)$ are given by :
$$
\dcb
\frac{1}{\sqrt{2\pi (1-t)}^3}(1-t)\left(\frac{W_t}{1-t}\exp\{-\frac{W_t^2}{2(1-t)}\}-\frac{(W_t-2\varkappa^{(t)}_{\rho(x)})}{1-t}\exp\{-\frac{(W_t-2\varkappa^{(t)}_{\rho(x)})^2}{2(1-t)}\}\right)dW_t\\
\\
\frac{1}{\sqrt{2\pi (1-t)}^3}(1-t)\left(\frac{W_t}{1-t}\exp\{-\frac{W_t^2}{2(1-t)}\}-\frac{(W_t-2\theta^{(t)}_{\rho(x)})}{1-t}\exp\{-\frac{(W_t-2\theta^{(t)}_{\rho(x)})^2}{2(1-t)}\}\right)dW_t
\dce
$$

Consider next the indicator functions. Because $t\rightarrow \varkappa^{(t)}_u$ is decreasing, the random set $$
\{t\in[u,v]: \rho(x)\leq u, W_{\rho(x)}=\varkappa^{(u)}_{\rho(x)}\geq 0, \varkappa^{(u)}_{\rho(x)}\leq \varkappa^{(t)}_u\}
$$ 
is a random interval. Let $\lambda_\kappa$ be its supremum ($\sup\emptyset=-\infty$ by definition). Likely, the set$$
\{t\in[u,v]: \rho(x)\leq u, W_{\rho(x)}=\theta^{(u)}_{\rho(x)}\leq 0, \theta^{(u)}_{\rho(x)}\geq \theta^{(t)}_u\}
$$
is a random interval. Let $\lambda_\theta$ be its supremum. Since$$
\dcb
\{\rho(x)\leq u,W_{\rho(x)}=\varkappa^{(u)}_{\rho(x)}\geq 0\}\in\mathcal{J}_u\\
\{\rho(x)\leq u,W_{\rho(x)}=\theta^{(u)}_{\rho(x)}\leq 0\}\in\mathcal{J}_u
\dce
$$ 
$\lambda_\kappa\vee u$ and $\lambda_\theta\vee u$ are $\mathbb{J}_+$ stopping times inferior to $v\leq \tau^u$. 

Knowing the decomposition of $\alpha^u\ind_{[u,\lambda_\kappa\vee u)}$ and of $\alpha^u\ind_{[u,\lambda_\theta\vee u)}$, the continuity of the brownian motion gives a easy computation of the brackets :$$
\dcb
&&\ind_{\{u<t\leq \lambda_\kappa \}}\frac{1}{\alpha^u_{t-}}d\cro{W,\alpha^u\ind_{[u,\lambda_\kappa\vee u)}}_t
=\ind_{\{u<t\leq \lambda_\kappa\}}\frac{\left(\frac{W_t}{1-t}\exp\{-\frac{W_t^2}{2(1-t)}\}-\frac{(W_t-2\varkappa^{(t)}_{\rho(x)})}{1-t}\exp\{-\frac{(W_t-2\varkappa^{(t)}_{\rho(x)})^2}{2(1-t)}\}\right)}{\left(-\exp\{-\frac{W_t^2}{2(1-t)}\}+\exp\{-\frac{(W_t-2\varkappa^{(t)}_{\rho(x)})^2}{2(1-t)}\}\right)} dt\\
\\
&&\ind_{\{u<t\leq \lambda_\theta\}}\frac{1}{\alpha^u_{t-}}d\cro{W,\alpha^u\ind_{[u,\lambda_\theta\vee u)}}_t
=\ind_{\{u<t\leq \lambda_\theta\}}\frac{\left(\frac{W_t}{1-t}\exp\{-\frac{W_t^2}{2(1-t)}\}-\frac{(W_t-2\theta^{(t)}_{\rho(x)})}{1-t}\exp\{-\frac{(W_t-2\theta^{(t)}_{\rho(x)})^2}{2(1-t)}\}\right)}{\left(-\exp\{-\frac{W_t^2}{2(1-t)}\}+\exp\{-\frac{(2\theta^{(t)}_{\rho(x)}-W_t)^2}{2(1-t)}\}\right)} dt
\dce
$$

\subsection{The local solutions before time 1 with their drifts}

Pull the above quantities back to the initial space. We have
$$
\dcb
&&(\ind_{\{u<t\leq \lambda_\kappa\}})\circ\phi
=(\ind_{\{\rho(x)\leq u,W_{\rho(x)}=\varkappa^{(u)}_{\rho(x)}\geq 0\}}\ind_{\{u<t\leq v, \varkappa^{(u)}_{\xi}\leq \varkappa^{(t)}_u\}})\circ\phi\\
&=&\ind_{\{\xi\leq u,W_{\xi}=\varkappa^{(u)}_{\xi}\geq 0\}}\ind_{\{u<t\leq v, \varkappa^{(u)}_{\xi}\leq \varkappa^{(t)}_u\}}
=\ind_{\{\xi\leq u,W_{\xi}\geq 0\}}\ind_{\{u<t\leq v\}}
=\ind_{\{\xi\leq u,W_{1}> 0\}}\ind_{\{u<t\leq v\}}
\dce
$$
and
$$
\dcb
&&(\ind_{\{u<t\leq \lambda_\theta\}})\circ\phi
=(\ind_{\{\rho(x)\leq u,W_{\rho(x)}=\theta^{(u)}_{\rho(x)}\leq 0\}}\ind_{\{ u<t\leq v, \theta^{(u)}_{\xi}\geq \theta^{(t)}_u\}})\circ\phi\\
&=&\ind_{\{\xi\leq u,W_{\xi}=\theta^{(u)}_{\xi}\leq 0\}}\ind_{\{ u<t\leq v, \theta^{(u)}_{\xi}\geq \theta^{(t)}_u\}}
=\ind_{\{\xi\leq u,W_{\xi}\geq 0\}}\ind_{\{u<t\leq v\}}
=\ind_{\{\xi\leq u,W_{1}< 0\}}\ind_{\{u<t\leq v\}}
\dce
$$
and
$
\varkappa^{(t)}_\xi=\frac{W_1}{2},\ \theta^{(t)}_\xi=\frac{W_1}{2}.
$

$W$ is a continuous function, hence is $\mathbb{F}^\circ$-locally bounded. We can now apply Theorem \ref{theorem:our-formula} to conclude that the $\ind_{(u,\lambda_\kappa(\phi)]}\centerdot W$ and $\ind_{(u,\lambda_\theta(\phi)]}\centerdot W$ ($0\leq u\leq v<1$) are $(\mathbb{P},\mathbb{G})$ special semimartingales with the following decompositions $$
\dcb
&&\ind_{\{\xi\leq u\}}\ind_{\{W_1> 0\}}\ind_{\{u<t\leq v\}}dW_t\\
&=& d_t(\mbox{$\mathbb{G}$ local martingale})+
\ind_{\{\xi\leq u\}}\ind_{\{W_1> 0\}}\ind_{\{u<t\leq v\}}\frac{\left(\frac{W_t}{1-t}\exp\{-\frac{W_t^2}{2(1-t)}\}-\frac{(W_t-W_1)}{1-t}\exp\{-\frac{(W_t-W_1)^2}{2(1-t)}\}\right)}{\left(-\exp\{-\frac{W_t^2}{2(1-t)}\}+\exp\{-\frac{(W_t-W_1)^2}{2(1-t)}\}\right)} dt\\
\\&&\ind_{\{\xi\leq u\}}\ind_{\{W_1< 0\}}\ind_{\{u<t\leq v\}}dW_t\\
&=& d_t(\mbox{$\mathbb{G}$ local martingale})+\ind_{\{\xi\leq u\}}\ind_{\{W_1< 0\}}\ind_{\{u<t\leq v\}} \frac{\left(\frac{W_t}{1-t}\exp\{-\frac{W_t^2}{2(1-t)}\}-\frac{(W_t-W_1)}{1-t}\exp\{-\frac{(W_t-W_1)^2}{2(1-t)}\}\right)}{\left(-\exp\{-\frac{W_t^2}{2(1-t)}\}+\exp\{-\frac{(W_t-W_1)^2}{2(1-t)}\}\right)} dt
\dce
$$
From these decompositions, we obtain also the decomposition formula for $\ind_{\{\xi\leq u\}}\ind_{(u,v]}\centerdot W$

\subsection{The random measure generated by the drifts of the local solutions and its distribution function}

Let us consider $\alpha^u$ for $u\geq 1$. In this case, $I\in\mathcal{F}_u$. We obtain $\tau^u=\infty$ and $\alpha^u_t=1, t\geq u$. Clearly, $\ind_{(u,\infty)}\centerdot M$ is a $(\mathbb{P},\mathbb{G})$ martingale with a null drift.

Let $\mathtt{B}_i, i\geq 1,$ denote the left random intervals $\{\xi\leq u\}\cap(u,v], u,v\in\mathbb{Q}_+\cap[0,1)$. Let $\mathtt{B}_{0}$ denote the interval $(1,\infty)$. Then, $\cup_{i\geq 0}\mathtt{B}_i=(\xi,1)\cup(1,\infty)$. The family of $\mathtt{B}_i$ satisfies the second condition of Assumption \ref{partialcovering}. 

The drift of $W$ on the left intervals $\mathtt{B}_i$ yields a random measure concentrated on $(\xi,1)$:$$
d\chi^U_t
=\ind_{\{\xi<t<1\}}\frac{\left(\frac{W_t}{1-t}\exp\{-\frac{W_t^2}{2(1-t)}\}-\frac{(W_t-W_1)}{1-t}\exp\{-\frac{(W_t-W_1)^2}{2(1-t)}\}\right)}{\left(-\exp\{-\frac{W_t^2}{2(1-t)}\}+\exp\{-\frac{(W_t-W_1)^2}{2(1-t)}\}\right)} dt
$$ 
Since $W$ is continuous, $W-W_\xi=\ind_{(\xi,1)}\centerdot W+\ind_{(1,\infty)}\centerdot W$. According to Lemma \ref{UB}, $W-W_\xi$ is a $(\mathbb{P},\mathbb{G})$ semimartingale, whenever $\mathsf{d}\chi^\cup$ has a distribution function.

We consider then the question of distribution function. We begin with three inequalities : for $\xi\leq u\leq t<1$,
$$
\dcb
&&\left(-\exp\{-\frac{W_t^2}{2(1-t)}\}+\exp\{-\frac{(W_t-W_1)^2}{2(1-t)}\}\right)\\
&=&\exp\{-\frac{W_t^2}{2(1-t)}\}\left(\exp\{-\frac{(W_t-W_1)^2-W_t^2}{2(1-t)}\}-1\right)\\
&=&\exp\{-\frac{W_t^2}{2(1-t)}\}\left(\exp\{-\frac{-2W_tW_1+W_1^2}{2(1-t)}\}-1\right)\\
&=&\exp\{-\frac{W_t^2}{2(1-t)}\}\left(\exp\{\frac{2(W_t-W_\xi)W_1}{2(1-t)}\}-1\right)\\
&\geq&\exp\{-\frac{W_t^2}{2(1-t)}\}\frac{2(W_t-W_\xi)W_1}{2(1-t)}>0\\
\dce
$$
and
$$
\dcb
&&k_{1-t}(W_t-\varkappa^{(t)}_{a},W_t-2\varkappa^{(t)}_{a})\\
&=&\frac{1}{\sqrt{2\pi (1-t)}^3}(1-t)\left(-\exp\{-\frac{W_t^2}{2(1-t)}\}+\exp\{-\frac{(W_t-2\varkappa^{(t)}_{a})^2}{2(1-t)}\}\right)\\
&\leq &\frac{1}{\sqrt{2\pi (1-t)}^3}(1-t)\left(-\frac{(W_t-2\varkappa^{(t)}_{a})^2}{2(1-t)}+\frac{W_t^2}{2(1-t)}\right)\\
&= &\frac{1}{\sqrt{2\pi (1-t)}^3}(1-t)\frac{4W_t\varkappa^{(t)}_{a}-4(\varkappa^{(t)}_{a})^2}{2(1-t)}\\
&= &\frac{2}{\sqrt{2\pi (1-t)}^3}(W_t-\varkappa^{(t)}_{a})\varkappa^{(t)}_{a}\\

\dce
$$
and similarly
$$
\dcb
k_{1-t}(\theta^{(t)}_{a}-W_t,2\theta^{(t)}_{a}-W_t)
&\leq&\frac{2}{\sqrt{2\pi (1-t)}^3}(W_t-\theta^{(t)}_{a})\theta^{(t)}_{a}
\dce
$$
We now compute, for any $0<\epsilon<1$
$$
\dcb
&&\mathbb{E}[\int_\xi^{1-\epsilon}\frac{1}{|W_\xi-W_t|}dt]\\
&=&\mathbb{E}[\int_0^{1-\epsilon}\frac{1}{|W_\xi-W_t|}\ind_{\{\xi<t\}}dt]\\
&=&\mathbb{E}[\int_0^{1-\epsilon}\mathbb{E}[\frac{1}{|W_\xi-W_t|}\ind_{\{\xi<t\}}|\mathcal{F}^\circ_t] dt]\\

&=&\mathbb{E}[\int_0^{1-\epsilon}\left(\int_0^t\frac{1}{|W_a-W_t|}\ind_{\{\varkappa^{(t)}_a\geq 0\}}k_{1-t}(W_t-\varkappa^{(t)}_a,W_t-2\varkappa^{(t)}_a)\ 2 d_a\varkappa^{(t)}_a\right.\\
&&+
\left.\int_0^t\frac{1}{|W_a-W_t|}\ind_{\{\theta^{(t)}_a\leq 0\}} k_{1-t}(\theta^{(t)}_a-W_t,2\theta^{(t)}_a-W_t)\ 2 (-d_a\theta^{(t)}_a)\right) dt]\\

&=&\mathbb{E}[\int_0^{1-\epsilon}\left(\int_0^t\frac{1}{|\varkappa^{(t)}_a-W_t|}\ind_{\{\varkappa^{(t)}_a\geq 0\}}k_{1-t}(W_t-\varkappa^{(t)}_a,W_t-2\varkappa^{(t)}_a)\ 2 d_a\varkappa^{(t)}_a\right.\\
&&+
\left.\int_0^t\frac{1}{|\theta^{(t)}_a-W_t|}\ind_{\{\theta^{(t)}_a\leq 0\}} k_{1-t}(\theta^{(t)}_a-W_t,2\theta^{(t)}_a-W_t)\ 2 (-d_a\theta^{(t)}_a)\right) dt]\\

&\leq&\frac{2}{\sqrt{2\pi (1-t)}^3}\mathbb{E}[\int_0^{1-\epsilon}\left(\int_0^t\ind_{\{\varkappa^{(t)}_a\geq 0\}}|\varkappa^{(t)}_{a}|\ 2 d_a\varkappa^{(t)}_a
+
\int_0^t\ind_{\{\theta^{(t)}_a\leq 0\}} |\theta^{(t)}_{a}|\ 2 (-d_a\theta^{(t)}_a)\right) dt]\\

&\leq&\frac{2}{\sqrt{2\pi (1-t)}^3}\mathbb{E}[\int_0^{1-\epsilon}\left(\  (\varkappa^{(t)}_t)^2\ind_{\{\varkappa^{(t)}_t>0\}}
+
(\theta^{(t)}_t)^2\ind_{\{\theta^{(t)}_t<0\}}\right) dt]\\

&=&\frac{2}{\sqrt{2\pi (1-t)}^3} \int_0^{1-\epsilon} \mathbb{E}[W_t^2]\ dt\\
&<&\infty
\dce
$$
This implies that $$
\int_\xi^{1}\frac{1}{|W_\xi-W_t|}dt\ <\infty
$$
Applying the first of the above three inequalities, we see that the random measure $\mathsf{d}\chi^\cup$ has effectively a distribution function $\chi^\cup$. Consequently, $W-W_\xi$ is a $(\mathbb{P},\mathbb{G})$ special semimartingale with drift $\chi^\cup$. 

\subsection{Ending remark}
To end this example, we write the drift in another form : for $\xi<t\leq 1$,
$$
\dcb
&&\frac{\frac{W_t}{1-t}\exp\{-\frac{W_t^2}{2(1-t)}\}-\frac{(W_t-W_1)}{1-t}\exp\{-\frac{(W_t-W_1)^2}{2(1-t)}\}}{-\exp\{-\frac{W_t^2}{2(1-t)}\}+\exp\{-\frac{(W_t-W_1)^2}{2(1-t)}\}}\\

&=& \frac{\frac{W_t}{1-t}\exp\{-\frac{W_t^2}{2(1-t)}\}}{-\exp\{-\frac{W_t^2}{2(1-t)}\}+\exp\{-\frac{(W_t-W_1)^2}{2(1-t)}\}}
+\frac{-\frac{(W_t-W_1)}{1-t}\exp\{-\frac{(W_t-W_1)^2}{2(1-t)}\}}{-\exp\{-\frac{W_t^2}{2(1-t)}\}+\exp\{-\frac{(W_t-W_1)^2}{2(1-t)}\}}\\
\\
&=&\frac{\frac{W_t}{1-t}\exp\{-\frac{W_t^2}{2(1-t)}\}-\frac{(W_t-W_1)}{1-t}\exp\{-\frac{W_t^2}{2(1-t)}\}}{-\exp\{-\frac{W_t^2}{2(1-t)}\}+\exp\{-\frac{(W_t-W_1)^2}{2(1-t)}\}}
+\frac{\frac{W_t-W_1}{1-t}\exp\{-\frac{W_t^2}{2(1-t)}\}-\frac{(W_t-W_1)}{1-t}\exp\{-\frac{(W_t-W_1)^2}{2(1-t)}\}}{-\exp\{-\frac{W_t^2}{2(1-t)}\}+\exp\{-\frac{(W_t-W_1)^2}{2(1-t)}\}}\\
\\
&=&\frac{\frac{W_1}{1-t}\exp\{-\frac{W_t^2}{2(1-t)}\}}{-\exp\{-\frac{W_t^2}{2(1-t)}\}+\exp\{-\frac{(W_t-W_1)^2}{2(1-t)}\}}-\frac{W_t-W_1}{1-t}\\
\\
&=& 
\frac{\frac{W_1}{1-t}}{-1+\exp\{-\frac{(W_t-W_1)^2}{2(1-t)}+\frac{W_t^2}{2(1-t)}\}}-\frac{W_t-W_1}{1-t}\\
\\
&=&\frac{1}{-1+\exp\{\frac{2W_tW_1-W_1^2}{2(1-t)}\}}\frac{W_1}{1-t}-\frac{W_t-W_1}{1-t}
\dce
$$
Be ware of the component $-\frac{W_t-W_1}{1-t}\ind_{\{\xi<t\leq 1\}}dt$ in the above formula. This term is intrinsically related with the brownian motion bridge as it is indicated in \cite{AM}. In fact, using brownian motion bridge, we can have a different and shorter proof of the above enlargement of filtration formula.

\

\section{Expansion with the future infimum process}\label{futur-infermum}

This is an example which does not fit the classification of initial enlargement or progressive enlargement. The example is a descendant of Pitman's theorem \cite{Pitman}. See \cite{Jeulin79, Jeulin80, MY, Rauscher, saisho, takaoka, Yor97}. The use of the theory of enlargement of filtration to study Pitman's theorem was initiated in \cite{Jeulin79}, and then, in \cite{Jeulin80, Yor97}. Later in \cite{Rauscher}, it was explained that the conditional expectation in the enlarged filtration could be computed directly. However, the techniques used in these papers stayed outside of general theory. In this section we will show the local solution method working on this example.

\subsection{The setting}

Consider a regular conservative diffusion process $Z$ taking values in $]0,\infty[$ defined on a polish probability space $\Omega$. Assume
\begin{enumerate}
\item
$\lim_{t\uparrow \infty} Z_{t} = \infty,$
\item
a scale function $e$ can be chosen such that $e(\infty)=0,e(0+)=-\infty$.
\end{enumerate}
Define the process $I = (I_{t})_{t>0}$ where $I_{t}=\inf_{s\geq
t} Z_{t},\forall t>0$. Recall that (cf.
\cite[Proposition(6.29)]{Jeulin80}), for any finite $\mathbb{F}$ stopping time $T$,
\[ P[I_{T}\leq a|{\mathcal{F}}_{T}] = e(Z_{T}) \int_{0}^{Z_{T}}\ d(e(c)^{-1})\ 1_{(0,a]}(c)
\]
We consider the process $I$ as a map from $\Omega$ into the space $E=C_0(\mathbb{R}_+,\mathbb{R})$ of all continuous functions on $\mathbb{R}_+$. Let $\mathbb{I}$ be the natural filtration on $E$ generated by the coordinates functions. Let $\mathbb{F}^\circ$ be the natural filtration of $Z$. Define the filtration $\mathbb{G}^\circ$ as in Subsection \ref{h-formula}, as well as the filtrations $\mathbb{F}$ and $\mathbb{G}$.   

\bl\label{lemma:Z=I}
For $0<s<t$, let  $U(s,t)=\inf_{s\leq v\leq t}Z_{v}$. Then, $\mathbb{P}[I_{t} = U(s,t)]=0$. Consequently, \[
\mathbb{P}[\mbox{\textup{there are at least two distinct points $v,v'\geq s$ such that $Z_{v}=I_{s},Z_{v'}=I_{s}$}}]=0.\]
\el 

\textbf{Proof.} To prove the first assertion, we write
\[\begin{array}{lcl}
\mathbb{P}[I_{t} = U(s,t)]&=& P[P[I_{t} = a|{\mathcal{F}}_{t}]_{a=U(s,t)}]
=\mathbb{E}[e(Z_{t}) \int_{0}^{Z_{t}} d(e(c)^{-1}) 1_{\{c=U(s,t)\}}]=0.
\end{array}\]
The second assertion is true because it is overestimated by $\mathbb{P}[\cup_{t\in \mathbb{Q}_{+},t>s}\{U(s,t)=I_{t}\}].$ \ok

As usual, for $x\in E$ and $s\geq 0$, the coordinate of $x$ at $s$ is denoted by $x_s$. We define $D_{u}(x) = \inf\{s>u: x_s 
\neq x_u\}$. For $0\leq u<t$, $I_u=U(u,t)\wedge I_t$. Therefore, $$
\{t\leq D_{u}(I)\}=\{I_u=I_t\}=\{U(u,t)\geq I_{t}\}
$$
and on the set $\{t\leq D_{u}(I)\}$, the stopped process $I^{t}$ coincides with the process $U^{u}\wedge I_{t}=U^{u}\wedge I_{u}=I^u$, where $U^{u}$ denotes the process $(U(s\wedge u,u))_{s\geq 0}$.

\subsection{Computation of $\alpha^u$}

Let $0\leq u<t$. Let $f$ be a bounded $\mathcal{I}_t$ measurable function. Consider the random measure $\pi_{u,I/I}$.
\[\begin{array}{lcl}
&&\int \pi_{u,I/I}(I,dx) f(x) \ind_{\{t\leq D_{u}(x)\}} 
=\mathbb{E}[f(I)\ind_{\{t\leq D_{u}(I)\}}|\sigma(I_s: 0\leq s\leq u)]\\
&=&\mathbb{E}[f(I^t)\ind_{\{t\leq D_{u}(I)\}}|\sigma(I_s: 0\leq s\leq u)]\ \mbox{ because $f\in\mathcal{I}_t$}\\
&=&\mathbb{E}[f(I^u)\ind_{\{t\leq D_{u}(I)\}}|\sigma(I_s: 0\leq s\leq u)]\\
&=&f(I^u)\mathbb{E}[\ind_{\{t\leq D_{u}(I)\}}|\sigma(I_s: 0\leq s\leq u)]\\
&=& f(I^{u}) q(u,t,I^{u}),
\end{array}\]
where $q(u,t,x)$ is a $\mathcal{I}_{u}$-measurable positive function such that $q(u,t,I)$ is a version of
$\mathbb{P}[t\leq D_{u}(I)|\sigma(I_s: 0\leq s\leq u)]$. The function $q(u,t,x)$ can be chosen l\`adc\`ag decreasing in $t$. 
Set $\eta(u,t,x) = \ind_{\{q(u,t,x)>0\}}, x\in E$. In the following computations, the parameters $u,t$ being fixed, we omit them from writing. We have 
$$
\int\pi_{u,I/I}(I,dx)(1-\eta(x)) \ind_{\{t\leq D_{u}(x)\}}
=\mathbb{E}[\ind_{\{t\leq D_{u}(I)\}}|\sigma(I_s: 0\leq s\leq u)](1-\eta(I)) 
= 0.
$$
Turn back to the random measure $\pi_{u,I/F}$ and $\pi_{u,I/I}$. 
\[\begin{array}{lcl}
&&\int \pi_{u,I/F}(\mathsf{i},dx'') \int \pi_{u,I/I}(x'',dx) f(x)\eta(x) \ind_{\{t\leq D_{u}(x)\}}
\\
&=&\int \pi_{u,I/F}(\mathsf{i},dx'')  f(x''^{u})\eta(x''^u) q(x''^{u})\\
&=&\mathbb{E}[f(I^{u})\eta(I^u) q(I^{u})|\mathcal{F}^\circ_u]\\
&=&\mathbb{E}[f(U^{u}\wedge I_u)\eta(U^{u}\wedge I_u) q(U^{u}\wedge I_u)|\mathcal{F}^\circ_u]\\
&=&\mathbb{E}[e(Z_{u}) \int_{0}^{Z_{u}}\ d(e(c)^{-1})\ f(U^{u}\wedge c)\eta(U^{u}\wedge c) q(U^{u}\wedge c)|\mathcal{F}^\circ_u]\\
&=&e(Z_{u}) \int_{0}^{Z_{u}}\ d(e(c)^{-1})\ f(U^{u}\wedge c)\eta(U^{u}\wedge c) q(U^{u}\wedge c)
\end{array}\]
Consider nextly the random measure $\pi_{t,I/F}$.
$$
\dcb
&&\int \pi_{t,I/F}(\mathsf{i},dx) f(x)\eta(x) \ind_{\{t\leq D_{u}(x)\}}\\
&=&\mathbb{E}[f(I)\eta(I) \ind_{\{t\leq D_{u}(I)\}}|\mathcal{F}^\circ_t]\\
&=&\mathbb{E}[f(U^u\wedge I_t)\eta(U^u\wedge I_t) \ind_{\{U(u,t)\geq I_t\}}|\mathcal{F}^\circ_t]\\
&=&\mathbb{E}[e(Z_{t})  \ \int_{0}^{Z_{t}}  d(e(c)^{-1}) \ f(U^{u}\wedge c)\eta(U^{u}\wedge c)\ind_{\{U(u,t)\geq c\}}|\mathcal{F}^\circ_t]\\
&=&e(Z_{t})  \ \int_{0}^{Z_{t}}  d(e(c)^{-1}) \ f(U^{u}\wedge c)\eta(U^{u}\wedge c)\ind_{\{U(u,t)\geq c\}}\\
&&\mbox{ noting that $c\leq U(u,t)\leq Z_u\wedge Z_t$}\\
&&\mbox{ and the coordinate of $(U^u\wedge c)$ at $u$ is $Z_u\wedge c$}\\
&=&\frac{e(Z_{t})}{e(Z_{u})}\ e(Z_{u}) \int_{0}^{Z_{u}}  d(e(c)^{-1}) \ f(U^{u}\wedge c)\frac{1}{q(U^{u}\wedge c)}\ind_{\{U(u,t)\geq (U^u\wedge c)_u\}}\eta(U^{u}\wedge c)q(U^{u}\wedge c)\\
&=&\frac{e(Z_{t})}{e(Z_{u})}\ \int \pi_{u,I/F}(\mathsf{i},dx'') \int \pi_{u,I/I}(x'',dx)f(x)
\eta(x)\ind_{\{t\leq D_{u}(x)\}}\frac{1}{q(x)}
\ind_{\{U(u,t)\geq x_u\}}
\dce
$$ 
For $x\in E$ and $\omega\in\Omega$, let \[\begin{array}{lcl}
W_{t}(x) &=& \frac{1}{q(u,t+,x)},\\
\psi(\omega,x) &=& \inf\{s\geq u: Z_{s}(\omega) < x_u\},\\
\lambda(x) &=& \inf\{s\geq u:   q(u,s,x) = 0\},\\
\rho(\omega,x) &=&
\psi(\omega,x)\wedge \lambda(x) \wedge D_{u}(x).
\end{array}\] 
According to Lemma \ref{alpha},
$\tau^{u}(\omega,x)\geq \lambda(x)\wedge D_{u}(x)$ and
\[
\alpha_{t}^{u}\ind_{\{\rho>t\}} = \frac{e(Z_{t})}{e(Z_{u})} W_{t}\ind_{\{\rho>t\}}, \ u\leq t.
\]

\subsection{Local and global solution}

The integration par parts formula gives the following identity under the probability $\nu^{u}$ on the interval $(u,\infty)$:
$$\begin{array}{lcl}
\alpha^{u}_{t}\ind_{\{t<\rho\}}& =& \frac{1}{e(Z_{u})}e(Z_{t\wedge \rho})W_{t}^{\rho-}
-\frac{1}{e(Z_{u})}e(Z_{\rho}) W_{\rho-}\ind_{\{\rho\leq t\}} \\
&=&\frac{1}{e(Z_{u})}e(Z_{u})W_{u}^{\rho-} + \frac{1}{e(Z_{u})}\int_{u}^{t\wedge \rho}W_{s-}de(Z)_s + \frac{1}{e(Z_{u})}\int_{u}^{t\wedge \rho}e(Z)_s dW^{\rho-}_s 
-\frac{1}{e(Z_{u})}e(Z_{\rho}) W_{\rho-}\ind_{\{\rho\leq t\}}
\end{array}
$$
Here $e(Z)$ is viewed as a $(\nu^{u},{\mathbb{J}}_{+})$ local martingale.

Let $\langle e(Z) \rangle$ be the quadratic variation of $e(Z)$ under $\mathbb{P}$ in the filtration $\mathbb{F}$. We notice
that  it is $\mathbb{F}$-predictable. It is straightforward to check that $\langle e(Z) \rangle$ is also a version of the
predictable quadratic variation of $e(Z)$ with respect to
$(\nu^{u},{\mathbb{J}}_{+})$. It results that
$$
\cro{e(Z),\alpha^{u}\ind_{[u,\rho)}}^{\nu^{u}\cdot{\mathbb{J}}_{+}}
=
e(Z_{u})^{-1}W_{-}\ind_{(u,\rho]}\centerdot \langle e(Z)
\rangle
$$ on the interval $[u,\infty)$. Pull these quantities back the the original space, we obtain$$
\dcb
\psi(\phi)&=&\inf\{s>u: Z_{s}<I_{s}\} = \infty\\
\lambda(I)&=&\inf\{s>u:   q(u,s,I) = 0\}\geq D_{u}(I)\ \mbox{ $\mathbb{P}$-a.s.}\\
\rho(\phi)&=&D_u(I)
\dce
$$
where the second inequality comes from the fact $
\mathbb{P}[q(u,t,I)=0,t<D_u(I)]=0.
$
Applying Theorem \ref{theorem:our-formula}, we conclude that $\ind_{(u,D_u(I)]}\centerdot e(Z)$ is a $(\mathbb{P},\mathbb{G})$ special semimartingale with drift given by:
\[
\mathsf{d}\chi^{(u,D_u(I)]} = \ind_{(u,D_u(I)]} \frac{1}{e(Z)} d\langle e(Z) \rangle
\]

Consider $\mathtt{B}_u=(u,D_u(I)]$ for $u\in\mathbb{Q}_+$. Note that $I$ increases only when $Z=I$. Therefore, according to Lemma \ref{lemma:Z=I}, $D_{u}(I) = \inf\{s>u; Z_{s}=I_{s}\}$. Hence, the interior of $\cup_{u\in\mathbb{Q}_+}(u, D_u(I)]$ is $\{Z>I\}$. The family of $\mathtt{B}_u$ satisfies Assumption \ref{partialcovering}. These random measures $\mathsf{d}\chi^{(u,D_u(I)]}$ yields a diffuse random measure $$
\mathsf{d}\chi^\cup=
\ind_{\{Z>I\}} \frac{1}{e(Z)} \mathsf{d}\langle e(Z) \rangle
$$ 
on $\cup_{u\in\mathbb{Q}_+}(u, D_u(I)]$ which has clearly a distribution function. 

\bl\label{g-decomposition}
Let $\mathtt{G}=\{Z>I\}$. $e(Z)$ is a $\mathcal{G}$-semimartingale. We have the semimartingale decomposition : \[
e(Z)=e(Z_{0})-e(I_{0})+ M + (\ind_{\mathtt{G}}\ \frac{1}{e(Z)})\cdot \langle e(Z) \rangle+V^++e(I),
\]
where $M$ is a continuous 
local $\mathbb{G}$-martingale whose quadratic variation is given by $\ind_{\mathtt{G}}\cdot \langle e(Z) \rangle$, and $V^+$ is the increasing process defined in Lemma \ref{Z-Ytheorem}. 
\el

\textbf{Proof.} The difference between $\mathtt{G}$ and $\cup_{u\in\mathbb{Q}_+}(u, D_u(I)]$ is contained in $\cup_{u\in\mathbb{Q}_+}[D_u(I)]$ which is a countable set. Let $\mathtt{A}$ be the set used in Lemma \ref{Z-Ytheorem}. We know also that the difference between $\mathtt{A}$ and $\cup_{u\in\mathbb{Q}_+}(u, D_u(I)]$ is a countable set. 

Now we can check that the process $e(Z)$ satisfies Assumption \ref{pieces} with respect to the family $(\mathtt{B}_u, u\in\mathbb{Q}_+)$ and with respect to the $(\mathbb{P},\mathbb{G})$ special semimartingale $e(I)$. Lemma \ref{Z-Ytheorem} is applicable. We make only a remark that $d\cro{e(Z)}$ does not charge $\mathtt{A}\setminus \mathtt{G}$ which is a countable set while $\cro{e(Z)}$ is continuous.
\ok

\subsection{Computation of $V^+$}

Lemma \ref{g-decomposition} solves the enlargement of filtration problem. Now we want to compute the process $V^+$. 

\bl\label{lemma:non-charge}
$\int 1_{\mathtt{G}^{c}} d\langle e(Z) \rangle = 0$. Consequently, if $X=e(Z)-e(I)$ and $l^{0}(X)$ denotes the local time of $X$ at 0, we have$$
V^+=\frac{1}{2}l^{0}(X)
$$
\el

\textbf{Proof.} By Lemme \ref{Z-Ytheorem}, $X=X_0+\ind_\mathtt{G}\centerdot X+V^+$. We can calculate the quadratic variation $[X]$ in $\mathcal{G}$ in two ways : 
\[
\dcb
d[X] = d[\ind_\mathtt{G}\centerdot X + V^+] = d[\ind_\mathtt{G}\centerdot X] = \ind_\mathtt{G}d\langle e(Z) \rangle\\
d[X] = d[e(Z)-e(I)] = d\langle e(Z) \rangle.
\dce
\]
This yields $\int 1_{\mathtt{G}^{c}} d\langle e(Z) \rangle = 0$. Using the formula in Lemma \ref{g-decomposition}, applying \cite[Chapter VI, Theorem 1.7]{Revuz-Yor}, $$
\frac{1}{2}l^0_t
=\ind_{\{X=0\}}\centerdot\left((\ind_{\mathtt{G}}\ \frac{1}{e(Z)})\cdot \langle e(Z) \rangle+V^+\right)
=\ind_{\{X=0\}}\centerdot V^+ = V^+
$$
\ok

\bl\label{lemma:tribu}
The processes of the form $Hf(I_{t})1_{(t,u]}$, where $0<t<u<\infty$, $H \in {\mathcal{F}}_{t}$, $f$ is a bounded continuous function on $]0,\infty[$, generate the $\mathbb{G}$-predictable $\sigma$-algebra.
\el

\textbf{Proof.} It is because $I_s=U(s,t)\wedge I_t$ for any $0\leq s\leq t$ and $U(s,t)\in\mathcal{F}_t$. \ok

For any $a>0,b>0$, set $\delta_b(a) = \inf\{s>b: Z_{s}\leq a\}$. Then, on the set $\{a< Z_{b},\delta_b(a)<\infty\}$, $I_{s}=I_{b}\leq a$ for $b\leq s\leq \delta_b(a)$,
and consequently, $e(I_{u\wedge\delta_b(a)})- e(I_{t\wedge\delta_b(a)}) = 0$ for $b\leq t\leq u$. On the other hand, on the set $\{a< Z_{b}, \delta_b(a) = \infty\}=\{a< I_{b}\}$. In sum, 
\[
\ind_{\{a< Z_{b}\}}(e(I_{u\wedge\delta_b(a)})- e(I_{t\wedge\delta_b(a)}))
=\ind_{\{a< I_{b}\}}(e(I_{u})- e(I_{t})), \ b\leq t\leq u.
\]
Note that $l^0$ is a continuous process increasing only when $Z=I$. The above argument holds also for $l^0$:
\[
\ind_{\{a< Z_{b}\}}(l^0_{u\wedge\delta_b(a)}- l^0_{t\wedge\delta_b(a)})
=\ind_{\{a< I_{b}\}}(l^0_{u}- l^0_{t}), \ b\leq t\leq u.
\]

\label{lemma:projection} 

\bl \label{corollary:projection}
For any $a>0$, for any finite $\mathbb{F}$ stopping time $\beta$, the $\mathbb{F}$-predictable dual projection of the increasing process $\ind_{\{a< I_{\beta}\}}\ind_{(\beta,\infty)}\centerdot e(I)$
is given by
$$
-\ind_{(\beta,\infty)}\ind_{\{U(\beta,\cdot)> a\}}\frac{1}{2e(Z)} \centerdot \langle e(Z) \rangle
$$
\el

\textbf{Proof.} Let $a>0$ to be a real number. Let $V$ be a bounded $\mathbb{F}$ stopping time such that, stopped at $V$, the two processes $(\ln(-e(a))-\ln(-e(Z_{s}))-1) \centerdot e(Z_{s})$ and $\frac{1}{2e(Z_{s})} \centerdot \langle e(Z) \rangle_{s}$ are uniformly integrable. Let $R,T$ be two others $\mathbb{F}$ stopping times such that $R\leq T\leq V$. Notice that, on the set $\{a\leq I_{R}\}$, $e(a)\leq e(I_T)\leq e(I_V)< 0$. Under these condition, we have the following computations:
\[\begin{array}{lcl}
&&\mathbb{E}[e(I_{V})\ind_{\{a< I_{R}\}}|{\mathcal{F}}_{T}]\\
&=& \mathbb{E}[\ind_{\{a< U(R,V)\}} \mathbb{E}[e(I_{V})\ind_{\{a< I_{V}\}}|{\mathcal{F}}_{V}] |{\mathcal{F}}_{T}]\\
&=& \ind_{\{a< Z_{R}\}}\mathbb{E}[\ind_{\{V< \delta_R(a)\}}e(Z_{V}) \int_{a}^{ Z_{V}}  e(c) \ d(e(c)^{-1})
|{\mathcal{F}}_{T}]\\ 
&=& \ind_{\{a< Z_{R}\}}\mathbb{E}[\ind_{\{V< \delta_R(a)\}} e(Z_{V})(\ln(-e(a))-\ln(-e(Z_{V})))
|{\mathcal{F}}_{T}]\\ 
&=& \ind_{\{a< Z_{R}\}}\mathbb{E}[e(Z_{V\wedge \delta_R(a)})(\ln(-e(a))-\ln(-e(Z_{V\wedge \delta_R(a)})))
|{\mathcal{F}}_{T}].
\end{array}\]
We have the same computation when $V$ is replaced by $T$. By Ito's formula (in the filtration $\mathbb{F}$), 
\[\begin{array}{lcl}
&&d\left( e(Z_{s})(\ln(-e(a))-\ln(-e(Z_{s})))\right)\\
&=&(\ln(-e(a))-\ln(-e(Z_{s}))-1) de(Z)_{s} - \frac{1}{2e(Z_{s})} d\langle e(Z) \rangle_{s}, \ s\geq 0.
\end{array}\]
Taking the difference of the above computations, we obtain
$$
\dcb
&&\mathbb{E}[\ind_{\{a< I_{R}\}}\int_T^Vde(I)_s|{\mathcal{F}}_{T}]
=\mathbb{E}[(e(I_{V})-e(I_T))\ind_{\{a< I_{R}\}}|{\mathcal{F}}_{T}]\\
&=& \ind_{\{a< Z_{R}\}}\mathbb{E}[\int_{T\wedge \delta_R(a)}^{V\wedge \delta_R(a)}(\ln(-e(a))-\ln(-e(Z_{s}))-1) de(Z)_{s} - \int_{T\wedge \delta_R(a)}^{V\wedge \delta_R(a)}\frac{1}{2e(Z_{s})} d\langle e(Z) \rangle_{s}
|{\mathcal{F}}_{T}]\\
&=& -\ind_{\{a< Z_{R}\}}\mathbb{E}[\int_{T}^{V}\ind_{\{s<  \delta_R(a)\}}\frac{1}{2e(Z_{s})} d\langle e(Z) \rangle_{s}
|{\mathcal{F}}_{T}]\\
&=& -\mathbb{E}[\int_{T}^{V}\ind_{\{U(R,s)> a\}}\frac{1}{2e(Z_{s})} d\langle e(Z) \rangle_{s}|{\mathcal{F}}_{T}]
\dce
$$
This identity yields that, for any positive bounded elementary $\mathbb{F}$ predictable process $H$,
$$
\dcb
\mathbb{E}[\ind_{\{a< I_{R}\}}\int_R^VH_sde(I)_s]

&=& -\mathbb{E}[\int_{R}^{V}H_s\ind_{\{U(R,s)> a\}}\frac{1}{2e(Z_{s})} d\langle e(Z) \rangle_{s}]
\dce
$$
It is to notice that we can find an increasing sequence $(V_n)_{n\geq 1}$ of $\mathbb{F}$ stopping times tending to the infinity, and each of the $V_n$ satisfies that assumption on $V$. Replace $V$ by $V_n$ and let $n$ tend to the infinity, monotone convergence theorem then implies
$$
\dcb
\mathbb{E}[\ind_{\{a< I_{R}\}}\int_R^\infty H_sde(I)_s]

&=& -\mathbb{E}[\int_{R}^{\infty}H_s\ind_{\{U(R,s)> a\}}\frac{1}{2e(Z_{s})} d\langle e(Z) \rangle_{s}]
\dce
$$
Notice that the left hand side is finite. Therefore, the right hand side is concerned by an integrable random variable. Now, replace $R$ by $\beta\wedge V_n$ and let $n$ tend to the infinity. The dominated convergence theorem implies 
$$
\dcb
\mathbb{E}[\ind_{\{a< I_{\beta}\}}\int_\beta^\infty H_sde(I)_s]
&=& -\mathbb{E}[\int_{\beta}^{\infty}H_s\ind_{\{U(\beta,s)> a\}}\frac{1}{2e(Z_{s})} d\langle e(Z) \rangle_{s}]
\dce
$$
This identity proves the lemma. \ok

\bethe
$l^{0}(X) = 2e(I)-2e(I_0)$.
\ethe

\textbf{Proof.} Using the occupation formula (cf. \cite{Revuz-Yor}), for
$0<t<u$, 
\[ l^{0}_{u}(X) - l^{0}_{t}(X)=\lim_{\epsilon \downarrow 0} \frac{1}{\epsilon} \int_{0}^{\epsilon} (l^a_u-l^a_t) da 
=\lim_{\epsilon \downarrow 0} \frac{1}{\epsilon} \int_{t}^{u} 1_{\{X_{v}\leq
\epsilon\}}d\langle X
\rangle_{v}
=
\lim_{\epsilon \downarrow 0} \frac{1}{\epsilon} \int_{t}^{u} 1_{\{e(Z_{v})-e(I_{v})\leq \epsilon\}}d\langle e(Z)
\rangle_{v}  ,
\]
almost surely. Let $a>0$. Let $0\leq R\leq T\leq V$ be three $\mathbb{F}$ stopping times such that, stopped at $V$, the processes $Z$ and $\cro{e(Z)}$ are bounded. Note that this condition implies that the family of random variables $\{l^a_{s\wedge V}: s\geq 0, 0\leq a\leq 1\}$ is uniformly integrable. By convex combination (Vallée-Poussin's theorem), the family $\{\frac{1}{\epsilon}\int_0^\epsilon (l^a_{V}-l^a_{T})da : 0< \epsilon\leq 1\}$ also is uniformly integrable. Consequently, we can write the identity:$$
E[Hf(e(I_{T})) \ind_{\{a< Z_{R}\}}(l^0_{V\wedge \delta_R(a)}-l^0_{T\wedge \delta_R(a)})]
=
\lim_{\epsilon\downarrow 0}E[Hf(e(I_{T})) \ind_{\{a< Z_{R}\}}\frac{1}{\epsilon} \int_{T\wedge \delta_R(a)}^{V\wedge \delta_R(a)} 1_{\{e(Z_{v})-e(I_{v})\leq
\epsilon\}}d\langle e(Z)
\rangle_{v} ]
$$
for any positive continuously differentiable function $f$ with compact support and any positive bounded random variable $H$ in ${\mathcal{F}}_{t}$. We now compute the limit at right hand side. 
\[\begin{array}{ll}
&E[Hf(e(I_{T})) \ind_{\{a< Z_{R}\}}\frac{1}{\epsilon} \int_{T\wedge \delta_R(a)}^{V\wedge \delta_R(a)} 1_{\{e(Z_{v})-e(I_{v})\leq
\epsilon\}}d\langle e(Z)
\rangle_{v} ]\\
=&E[H \ind_{\{a< Z_{R}\}} \frac{1}{\epsilon} \int_{T\wedge \delta_R(a)}^{V\wedge \delta_R(a)} d\langle e(Z)
\rangle_{v}\ f(e(U(T,v))\wedge
e(I_{v}))1_{\{e(Z_{v})-\epsilon\leq
e(I_{v})\}} ]\\
=&E[H \ind_{\{a< Z_{R}\}} \frac{1}{\epsilon} \int_{T\wedge \delta_R(a)}^{V\wedge \delta_R(a)} d\langle e(Z)
\rangle_{v} \ e(Z_{v}) \int_{e^{-1}(e(Z_{v})-\epsilon)}^{Z_{v}}  d(e(c)^{-1})f(e(U(T,v))\wedge
e(c)) ]\\
=&E[H \ind_{\{a< Z_{R}\}} \frac{1}{\epsilon} \int_{T}^{V} \ind_{\{v<\delta_R(a)\}}d\langle e(Z)
\rangle_{v} \ e(Z_{v}) \int_{e^{-1}(e(Z_{v})-\epsilon)}^{Z_{v}}  d(e(c)^{-1})f(e(U(T,v))\wedge
e(c)) ]\\ 
=&-2E[H\ind_{\{a< I_{R}\}} \frac{1}{\epsilon} \int_{T}^{V} de(I_{v}) e(Z_{v})^{2}
 \int_{e^{-1}(e(Z_{v})-\epsilon)}^{Z_{v}} d(e(c)^{-1})\  f(e(U(T,v))\wedge
e(c)) ]\\
&\mbox{ according to Lemma \ref{corollary:projection}}\\
=&-2E[H\ind_{\{a< I_{R}\}} \frac{1}{\epsilon} \int_{T}^{V} de(I_{v}) e(I_{v})^{2}
 \int_{e^{-1}(e(I_{v})-\epsilon)}^{I_{v}} d(e(c)^{-1})\  f(e(U(T,v))\wedge
e(c)) ]\\
\end{array}\]
The above last term is divided into two parts
$$
\dcb
-2E[H\ind_{\{a< I_{R}\}} \frac{1}{\epsilon} \int_{T}^{V} de(I_{v}) e(I_{v})^{2}
 \int_{e^{-1}(e(I_{v})-\epsilon)}^{I_{v}} d(e(c)^{-1})\  f(e(U(T,v))\wedge
e(I_v)) ]\\
\\
-2E[H\ind_{\{a< I_{R}\}} \frac{1}{\epsilon} \int_{T}^{V} de(I_{v}) e(I_{v})^{2}
 \int_{e^{-1}(e(I_{v})-\epsilon)}^{I_{v}} d(e(c)^{-1})\  (f(e(U(T,v))\wedge
e(c))-f(e(U(T,v))\wedge
e(I_v))) ]
\dce
$$
The first part is computed as follows :
$$
\dcb
&&-2E[H\ind_{\{a< I_{R}\}} \frac{1}{\epsilon} \int_{T}^{V} de(I_{v}) e(I_{v})^{2}
 \int_{e^{-1}(e(I_{v})-\epsilon)}^{I_{v}} d(e(c)^{-1})\  f(e(I_T)) ]\\
&=&2E[H\ind_{\{a< I_{R}\}} f(e(I_T)) \frac{1}{\epsilon} \int_{T}^{V} de(I_{v}) e(I_{v})^{2}
 (-\frac{1}{e(I_v)}+\frac{1}{e(I_v)-\epsilon})]\\
&=&2E[H\ind_{\{a< I_{R}\}} f(e(I_T)) \frac{1}{\epsilon} \int_{T}^{V} de(I_{v}) e(I_{v})^{2}
 \frac{\epsilon}{e(I_v)(e(I_v)-\epsilon)}]\\
&=&2E[H\ind_{\{a< I_{R}\}} f(e(I_T)) \int_{T}^{V} de(I_{v}) e(I_{v})^{2}\frac{1}{e(I_v)(e(I_v)-\epsilon)}]\\
\dce
$$
When $\epsilon$ decreases down to zero, this quantity increases to $$
2E[H\ind_{\{a\leq I_{R}\}} f(e(I_T))(e(I_{V})-e(I_T))]
$$
The second part is overestimated by
$$
\dcb
&&\left|-2E[H\ind_{\{a< I_{R}\}} \frac{1}{\epsilon} \int_{T}^{V} de(I_{v}) e(I_{v})^{2}
 \int_{e^{-1}(e(I_{v})-\epsilon)}^{I_{v}} d(e(c)^{-1})\  (f(e(U(T,v))\wedge
e(c))-f(e(U(T,v))\wedge
e(I_v))) ]\right|\\

&\leq&-2E[H\ind_{\{a< I_{R}\}} \frac{1}{\epsilon} \int_{T}^{V} de(I_{v}) e(I_{v})^{2}
 \int_{e^{-1}(e(I_{v})-\epsilon)}^{I_{v}} d(e(c)^{-1})\|f'\|_\infty  |e(c)-e(I_v)| ]\\

&\leq&-2E[H\ind_{\{a< I_{R}\}} \frac{1}{\epsilon} \int_{T}^{V} de(I_{v}) e(I_{v})^{2}
 \int_{e^{-1}(e(I_{v})-\epsilon)}^{I_{v}} d(e(c)^{-1})\|f'\|_\infty  \epsilon ]\\

&\leq&-2E[H\ind_{\{a< I_{R}\}} \int_{T}^{V} de(I_{v}) e(I_{v})^{2}
\int_{e^{-1}(e(I_{v})-\epsilon)}^{I_{v}} d(e(c)^{-1})\|f'\|_\infty]\\

&=&2E[H\ind_{\{a< I_{R}\}} \int_{T}^{V} de(I_{v}) e(I_{v})^{2}
(-\frac{1}{e(I_v)}+\frac{1}{e(I_v)-\epsilon})\|f'\|_\infty]\\

&=&2E[H\ind_{\{a< I_{R}\}} \int_{T}^{V} de(I_{v}) e(I_{v})^{2}
\frac{\epsilon}{e(I_v)(e(I_v)-\epsilon)}\|f'\|_\infty]\\

&\leq&\epsilon 2E[H\ind_{\{a< I_{R}\}} \int_{T}^{V} de(I_{v}) e(I_{v})^{2}
\frac{1}{e(I_v)^2}\|f'\|_\infty]\\

&\leq&\epsilon 2E[H\ind_{\{a< I_{R}\}} (e(I_{V})-e(I_T))\|f'\|_\infty]\\
 
\dce
$$
which tends to zero when $\epsilon$ goes down to zero. In sum, we have
$$
E[Hf(e(I_{T})) \ind_{\{a< Z_{R}\}}(l^0_{V\wedge \delta_R(a)}-l^0_{T\wedge \delta_R(a)})]
=
2E[H\ind_{\{a< I_{R}\}} f(e(I_T))(e(I_{V})-e(I_T))]
$$
Using the formula preceding Lemma \ref{corollary:projection}, we write$$
E[Hf(e(I_{T})) \ind_{\{a< I_{R}\}}(l^0_{V}-l^0_{T})]
=
2E[H\ind_{\{a< I_{R}\}} f(e(I_T))(e(I_{V})-e(I_T))]
$$
or equivalently
$$
E[H\ind_{\{a< I_{R}\}}f(e(I_{T}))\int_T^V dl^0_{s}]
=
2E[H\ind_{\{a\leq I_{R}\}} f(e(I_T))\int_T^V de(I)_s]
$$
Recall that, according to Lemma \ref{lemma:tribu}, the family of all processes of the form $\ind_{\{a\leq I_R\}}Hf(I_T)\ind_{(T,V]}$ generates all $\mathbb{G}$ predictable processes. The above identity means then that, on $\mathcal{P}(\mathbb{G})$, the Dolean-Dade measure of $l^0$ and of $2e(I)$ coincide. That is the theorem. \ok

\bethe
The process $e(Z)$ is a $\mathbb{G}$-semimartingale. Its canonical decomposition is given by :\[
e(Z) = e(Z_0) +M + 2e(I)-2e(I_0) + \frac{1}{e(Z)}\cdot \langle e(Z) \rangle,
\]
where $M$ is a $\mathbb{G}$ local martingale. 
\ethe

Applying Ito's formula, we obtain also an equivalent result

\bcor
The process
$\frac{1}{e(Z)}-\frac{2}{e(I)}$ is a $\mathbb{G}$ local  martingale.
\ecor

\

\

\end{document}